\documentclass{article}

\pdfoutput=1

\usepackage{fohl}

\title{First-order homotopical logic}
\author{Joseph Helfer}
\date{}

\usepackage{fancyhdr}
\begin{document}
\maketitle

\begin{abstract}
\noindent
We introduce a homotopy-theoretic interpretation of intuitionistic first-order logic based on ideas from Homotopy Type Theory. We provide a categorical formulation of this interpretation using the framework of Grothendieck fibrations. We then use this formulation to prove the central property of this interpretation, namely \emph{homotopy invariance}. To do this, we use the result from \cite{sentai} that any Grothendieck fibration of the kind being considered can automatically be upgraded to a 2-dimensional fibration, after which the invariance property is reduced to an abstract theorem concerning pseudonatural transformations of morphisms into 2-dimensional fibrations.
\end{abstract}

\vspace{25pt}
\tableofcontents
\vspace{5pt}

\section{Overview}\label{sec:overview}
The goal of this paper is to introduce a ``homotopy-invariant'' interpretation of first-order logic with equality, to give a description of this interpretation within the framework of categorical logic, and to give an abstract formulation and proof of the homotopy-invariance property within this framework.

The interpretation can be concisely described by the following commutative diagram.
\begin{equation}
  \begin{tikzcd}
    \text{MLTT}
    \ar[r, "\text{Voevodsky-Awodey-Warren-Kapulkin-Lumsdaine}"]&[190pt]\text{Simplicial Sets}\\
    \text{IFOL}\ar[u, "\text{Martin-L\"of 1972}"]
    \ar[ru, "\text{First-order homotopical logic}"', sloped]
  \end{tikzcd}\label{eq:interpretation-diagram}
\end{equation}

On the bottom left we have intuitionistic first-order logic, on the top-left we have Martin-L\"of type theory, and the vertical arrow is the interpretation of IFOL into MLTT which was described in Martin-L\"of's original paper \cite{ml72}. The long horizontal arrow is the homotopy-theoretic interpretation of type theory \cite{awodeywarrenid,kapulkinlumsdaine,warrenthesis} which initiated the subject of Homotopy Type Theory (HoTT). Hence, composing these two interpretations, one obtains a homotopy-theoretic interpretation of first-order logic, which it is the purpose of this paper to elaborate.

In fact, one does not need to go through these two interpretations, as the homotopical semantics for first-order logic can be described directly and very simply -- it is essentially as simple as the ordinary (Tarskian) semantics for first-order logic.
Hence, this interpretation gives us a simpler -- but already interesting -- context in which to consider the ideas involved in HoTT.
A second motivation to study this interpretation is the inherent interest in first-order logic, and in particular in semantics for intuitionistic first-order logic.

Our main result is the following fundamental ``homotopy invariance'' theorem saying, essentially, that homotopy-equivalent structures have the same (first-order) homotopical properties. (See \S\S\ref{subsec:interpretations}-\ref{subsec:htpy-equiv} for the relevant notation and terminology.)

\makebox[\textwidth][c]{
  \begin{minipage}[center]{0.9\linewidth}
    \thm[ \ref{thm:special-invariance}]\label{thm:special-invariance-in-intro}
    Let $M,N\colon\sigma\to\Kan$ be two interpretations of an algebraic signature $\sigma$ in $\Kan$. Let $(\phi,\vec{x})$ be a formula-in-context over $\sigma$, and let $\widehat M_{\vec{x}}(\phi)\in\Kan/M(\vec{x})$ and $\widehat N_{\vec{x}}(\phi)\in\Kan/N(\vec{x})$ be homotopical interpretations of $(\phi,\vec{x})$.

    Then for any homotopy equivalence $\alpha\colon M\to N$, there exists a homotopy equivalence
    \[
      \widehat M_{\vec{x}}(\phi)\to\widehat N_{\vec{x}}(\phi)
    \]
    lying over $\alpha_{\tp\vec{x}}\colon M(\vec{x})\to N(\vec{x})$. In particular, if $\phi$ is a closed formula (i.e., $\vec{x}=\emptyset$), then $M\vDash\phi$ if and only if $N\vDash\phi$.
  \end{minipage}
}

To prove this theorem, we first provide an ``algebraic'' or ``functorial'' version of the homotopical semantics, in the usual style of categorical logic. The above invariance theorem (which we call the ``special'' or ``syntactic'' invariance theorem) is then deduced from the following purely categorical ``abstract'' invariance theorem. (See \S\S\ref{subsec:2-cat-prelims}-\ref{subsec:abstract-invariance-theorem} for the relevant terminology.).

\makebox[\textwidth][c]{
  \begin{minipage}[center]{0.9\linewidth}
    \thm[ \ref{thm:abstract-invariance}]
    Let $\fib{C}$ be a free $h^=$-fibration, let $\fib{C'}$ be a 1-discrete 2-fibration which is also an $h^=$-fibration, and let $(\Phi,\phi),(\Psi,\psi)\colon\fib{C}\to\fib{C'}$ be morphisms of $h^=$-fibrations. Then for any pseudonatural equivalence $\alpha\colon\phi\to\psi$, there exists a pseudonatural equivalence $\Phi\to\Psi$ lying over $\alpha$.
  \end{minipage}
}

There is also a ``second part'' to the abstract invariance theorem, Theorem~\ref{thm:free-htpy-equiv}, which we discuss in the introduction to Part~\ref{sec:invariance}.

We now give some general background on homotopical and categorical semantics, in order to elucidate the significance of the above theorems.

{\bf Acknowledgments:}
I am grateful to McGill's Logic, Category Theory, and Computation seminar and Carnegie Mellon's Homotopy Type Theory seminar for allowing me to speak about this project at a very early stage, and to Steve Awodey for encouraging me to write it up, and also for pointing out an interesting connection between our invariance theorem and the Univalence Axiom (which was also independently observed by Ulrik Buchholtz).

I would also like to acknowledge the various contributors to \href{http://ncatlab.org}{\texttt{ncatlab.org}}, which is an invaluable resource in general, and has been no less for this project.

\subsection{Homotopical structures}\label{subsec:homotopical-structures}
We do not have the space here to explain the somewhat complicated and many-faceted history behind HoTT, and must assume the reader already has some familiarity with it. We will say only this much about it: the central concept of interest is that of \emph{homotopy type} which roughly means ``topological space up to homotopy equivalence'' (though, importantly, there are other equivalent ``models'' for the notion of homotopy type), and the main idea is that it is homotopy types, and not sets, that are the fundamental building blocks in mathematics, with sets arising as special cases (namely the ``discrete'' homotopy types). Moreover, the notion of ``two elements of a set being equal'' is to be replaced by ``two points in a space being joined by a path'' (note, for instance, that in the case of discrete spaces -- i.e., sets -- these two notions coincide).

The fundamental role of sets in mathematics is as \emph{carriers for structures}. We will be in particular concentrating on the \emph{algebraic structures}, such as groups and rings. Hence, if homotopy types are to take the role of sets, there should be a notion of structures on homotopy types.

In fact, putting it this way distorts the history, as the appearance of such ``homotopical structures'' was one of the first developments in the history leading up to HoTT, and was, for example, studied in depth in the text \cite{boardmanvogt} titled ``Homotopy Invariant Algebraic Structures on Topological Spaces''. Here, the emphasis is on \emph{homotopy-invariance}, meaning roughly that if some topological space $X$ carries such a structure, so should any homotopy equivalent space $Y$.

For example, if one considers the structure of a (continuous) associative binary operation on $X$ (i.e., the structure of a \emph{semigroup}), we find that this indeed induces a binary operation on $Y$, but it is no longer associative. Rather, it is \emph{homotopy-associative}, in the sense that the two resulting maps $Y\times Y\times Y\to Y$ are not \emph{equal}, but \emph{homotopic}. In particular, for each $x,y,z\in Y$, we find that the elements $(x\cdot y)\cdot z$ and $ x\cdot(y\cdot z)\in Y$ are not \emph{equal} but rather \emph{joined by a path}; this is the first instance of ``joined by a path'' replacing ``equal''. Such a homotopy-associative binary operation is the first instance of a homotopy-invariant algebraic structure.

Before proceeding, we would like to point out as soon as possible that the main homotopy-invariant algebraic structures of interest in topology are the so-called ``higher structures'' -- for example ``$A_\infty$'' and ``$E_\infty$'' structures. These result from, for example, not just demanding that the two maps $X\times X\times X\to X$ are homotopic, but \emph{specifying} such a homotopy as part of the structure, and then demanding that certain induced maps $X\times X\times X\times X\to X$ are homotopic (and then specifying such a homotopy as part of the structure, and so on). HoTT is able (to some extent) to deal with such higher structures, whereas the First-Order Homotopical Logic considered here cannot, which is a serious limitation.

Returning to our discussion of homotopical structures, we next note that, in the subject of model theory, the way that structures satisfying certain properties are studied is via a suitable \emph{language}. For example, given a set $X$ with a binary operation ``$\cdot$'', demanding that the operation be associative is tantamount to requiring that the structure $(X,\cdot)$ \emph{satisfy} the sentence
\[
  \forall x,y,z\ [(x\cdot y)\cdot z=x\cdot(y\cdot z)],
\]
formulated in first-order logic over the algebraic signature (i.e., specification of a set of sorts and of operations with given arities -- see Definition~\ref{defn:signature}) consisting of a single sort with a single binary operation ``$\cdot$''. Of course, once one has such a definition in place, one can consider sentences of arbitrary complexity, and structures satisfying them.

Our goal here is to establish a similar framework for homotopical structures -- i.e., define what it means for a given homotopical structure to satisfy an arbitrary first-order sentence. Moreover, it should be homotopy invariant in the sense that homotopy equivalent structures satisfy the same sentences. For example, a space with a binary operation should satisfy the above sentence if and only if the operation is homotopy-associative.

As another example (both of these are worked out in the examples section \S\ref{subsec:examples}), the sentence
\[
  \exists x\forall y(x=y),
\]
over the empty signature, which would normally mean ``there is a single element'', should be a satisfied by a space $X$ if and only if $X$ is contractible.

Before proceeding to describe how the homotopical semantics are defined, let us make two general comments.

The first is that, for reasons we will discuss below, it turns out to be difficult and unnatural to define the homotopical semantics using topological spaces. Instead, we use simplicial sets (more precisely, Kan complexes) -- as indicated in (\ref{eq:interpretation-diagram}) -- which are another of the ``models of homotopy types'' mentioned above (it is possible in principle to use other categories as well -- see \S\ref{subsec:h-fib-of-spaces}). The case of topological spaces is then treated via that of Kan complexes (see \S\ref{subsec:top-spaces}). This is somewhat disappointing given the historical importance and intuitive appeal of topological spaces; however, it can be argued that Kan complexes are a more natural setting for homotopy theory in the context of HoTT, specifically in light of the so-called ``Homotopy Hypothesis'' of Grothendieck (formulated in \cite{pursuingstacks} -- the name seems to be due to John Baez \cite{baez-homohypo}).

The second is the role of intuitionistic logic, about which we have not said anything yet. An important (really, the central) feature of any interpretation of first-order logic is \emph{soundness}: if a structure satisfies a sentence $\phi$, it should also satisfy any sentence $\psi$ which is a logical consequence of $\phi$. Accordingly, our homotopical semantics will satisfy this property, but only -- for reasons we discuss below -- for intuitionistic logic, and not all of classical logic.

\subsection{Propositions-as-spaces}\label{subsec:propositions-as-spaces}
The key to the homotopical interpretation of first-order logic (as well as the interpretation of IFOL into MLTT -- i.e., the vertical arrow in (\ref{eq:interpretation-diagram})) -- is the idea, inspired by the ``BHK (Brouwer-Heying-Kolomogorov)'' interpretation of intuitionistic first-order logic (see \cite[\S3.1]{troelstra-van-dalen-cim}), that formulas should not be interpreted just as a \emph{truth-value}, but rather a \emph{set} (or, here, \emph{space}), to be thought of as the ``set of its proofs'', so that the proposition is interpreted as true just if this set is inhabited. The various logical connectives are then each interpreted as a certain operation on sets/spaces; the homotopy-invariance of the interpretation could then be ensured by having each of these operations on spaces be homotopy invariant.

The formalization of this -- with one complication, which we return to below -- can be obtained as a special case of the general ``fibrational'' (or ``hyperdoctrinal'') interpretation of first-order logic (which we assume familiarity with -- see \cite{makkai-lauchli1,jacobscatlogic}).

Fix an algebraic signature $\sigma$, a category $\B$ with finite products, and an interpretation $M\colon\sigma\to B$ of $\sigma$ in $\B$ (in the usual sense of categorical logic -- see Definition~\ref{defn:signature}). Then, given any \emph{$h^=$-fibration}\footnote{This (as well as ``$h$-fibration'' for the version without equality) is the name used in \cite{makkai-lauchli1}; the ``h'' stands for ``Heyting'' and is unrelated to the notion of ``Hurewicz fibration'' as in \cite[p.~340]{may-ponto}. The corresponding hyperdoctrinal notion is simply called a ``hyperdoctrine'' in \cite{seely-hyperdocrtines}, and the version with preorder fibers is called a ``first-order fibration'' in \cite{jacobscatlogic}. The original notion of hyperdoctrine (which is adapted to higher-order rather than first-order logic) is due to Lawvere \cite{lawvereadj,lawvereequality}.} (see Definition~\ref{defn:h-fibration}) $\fibr{C}{C}{B}$, we can consider the interpretation of first-order formulas over $\sigma$ in $\fib{C}$ (see Definition~\ref{defn:interpretation-over}).

Thus, different choices of $h^=$-fibrations give us different semantics for intuitionistic first-order logic, and it is now just a question of finding one which will give us the homotopical semantics.

For instance, taking an appropriate fibration $\fib{P(\Set)}$ in which the fiber over a set $X$ is the power set $\mathcal{P}(X)$, we obtain the classical, set-theoretic semantics. Next, if we take the codomain fibration $\fibr{F(\Set)}{\Set^\to}{\Set}$ whose fiber over $X$ is the slice category $\Set/X$ (i.e., the category of $X$-indexed families of sets), then we obtain the ``propositions-as-sets'' interpretation alluded to above. We note that these examples are in fact not just $h$-fibrations but \emph{Boolean} fibrations, i.e., the corresponding interpretations are sound for classical, and not just intuitionistic, logic.\footnote{As remarked in \cite[p.34,~Exercise~1.3.4]{troelstra-van-dalen-cim}, this shows that this formalization of the ``propositions-as-sets'' interpretation is not really faithful to the BHK-interpretation, which is meant to be an interpretation only of \emph{intuitionistic} logic. Note, however, that Läuchli \cite{lauchlisemantics} uses a variant of this interpretation to give a complete semantics for intuitionistic first-order logic. A fibrational interpretation of Läuchli's result is worked out in \cite{makkai-lauchli1}.}

The last example immediately generalizes: for any category $\B$ which is \emph{locally bicartesian closed} (i.e., locally cartesian closed with finite coproducts), the fibration $\fib{F(\B)}$ will be an $h^=$-fibration, resulting in a ``propositions-as-objects-in-$\B$'' interpretation (\textit{cf}. \cite{palmgren-bhk}). For general $\B$, this will be a genuinely intuitionistic interpretation, i.e., $\fib{F(\B)}$ will not be a Boolean fibration.

Hence, it might seem that our ``propositions-as-spaces'' semantics should simply result by taking $\B=\sSet$, the category of simplicial sets, which is locally bicartesian closed (as opposed to the category of topological spaces -- this is one reason we cannot use the latter for the homotopical semantics). This ``almost'' works, except that \emph{equality} is interpreted the wrong way.

To explain this, let us recall, in concrete terms, how the ``propositions-as-sets'', and more generally, the ``propositions-as-objects-of-$\B$'' interpretation is defined. We have the following table, with logical connectives in the first row, and operations on (families of) sets (or more generally, objects in a locally bicartesian closed category) in the second (note that negation is omitted as it is equivalent to $(-)\To\bot$):
\begin{equation}\def\arraystretch{1.5}
  \begin{array}{c||c||c||c||c||c||c||c}
    \top&\bot&\vee&\wedge&\To&\forall x\in{}A_i&\exists x\in{}A_i&=\\
    \hline
    \tm&\emptyset&+&\times&(-)^{(-)}&\prod_{x\in{}A_i}&\sum_{x\in{}A_i}&\Delta\\
  \end{array}\label{eq:props-as-types-chart}
\end{equation}
The first two items in the bottom row (one-element set and empty set) are nullary operations, and the following three (disjoint union, product, and exponential (set of functions)) are binary operations. All of these act on families ``point-wise'', e.g., $\set{X_i}_{i\in I}\times\set{Y_i}_{i\in I}:=\set{X_i\times Y_i}_{i\in I}$.

The next two, the indexed product and disjoint union, take a family of sets indexed by $A_1\times\cdots\times A_n$ and return a family of sets indexed by $A_1\times\cdots\times\widehat A_{i}\times\cdots\times A_n$ (where $\widehat A_i$ indicates that the $i$-th entry is omitted).

We will return to the last entry shortly.

Now, given an interpretation $M\colon\sigma\to\Set$, the propositions-as-sets semantics should assign a set $M(\phi)$ to each sentence (i.e., closed formula) $\phi$, and more generally, to a formula $\phi$ with free variables $\vec{x}=\seq{x_1,\ldots,x_n}$ with sorts $\vec{A}=\seq{A_1,\ldots,A_n}$, it should assign a family $M(\phi)$ of sets indexed by the product $M(\vec{A})\underset{\mathrm{def}}=M(A_1)\times\cdots\times M(A_n)$ of the interpretations $M(A_i)$ under $M$ of the sorts $A_i$.

This is done, of course, by recursion on the complexity of $\phi$, and the various recursive clauses are handled according to the table (\ref{eq:props-as-types-chart}). It is then straightforward to verify the desired property that the set $M(\phi)$ is inhabited if and only if $M\vDash\phi$, i.e., the structure $M$ satisfies $\phi$ in the classical sense.

This direct description of the propositions-as-sets semantics is precisely what comes out in the case $\fib{C}=\fib{F(\Set)}$ from the general notion of ``interpretation in the $h^=$-fibration $\fib{C}$''. In the case $\fib{C}=\fib{F(\B)}$ for a general locally bicartesian category $\B$, the description is the same: we observe each of the set-theoretic operations in (\ref{eq:props-as-types-chart}) is governed by a certain universal property in $\Set$, and we take the corresponding operation in the category $\B$ (which exist by locally bicartesian closedness). In this case, ``family of sets indexed by $M(\vec A)$'' becomes ``object in the slice category $B/M(\vec A)$''.

In the general case, the reason that this interpretation is sound for intuitionistic logic is that under the Lambek-Lawvere axiomatization of intuitionistic logic (see \S\ref{subsec:deductions} in the \hyperref[sec:appendix]{appendix}), the logical axioms governing each connective correspond precisely to the universal property defining the operations (\ref{eq:props-as-types-chart}) in a locally bicartesian closed category -- and more generally in an $h^=$-fibration.

We have yet to discuss the interpretation, both for $\Set$ and for general $\B$, of \emph{atomic} formulas -- where in this paper, we consider only \emph{algebraic} signatures (meaning those containing no primitive relations), so that the only atomic formulas are equalities. This brings us to the last column in the table (\ref{eq:props-as-types-chart}).

In the case of propositions-as-sets, we interpret equality in the obvious way: given two terms $s$ and $t$ of sort $A$ which have been interpreted as some elements $M(s)$ and $M(t)$ of the set $M(A)$, we define $M(s=t):=\set{x\in\tm\mid M(s)=M(t)}$, i.e. (assuming classical logic):
\[
  M(s=t)=
  \begin{cases}
    \tm&\text{if }s=t\\
    \emptyset&\text{otherwise}
  \end{cases}
\]

The categorical description of this is as follows. If $s$ and $t$ (still assumed to be of sort $A$) have free variables $\vec{x}$ of sorts $\vec{A}$, then $M(s)$ and $M(t)$ are morphisms $M(\vec{A})\to M(A)\times M(A)$, and $M(s=t)$ is the object in $\B/M(\vec{A})$ given by the pullback
\[
  \begin{tikzcd}
    \cdot\ar[r, ""]\ar[d, "M(s=t)"']
    \ar[rd, "\lrcorner", pos=0.05, phantom, shift right=4pt]&[30pt]
    M(A)\ar[d, "\Delta"]\\
    M(\vec{A})\ar[r, "\br{M(s),M(t)}"]&M(A)\times M(A)
  \end{tikzcd}
\]
of the diagonal morphism $\Delta\colon M(A)\to M(A)\times M(A)$ (hence the ``$\Delta$'' in the table (\ref{eq:props-as-types-chart})).

In particular, this is how equality is interpreted in the $h^=$-fibration $\fib{F(\sSet)}$, and amounts to strict equality of elements of a simplicial set.

\emph{However, for the homotopical semantics, this is not how we want to interpret equality!} Indeed, the whole point of the homotopical semantics is that we want to interpret equality as the \emph{space of paths}. This is why the homotopical semantics cannot literally be given by the $h^=$-fibration $\fib{F(\sSet)}$.

In order to achieve the desired notion of equality, we need to replace the diagonal morphism with the morphism $M(A)^I\to M(A)\times M(A)$, where $I$ is the ``simplicial interval'', making $M(A)^I$ the space of paths in $M(A)$, and where the displayed morphism takes each path to its two endpoints (it is given by the two morphisms $M(A)^I\to M(A)^{\mathrm{pt}}\cong M(A)$ induced by the endpoint inclusions $\mathrm{pt}\to I$).

Thus, the definition of the homotopical semantics (stated in \S\ref{subsec:interpretations}) is: given an interpretation $M\colon\sigma\to\Kan$,\footnote{From what we have said so far, it seems we can use any interpretation in $\sSet$, and not only in $\Kan$. However, this will give the ``wrong'' semantics in general: it will not have the homotopy-invariance property.} we interpret formulas as in the fibration $\fib{F(\sSet)}$ -- i.e., according to the rules in (\ref{eq:props-as-types-chart}) -- except that we interpret equality using the path-space instead of the diagonal.

A priori, this is no longer an instance of the general notion of ``interpretation in an $h^=$-fibration''; but in fact it is! Namely, for any model category $\C$, there is a certain fibration $\fib{HoF_\fb(\C_\fb)}$ obtained from $\fib{F(\C)}$ -- which we show in \S\ref{subsec:h-fib-of-spaces} is a $h^=$-fibration for suitable model categories $\C$ such as $\sSet$ -- such that the just-described homotopical semantics corresponds precisely to the interpretation in $\fib{HoF_\fb(\sSet_\fb)}$, at least up to homotopy equivalence (see \S\ref{subsec:interpretations}). This fibration was studied in \cite{sentai} (and had been previously introduced independently in \cite{cagnethesis}), where a certain 2-categorical structure on it (and in fact, on any ``$\wedgeq$-fibration'') is introduced, whose significance we will return to below in \S\ref{subsec:invariance-disc}.

There are certain advantages to having our semantics be given by an $h^=$-fibration -- for example, as mentioned above, it automatically gives us the soundness of the interpretation with respect to intuitionistic logic.

But also, as we explain next, the fibrational formulation is what we use to prove the homotopy invariance property.

\subsection{Functorial semantics and invariance}\label{subsec:invariance-disc}
To prove the homotopy-invariance property of the homotopical semantics, we must bring in an important aspect of the fibrational semantics that we have yet to discuss, namely that the syntax of first-order logic over the signature $\sigma$ can itself be organized into an $h^=$-fibration $\fib{Pf_\sigma}$ in such a way that the resulting interpretation in a general $h^=$-fibration $\fib{C}$ is then mediated by a morphism of $h^=$-fibrations $\fib{Pf_\sigma}\to\fib{C}$.

This is the notion of ``functorial semantics'' introduced by Lawvere \cite{lawverethesis} at the dawn of categorical logic. Of course, the origin of this idea is in the ``algebraic semantics'' for propositional logic using Boolean (or, in the intuitionistic case, Heyting) algebras.

The fibration $\fibr{Pf_\sigma}{Pf_\sigma}{Tm_\sigma}$ (which is constructed in the \hyperref[sec:appendix]{appendix}) is given roughly as follows. The base category $\TM_\sigma$ is the finite product category associated by Lawvere to (the empty theory over) $\sigma$: the objects are ``contexts'' -- i.e., finite sequences of sorts of $\sigma$ -- and the morphisms are given by sequences of terms of $\sigma$. The objects of $\Pf_\sigma$ are first-order formulas over $\sigma$, and in particular, the fiber $\fib{Pf_\sigma}^{\vec{A}}$ over a context $\vec{A}$ has as objects the formulas with free variables in the context $\vec{A}$, and the morphisms $\phi\to\psi$ are certain equivalence classes of intuitionistic deductions (proofs) of $\phi\To\psi$.\footnote{There is a simpler version of this fibration in which the fibers are replaced by their \emph{posetal reflections}, so that they are instead Heyting algebras, with the ordering given by implication. This version suffices for interpretations in $h^=$-fibrations $\fib{C}$ which are themselves ``posetal'' in this sense, but this is not the case for the fibration $\fib{HoF_\fb(C_\fb)}$ of interest.}

The central feature of the structure $\fib{Pf_\sigma}$ (and more generally, of the algebraic structure constructed out of the syntax in the various flavours of algebraic/functorial semantics) is that it is \emph{free} in an appropriate sense. The freeness precisely captures the desired property that each interpretation $\sigma\to\fib{C}$ in an $h^=$-fibration induces a unique (up to isomorphism) morphism $\fib{Pf_\sigma}\to\fib{C}$ of $h^=$-fibrations. (In fact, for our purposes, we will only need the abstract existence of such a free $h^=$-fibration, and not the particular description in terms of syntax given above.)

However, as is typical of categorical structures, the freeness $\fib{Pf_\sigma}$ has an additional element: given an $h^=$-fibration $\fibr{C}CB$ and two interpretations $M_1,M_2\colon\sigma\to\B$, any \emph{isomorphism} of interpretations $M_1\toi M_2$ gives rise to an isomorphism between the induced morphisms $\fib{Pf_\sigma}\to\fib{C}$. This immediately implies the \emph{isomorphism invariance} of the interpretation in $\fib{C}$; for instance, when $\fib{C}=\fib{F(\Set)}$, this says that for given an isomorphism $M_1\toi M_2$ of set-based structures $M_1,M_2\colon\sigma\to\Set$ (giving, in particular, a bijection $\alpha_{\vec{A}}\colon M_1(\vec{A})\to M_2(\vec{A})$ for each sequence of sorts $\vec{A}$), there is an induced bijection $M_1(\phi)\to M_2(\phi)$ lying over $\alpha_{\vec{A}}$ for each formula $\phi$ with free variables $\vec{x}$ with sorts $\vec{A}$. In particular, $M_1$ and $M_2$ satisfy the same sentences.

This is the recipe for our proof of the homotopy invariance, however things are somewhat more subtle in the latter case. The point is that we are now interested in a \emph{homotopy equivalence} of two structures $M_1,M_2\colon\sigma\to\Kan$, which is no longer a purely category-theoretic notion, and therefore cannot be expected to induce an isomorphism (or rather, homotopy equivalence) of the corresponding morphisms $\fib{Pf_\sigma}\to\fib{HoF_\fb(\Kan)}$ as above.

On the other hand, a homotopy equivalence is naturally formulated as a \emph{2-categorical} notion with respect to a natural 2-categorical structure on $\Kan$ -- namely, it is given by a ``pseudonatural equivalence'' of functors $\TM_\sigma\to\Kan$.

Now, in \cite{sentai}, we showed that for \emph{any} $h^=$-fibration (or more generally $\wedgeq$-fibration) $\fibr{C}CB$, the base category $\B$ automatically inherits a 2-categorical structure -- and in fact, the whole fibration $\fib{C}$ becomes a ``1-discrete 2-fibration'' -- in such a way that this recovers the usual 2-categorical structure on $\Kan$ (or more generally the category of cofibrant-fibrant objects in any model category).

With this preparation, we are able to formulate and prove the homotopy-invariance property as a purely abstract, categorical theorem (Theorem~\ref{thm:abstract-invariance}, stated in the introduction) concerning morphisms of $h^=$-fibrations from a free $h^=$-fibration into a 1-discrete 2-fibration.

Finally, we note that one needn't necessarily deduce the homotopy-invariance property of the homotopical semantics from the abstract invariance theorem; there is presumably a simpler, direct proof by induction. However, the abstract invariance theorem is interesting in itself, and helps to situate the syntactic invariance theorem in a broader context.

\section{Homotopical semantics}
In this section, we give the definition of the homotopical semantics that was sketched in the introduction. As stated there, this is essentially (i.e., up to homotopy equivalence) obtained as a special case of interpreting logic in an $h^=$-fibration, for a particular $h^=$-fibration $\fib{HoF_{\fb}(\Kan)}$.

We first define the latter fibration in \S\ref{subsec:h-fib-of-spaces}. Then, in \S\ref{subsec:interpretations}, we introduce the general notion of interpreting first-order logic in an $h^=$-fibration, define the homotopical semantics and describe its relationship to the interpretation in $\fib{HoF_{\fb}(\Kan)}$. In \S\ref{subsec:htpy-equiv}, we provide the remaining notions needed to state the special invariance theorem, in particular that of homotopy equivalence of interpretations. Finally, in \S\ref{subsec:top-spaces}, we discuss how topological spaces (as opposed to simplicial sets/Kan complexes) can be handled.

\subsection{Preliminaries on $h^=$-fibrations}
We will use the definitions and notation concerning fibrations from \cite{sentai}. In particular, we take everything from \S1 of \emph{op. cit.} for granted, and our present discussion of fibrations will, so to speak, continue from there. We will also use material from other parts of \emph{op. cit.}, but we will indicate when we do so.

\subsubsection{}\defn
A fibration $\fib{C}$ has \emph{fiberwise finite coproducts} if every fiber of $\fib{C}$ has finite coproducts. We have the notion of a coproduct diagram in a fiber of $\fib{C}$ being \emph{stable} (under pullbacks), analogous the corresponding notion for products (see \cite[\S1.5]{sentai}).

We follow the conventions concerning coproducts from \cite[\S10.2]{sentai}, except that we use the symbols $\vee$ and $\bot$ instead of $+$ and $\init$ when the category under consideration is the fiber of some fibration. Also, we denote by $\inn_1$ and $\inn_2$ the coprojections into a coproduct.

Next, given objects $B,C$ in a category $\C$, an \emph{exponential diagram based on $B$ and $C$} consists of an object $C^B$ (the corresponding \emph{exponential object}), a product diagram $C^B\xot{\pi_1}C^B\times B\tox{\pi_2}B$ and a morphism (the \emph{evaluation morphism}) $\varepsilon\colon C^B\times B\to C$ satisfying the usual universal property (see, e.g., \cite[p.~45]{elephant1}).

A category is \emph{cartesian closed} if it has finite products, and there is an exponential diagram based on each pair of objects. It is \emph{bicartesian closed} if it is cartesian closed and has finite coproducts. A functor between cartesian closed categories is \emph{cartesian closed} if it preserves finite products and takes exponential diagrams to exponential diagrams, and a \emph{bicartesian closed} functor is defined similarly.

We will generally use the above notation for exponential objects, except when the category in question is the fiber of some fibration, in which case we will write $B\To{}C$ instead of $C^B$.

A $\wedge$-fibration $\fib{C}$ has \emph{fiberwise exponentials} if each fiber of $\fib{C}$ is cartesian closed. We have the notion of \emph{stability (under pullbacks)} of exponential diagrams in fibers, analogous to that of product and coproduct diagrams.

\subsubsection{}
Let $\fibr{C}CB$ be a fibration. We recall the notion of indexed products and sums in $\fib{C}$.

Given a morphism $f\colon{}A\to{}B$ in $\B$ and $P\in\fib{C}^A$, a \emph{$\prod_f$-diagram based on $P$} consists of an object $\prod_fP\in\fib{C}^B$ (the corresponding \emph{indexed product object}), a cartesian morphism $\ct\colon f^*\prod_fP\to P$ over $f$ and an \emph{evaluation morphism} $\varepsilon\colon f^*\prod_fP\to P$ over $\id_A$, satisfying the appropriate universal property (see \cite[p.~341]{makkai-lauchli1}).

We refer to \emph{loc. cit.} for the notion of an indexed product $\prod_fP$ (or more precisely, of the corresponding $\prod_f$-diagram) being \emph{stable with respect to a pullback square}
\[
\begin{tikzcd}
C\ar[r, "g"]\ar[d, "h"']\ar[rd, phantom, "\lrcorner", pos=0]&D\ar[d, "k"]\\
A\ar[r, "f"]&B
\end{tikzcd}
\]
(this is also known as the \emph{Beck-Chevalley condition}, and amounts to a certain morphism $\prod_g(h^*P)\to k^*\prod_fP$ being an isomorphism).

We say that $\prod_fP$ is \emph{stable along a morphism} $k\colon D\to B$ if it is stable along each pullback square based on $k$ as above.

By an \emph{indexed sum object} $\sum_fP$, we simply mean the codomain of a cocartesian morphism over $f$ with domain $P$. As in \cite{sentai}, when $f$ is a diagonal morphism and $P$ a terminal object, we will usually write $\Eq_A$ in place of $\sum_fP$.

The stability or Beck-Chevalley condition for indexed sums is the same as the one for cocartesian morphisms, which was explained in \cite[\S1.6]{sentai} (cf. \cite[p.~342]{makkai-lauchli1}).

\subsubsection{}\label{defn:h-fibration}\defn
A fibration $\fibr{C}CB$ is a \emph{$h$-fibration} if it satisfies the following four conditions.
\begin{enumerate}[(i),itemsep=0pt,topsep=-3pt]
\item $\fib{C}$ has stable fiberwise finite products and coproducts and exponentials.
\item $\B$ has finite products.
\item\label{item:hfib-def-lift-cond}  For any product projection $\pi\colon{}A\times{}B\to{}B$ in $\B$ and any $P\in\fib{C}^{A\times{}B}$, there exists indexed sums and products $\sum_{\pi}P$ and $\prod_{\pi}P$, and these (i.e., the corresponding cocartesian morphisms and $\prod_\pi$-diagrams) are stable along all morphisms $k\colon D\to B$.
\end{enumerate}

$\fib{C}$ is an $h^=$-fibration if, in addition
\begin{enumerate}[(i),resume,topsep=0pt]
\item For each $B\in\B$, there exists an equality object $\Eq_B$ (i.e., a cocartesian lift of a diagonal $\Delta_B\colon{}B\to{}B\times{}B$ with domain a terminal object $\top_{B}\in\fib{C}^B$).
\end{enumerate}

\subsubsection{}\label{defn:morphism-of-fibs}\defn
Given prefibrations $\fibr{C}CB$ and $\fibr{C'}{C'}{B'}$, a \emph{morphism of prefibrations} $\fib{C}\to\fib{C}'$ is a pair $(\Phi,\phi)$, where $\phi\colon\B\to\B'$ and $\Phi\colon\C\to\C'$ are functors such that the square
\[
  \begin{tikzcd}
    \C\ar[r, "\Phi"]\ar[d, "\fib{C}"']&\C'\ar[d, "\fib{C'}"]\\
    \B\ar[r, "\phi"]&\B'
  \end{tikzcd}
\]
commutes (strictly). We say that $(\Phi,\phi)$ is a morphism of prefibrations \emph{over} $\phi$. If $\B=\B'$ and $\phi=\id_\B$, we may just write $\Phi$ instead of $(\Phi,\id_B)$, and we say in this case that $\Phi$ is \emph{over $\B$}.\footnote{In \cite[Definition~11.5]{sentai}, what we now call a ``morphism of prefibrations over $B$'' was just called a ``morphism of prefibrations''. We now adopt the new terminology.} If $\fib{C}$ and $\fib{C'}$ are fibrations, then $(\Phi,\phi)$ is a \emph{morphism of fibrations} if $\Phi$ takes cartesian morphisms to cartesian morphisms.

Note that for each $A\in\B$, $(\Phi,\phi)$ induces a functor $\Phi\colon\fib{C}^A\to(\fib{C'})^{\phi{}A}$.

If $\fib{C}$ and $\fib{C'}$ are $*$-fibrations (where $*$ is one of $\wedge$, $\wedgeq$, $h$, $h^=$ -- see \cite[\S1]{sentai} for the first two notions), we say that $(\Phi,\phi)$ is a \emph{morphism of $*$-fibrations} if it preserves the relevant structure: (i) in all cases, the induced functors on fibers should be f.p.\ (ii) if $*$ is $h$ or $h^=$, the induced functors on fibers should moreover be bicartesian closed, and $\Phi\colon\C\to\C'$ should preserve $\prod_\pi$-diagrams and cocartesian morphisms over product projections $\pi$ (iii) if $*$ is $\wedgeq$, $h$, or $h^=$, $\phi\colon\B\to\B'$ should preserve finite products (iv) if $*$ is $h^=$ or $\wedgeq$, $\Phi$ should preserve cocartesian lifts of diagonal morphisms with domain a terminal object.

Next, given two morphisms $(\Phi,\phi),(\Psi,\psi)\colon\fib{C}\to\fib{C'}$ of pre-fibrations, a \emph{natural transformation} $(\Phi,\phi)\to(\Psi,\psi)$ is a pair $(\tilde \alpha,\alpha)$ of natural transformations $\alpha\colon\phi\to\psi$ and $\tilde\alpha\colon\Phi\to\Psi$ such that $\fib{C'}\circ \tilde{\alpha}=\alpha\circ\fib{C}$. Again, we say that $(\tilde \alpha,\alpha)$ \emph{lies over} $\alpha$, or \emph{over} $\B$ if $\B=\B'$, $\phi=\id_\B$, and $\alpha=\id_{\id_B}$. We say that $(\tilde\alpha,\alpha)$ is an \emph{equivalence} if $\alpha$ and $\tilde\alpha$ are both equivalences.

\subsubsection{}\prop\label{prop:hfib-is-wedgeq}
Every $h^=$-fibration is a $\wedgeq$-fibration.
\pf
Referring to the definition of $\wedgeq$-fibration from \cite[\S1.7]{sentai}, we see that we are only missing the Frobenius reciprocity and stability along product projections for equality objects. These follow from certain well-known facts, which we leave to the reader, namely (i) that the existence of stable exponential objects implies Frobenius reciprocity (see, e.g., \cite[p.~6]{lawvereequality}\cite[p.~343]{makkai-lauchli1}\cite[p.~102]{jacobscatlogic}) and (ii) that, given $f\colon A\to B$ in the base category (here, we take $f$ a product projection), if indexed products $\prod_fP$ exist for all $P$ lying over $A$, then all indexed sums are stable along $f$ (see, e.g., \cite[p.~343]{makkai-lauchli1}\cite[Lemma~1.9.7]{jacobscatlogic}).
\qed

\subsubsection{}\label{defn:lccc}\defn
A category $\C$ is \emph{locally cartesian closed} if each slice category $\C/X$ is cartesian closed and $\C$ has a terminal object (so in particular, $\C\cong{}\C/\tm$ is itself cartesian closed), and it is \emph{locally bicartesian closed} if it is locally cartesian closed and has finite coproducts (equivalently, the slices are bicartesian closed and $\C$ has a terminal object).

It is well-known that $\fib{F(\C)}$ is an $h$-fibration (in fact, an $h^=$-fibration) if and only if $\C$ is locally bicartesian closed (see, e.g, \cite[p.~345]{makkai-lauchli1}\cite[\S2.4]{seely-loc-cart-closed}\cite[p.~81]{jacobscatlogic}).

\subsubsection{}\label{rmk:adjoints}\rmk
It is easy to see (and well-known -- see \cite{makkai-lauchli1}) that, given a fibration $\fibr{C}CB$ and a morphism $f\colon{}A\to{}B$ in $\B$, the existence of indexed sums or products over $f$ is equivalent to the existence of left and right adjoints $\sum_f$ and $\prod_f$ to ``the'' pullback functor $f^*\colon\fib{C}^B\to\fib{C}^A$.

In general, this requires the axiom of choice.\footnote{We note that this and the various other uses of the axiom of choice which arise in this paper are needed only when considering fibrations in general, but are not needed in relation to our specific fibration of interest $\fib{HoF_{\fb}(\Kan)}$, in which we can uniformly specify all the required operations.} One way to formulate this which does not involve any choices is to say there is always a canonical \emph{anafunctor} (see \cite{avoiding-choice}) $f^*\colon\fib{C}^B\to\fib{C}^A$, and the existence of indexed sums or products over $f$ is equivalent to the existence of a left or right adjoint anafunctor $\sum_f$ and $\prod_f$ (which is then also canonically defined).

We also note that in the case $\fib{C}=\fib{F(\C)}$, there is always an explicit indexed sum functor $\sum_f\colon\C/A\to\C/B$ taking $(X,x)$ to $(X,fx)$.

\subsection{The $h^=$-fibration of spaces}\label{subsec:h-fib-of-spaces}
We now introduce the $h^=$-fibration $\fib{HoF_{\fb}(\Kan)}$; or rather, we show that it is an $h^=$-fibration, as it was already introduced in \cite{sentai}. There, it was already shown that it was a $\wedgeq$-fibration. Hence, to see that it is an $h^=$-fibration, it remains to show that this fibration has the necessary extra structure -- namely, that it supports the ``logical'' operations $\vee,\To,\forall,\exists$.

In fact, in \emph{op. cit.}, it was shown more generally that the fibration $\fib{HoF_{\fb}(\C_\fb)}$ is a $\wedgeq$-fibration for \emph{any} model category $\C$. We recall that the first subscript ``$\fb$'' means we are restricting to fibrations, and the second that we are restricting to fibrant objects (removing the first ``$\fb$'' actually results in an equivalent fibration, in contrast to the situation with the related fibration $\fib{F_\fb(\C_\fb)}$). We have $(\sSet)_\fb=\Kan$, which is why we write $\fib{HoF_\fb(\Kan)}$.

It is not true that $\fib{HoF_{\fb}(\C_\fb)}$ is a $h^=$-fibration for any model category $\C$; we must put certain restrictions on $\C$ (satisfied, of course, by $\sSet$). We encapsulate the needed conditions in Definition~\ref{defn:suitable} (``suitable model category''), similar in spirit to the notion of ``type-theoretic model category'' and related notions (see, e.g., \cite[p.~166]{lumsdaine-shulman-hits}).
We do this just to clarify what assumptions we are using; we do not actually consider any suitable model categories besides $\sSet$, though we suspect there are other interesting ones -- we have in mind the so-called \emph{Cisinski} model structures, which always satisfy conditions \ref{item:suitable-lcc} and \ref{item:suitable-LAST-cof-monos} of the definition.

The approach taken here is the same as in \cite[\S\S12-13]{sentai}. We first show that the fibration $\fib{F_\fb(\C_\fb)}$ -- from which $\fib{HoF_{\fb}(\C_\fb)}$ is obtained by passing to the homotopy category in each fiber -- is a $h$-fibration (in fact, it is an $h^=$-fibration, but we are not interested in its equality structure), and then show that this structure is preserved upon passing to homotopy categories. Along the way, we deduce the important fact that the canonical morphism $\fib{F_\fb(\C_\fb)}\to\fib{HoF_\fb(\C_\fb)}$, as well as the inclusion $\fib{F_\fb(\C_\fb)}\to\fib{F(\C)}$, are morphisms of $h$-fibrations, which gives us an explicit description of the $h$-fibration structure in $\fib{HoF_\fb(\C_\fb)}$.

\subsubsection{}\defn
The category $\sSet$ of simplicial sets is by definition the category $\Set^{\cat{\Delta}^\op}$ of presheaves on the category $\Delta$ of non-empty finite ordinals and non-decreasing maps. There is a standard (``Quillen'') model structure on $\sSet$ (see \cite[II.3.14]{quillen-ha}). Its fibrant objects are called \emph{Kan complexes}, and we set $\Kan:=\sSet_\fb$ for the category of Kan complexes.

Outside of the discussion in \S\ref{subsec:top-spaces} and the examples in \S\ref{subsec:examples}, we will not really use anything about this model category, except what little is needed to show that it is suitable in the sense of Definition~\ref{defn:suitable}. However, we certainly expect the reader to be familiar with it. For example, it is good to know that there is an adjunction (in fact, ``Quillen equivalence'') $\sSet\leftrightarrows\Top$ (the functors being ``geometric realization'' and ``singular simplicial set'') inducing an equivalence $\Ho(\sSet)\cong\Ho(\Top)$. It might also be good to know (for the sake of general cultural context) that Kan complexes can justifiably be said to provide a formal definition of the (informal) notion of ``weak infinity-groupoid''.

\subsubsection{}\label{defn:suitable}\defn
A model category $\C$ is \emph{suitable} if the following four conditions are satisfied.
\begin{enumerate}[(i),itemsep=0pt,topsep=-3pt]
\item\label{item:suitable-FIRST-rp} $\C$ is \emph{right-proper}, i.e., weak equivalences are closed under pullbacks along fibrations.\footnote{Actually we will only need that trivial cofibrations are closed under pullbacks along $!_A\colon A\to\tm_\C$ for $A$ fibrant.}
\item\label{item:suitable-lcc} $\C$ is locally cartesian closed (as a category).
\item\label{item:suitable-fib-coprods} $[f,g]\colon{}A+B\to{}C$ is a fibration whenever $f\colon{}A\to{}C$ and $g\colon{}B\to{}C$ are, and the unique morphism $0\to{}A$ is always a fibration.
\item\label{item:suitable-LAST-cof-monos}
  The cofibrations of $\C$ are precisely the monomorphisms.
\end{enumerate}

We note that the Quillen model structure on simplicial sets is suitable. Condition~\ref{item:suitable-LAST-cof-monos} holds by definition. Condition~\ref{item:suitable-FIRST-rp} is non-trivial, but well-known, and follows from the existence of pullback and fibration preserving fibrant replacement functors (see, e.g. \cite[p.~370]{may-ponto}). As for \ref{item:suitable-lcc}, it is well-known that any presheaf category is locally cartesian closed (see \cite[p.~48]{elephant1}).

To see \ref{item:suitable-fib-coprods}, note that the ``horns'' $\Lambda^n_k$ are all \emph{connected}, in the sense that any two vertices are connected by a path of edges (in fact, for $n>2$, by a single edge). It follows that any morphism $\Lambda^n_k\to{}A+B$ to a coproduct must factor through one of the summands $A,B$. Now, for $[f,g]\colon{}A+B\to{}C$ to be a fibration, it must lift against each horn inclusion $\Lambda^n_k\to{}\Delta^n$. But if $f$ and $g$ are each a fibration, this follows immediately from the fact that any given morphism $\Lambda^n_k\to{}A+B$ factors through $A$ or $B$. The corresponding lifting problem for morphisms $0\to{}A$ is trivial, since there are no morphisms $\Lambda^n_k\to0$.

\subsubsection{}\label{prop:slice-suitable}\prop
If $\C$ is a suitable model category, then any slice category $\C/A$ of $\C$, with its induced model structure, is also suitable (and in fact, this is true separately for each of the conditions in Definition~\ref{defn:suitable}).
\pf
It is well-known (and easy to see) that each slice of a locally cartesian closed category is locally cartesian closed; this follows from the existence of the canonical isomorphisms $(\C/A)/B\cong\C/B$.
That Conditions~\ref{item:suitable-fib-coprods}~and~\ref{item:suitable-LAST-cof-monos} hold in $\C/A$ follows from the fact that a morphism $(p,\id_A)$ in $\C/A$ is a monomorphism, cofibration, or fibration if and only if $p$ is, and the fact that the forgetful functor $\C/A\to\C$ preserves coproducts.
Similarly, Condition~\ref{item:suitable-FIRST-rp} holds since a square in a slice category $\C/A$ is a pullback square if and only if its image under the forgetful functor $\C/A\to\C$ is.
\qed

\subsubsection{}\prop\label{prop:suitable-all-cofib}
In a suitable model category $\C$, every object is cofibrant. Hence $\C_\cfb=\C_\fb$, and $\Ho(\C_\fb)=\pi(\C_\cfb)$.
\pf
Since the cofibrations are the monomorphisms, this amount to checking that each morphism from the initial object is a monomorphism. It is well-known that this holds in any cartesian closed category (see \cite[p.~61]{elemelem}).
\qed

\subsubsection{}\label{prop:suitable-fib-exp}\prop
Let $\C$ be a suitable model category, let $p\colon{}X\to{}Y$ be a fibration, and let $C$ be a fibrant object. Then the induced map $p^C\colon{}X^C\to{}Y^C$ is a fibration.
\pf
We must show that for any solid commutative diagram
\[
\begin{tikzcd}
A\ar[d, "i"']\ar[r]&X^C\ar[d, "p^C"]\\
B\ar[r]\ar[ru, dashed]&Y^C
\end{tikzcd}
\]
in $\C$ with $i$ a trivial cofibration, there exists a dashed morphism making the diagram commute.

Using the adjunction $(-\times{}C)\dashv(-)^C$, this is seen to be equivalent to the existence of a dashed morphism making the corresponding diagram
\[
\begin{tikzcd}
A\times{}C\ar[d, "i\times\id_C"']\ar[r]&X\ar[d, "p"]\\
B\times{}C\ar[r]\ar[ru, dashed]&Y
\end{tikzcd}
\]
commute. For this, it suffices that $i\times\id_C$ be a trivial cofibration. Noting that it is the pullback of $i$ along the projection $B\times{}C\to{}B$ (which is a fibration since $C$ is fibrant), we have that $i\times\id_\C$ is a cofibration since monomorphisms are stable under pullback, and it is a weak equivalence since $\C$ is right-proper.
\qed

\subsubsection{}\label{prop:suitable-fb-bcc}\prop
Let $\C$ be a suitable model category. We know from \cite[Proposition~12.2]{sentai} that $\C_\fb$ is an f.p.\ category and that the inclusion $\C_\fb\hookrightarrow\C$ is an f.p.\ functor.

We now claim that $\C_{\fb}$ and the inclusion $\C_{\fb}\hookrightarrow\C$ are bicartesian closed.
\pf
Since $\C_\fb$ is a full subcategory of $\C$, it suffices to show that the fibrant objects in $\C$ are closed under exponentials and coproducts. Condition~\ref{item:suitable-fib-coprods} of ``suitable'' (Definition~\ref{defn:suitable}) implies that the fibrant objects are closed under coproducts. That they are closed under exponentials follows from Proposition~\ref{prop:suitable-fib-exp}.
\qed

\subsubsection{}\label{defn:fib-restrictions}\defn
In \cite[Definition~13.6]{sentai}, for a model category $\C$ and a full subcategory $\D\subseteq\C$, we defined $\fib{HoF_*(\D)}$ to be the restriction of $\fib{HoF_*(\C)}$ to $\D$ (where $*\in\set{\cf,\fb,\cfb}$). We similarly define $\fib{F_*(\D)}$ to be the restriction of $\fib{F_*(\C)}$ to $\D$.

We note that, in general, the restriction of a $*$-fibration (with $*$ one of $\wedge$, $\wedgeq$, $h$, $h^=$) to any full subcategory having finite products is again an $*$-fibration.

\subsubsection{}\label{prop:suitable-mfib-hfib}\prop
Let $\C$ be a suitable model category. Since $\C$ is locally bicartesian closed, we know that $\fib{F(\C)}$, and hence $\fib{F(\C_\fb)}$, is an $h$-fibration (see Definition~\ref{defn:lccc}), and by \cite[Proposition~12.3]{sentai}, we know that $\fib{F_\fb(\C_\fb)}$ is a $\wedge$-fibration, and that the inclusion $\fib{F_\fb(\C_\fb)}\hookrightarrow\fib{F(\C_\fb)}$ is a morphism of $\wedge$-fibrations.

We now claim that $\fib{F_\fb(\C_\fb)}$ is an $h$-fibration and the inclusion $\fib{F_\fb(\C_\fb)}\hookrightarrow\fib{F(\C_\fb)}$ is a morphism of $h$-fibrations.

\pf
It follows from Propositions~\ref{prop:slice-suitable}~and~\ref{prop:suitable-fb-bcc} that the fibers of $\fib{F_\fb(\C_\fb)}$, and the functors on fibers induced by the inclusion, are bicartesian closed.

Next, we need to check that, given a product projection $\pi_2\colon{}A\times{}B\to{}B$ in $\C_{\fb}$ and a cocartesian morphism $(p,\pi_2)\colon(X,A\times{}B,x)\to{}(Y,B,y)$ in $\C^\to$ lying over $\pi_2$, if $(X,x)$ is in $\fib{F_\fb(\C_\fb)}^{A\times{}B}$, then $(Y,y)$ is in $\fib{F_\fb(\C_\fb)}^B$ -- i.e., if $x$ is a fibration, then so is $y$. But the product projection $\pi_2$ is a fibration since $A$ is fibrant, and $p$ is an isomorphism (this being equivalent to the cocartesianness of $(p,\pi_2)$), hence $y=\pi_2xp\I$ is a fibration as well.

Similarly, we need to check that if $(X,x)$ is a fibration, then $\prod_{\pi_2}(X,x)$ is a fibration. In the case where $B\cong\tm_\C$, it is well-known (e.g., \cite[Lemma~1.5.2]{elephant1}) that $\prod_{\pi_2}(X,x)$ can be computed as a pullback along a morphism $\tm_\C\to A^A$ of $x^A\colon X^A\to A^A$, which is a fibration by Proposition~\ref{prop:suitable-fib-exp}.

We now reduce the general case to this one by using the canonical isomorphisms $(\C/A\times{}B)\cong(\C/B)/(A\times{}B,\pi_2)$ and $\C/B\cong(\C/B)/(B,\id_B)$ as follows. We wish to show that the (ana-)functor $\prod_{\pi_2}\colon\C/(A\times B)\to\C/B$ (see Remark~\ref{rmk:adjoints}) preserves fibrant objects, which we know in the case when $B$ is terminal. Consider the following diagram on the left.
\[
  \begin{tikzcd}
    \C/(A\times B)\ar[r, "\prod_{\pi_2}"]&\C/B\\
    (\C/B)/(A\times B,\pi_2)\ar[r, "\prod_{\pi_2}"']\ar[u, "\sim" sloped, -]&
    (\C/B)/(B,\id_B)\ar[u, "\sim" sloped, -]
  \end{tikzcd}
  \quad
  \begin{tikzcd}
    \C/(A\times B)\ar[r, "\sum_{\pi_2}"]&\C/B\\
    (\C/B)/(A\times B,\pi_2)\ar[r, "\sum_{\pi_2}"']\ar[u, "\sim" sloped, -]&
    (\C/B)/(B,\id_B)\ar[u, "\sim" sloped, -]
  \end{tikzcd}
\]
The bottom row preserves fibrant objects since $(B,\id_B)$ \emph{is} terminal in $\C/B$. Hence, (since the vertical isomorphisms clearly preserve fibrant objects,) we will be done if we can show that the diagram commutes. But note that the diagram to the right clearly commutes (with $\Sigma_{\pi_2}$ the explicit indexed sum functor from Remark~\ref{rmk:adjoints}), hence the left diagram as well, since $\Pi_{\pi_2}$ is a right (ana-)adjoints to a right (ana-)adjoint $\pi_2^*$ of $\Sigma_{\pi_2}$.

It remains to see that all the operations are ``stable'', i.e., that the pullback functors are bicartesian closed, and that the cocartesian morphisms and $\prod_\pi$-diagrams over product projections $\pi$ are stable. In each case, this follows immediately from the corresponding fact in $\fib{F(\C_\fb)}$.
\qed

\subsubsection{}\label{prop:suitable-hofb-bcc}\prop
Let $\C$ be a suitable model category. We know from Proposition~\ref{prop:suitable-fb-bcc} that $\C_{\fb}$ is bicartesian closed, and we know from \cite[Proposition~12.4]{sentai} that $\Ho(\C_\fb)$ is an f.p.\ category and that the functor $\gamma\colon\C_\fb\to\Ho(\C_\fb)$ is an f.p.\ functor.

We now claim that the category $\Ho(\C_{\fb})$ and the functor $\gamma\colon\C_\fb\to\Ho(\C_\fb)$ are bicartesian closed.
\pf
That $\Ho(\C_\fb)$ has, and $\gamma$ preserves, finite coproducts, follows from an argument dual to the one given in \cite[Proposition~12.4]{sentai} (using that all objects in $\C$ are cofibrant).

We next turn to exponentials. Consider an exponential diagram
\[
\begin{tikzcd}[column sep=6pt, row sep=10pt]
&C^B\times{}B\ar[dl, "\pi_1"]\ar[dr, "\pi_2"']\ar[rr, "\expev"]&&C\\
C^B&&B.&
\end{tikzcd}
\]
We already know that $C^B\times{}B$ is still a product in $\Ho(\C_\fb)$, so it remains (by Proposition~\ref{prop:suitable-all-cofib}) to see that for any $A\in\C_\fb$ and product $A\xot{\pi_1}A\times{}B\tox{\pi_2}B$, the composite
\[
\pi(A,C^B)\tox{(-)\times\id_B}\pi(A\times{}B,C^B\times{}B)\tox{(\expev\circ\text{--})}\pi(A\times{}B,C)
\]
is a bijection. That this map is surjective is clear, since it is already surjective before passing to homotopy classes.
To show injectivity, we need to show that if two morphisms $f_1,f_2\colon{}A\times{}B\to{}C$ are homotopic, then the corresponding morphisms $\expind{f_1},\expind{f_2}\colon{}A\to{}C^B$ are.

Let $A+A\tox{[\partial_1,\partial_2]}A\times{}I\tox{\sigma}A$ be a cylinder object for $A$. Because the functor $(-\times{}B)$ is a left \hbox{(ana-)adjoint}, it preserves coproducts, and hence the canonical morphism $[\br{\inn_1\pi_1,\pi_2},\br{\inn_2\pi_1,\pi_2}]\colon{}A\times{}B+A\times{}B\to(A+A)\times{}B$ is an isomorphism. Applying $(-\times{}B)$ to our cylinder object for $A$, we have a sequence of morphisms
\[
A\times{}B+A\times{}B\tox{\sim}
(A+A)\times{}B\tox{[\partial_1,\partial_2]\times\id_B}
(A\times{}I)\times{}B\tox{\sigma\times\id_B}
A\times{}B
\]
and we claim that this exhibits $(A\times{}I)\times{}B$ as a cylinder object for $A\times{}B$. Indeed, the composite is clearly equal to $\nabla_{A\times{}B}$, and $\sigma\times\id_B$ is a weak equivalence by the right-properness of $\C$, since it the pullback of a weak equivalence along the projection $(A\times{}I)\times{}B\to{}A\times{}I$, which is a fibration since $B$ is fibrant. Moreover, the first two morphisms are cofibrations (the first being an isomorphism and the second being the pullback of a monomorphism and hence a monomorphism).

Hence, by \cite[Proposition~10.5~(ii)]{sentai}, given two homotopic maps $f_1,f_2\colon{}A\times{}B\to{}C$, we obtain a left-homotopy $h\colon(A\times{}I)\times{}B\to{}C$ between them, and hence a morphism $\expind{h}\colon{}A\times{}I\to{}C^B$. It remains to see that this is a homotopy between $\expind{f_1}$ and $\expind{f_2}$, i.e. that $\expind{h}\partial_i=\expind{f_i}\colon{}A\to{}C^B$. It suffices to see that $\expev\cdot((\expind{h}\partial_i)\times\id_B)=f_i$, which follows from the definition of $\tilde h$.
\qed

\subsubsection{}\label{thm:suitable-hfib-hfib}\thm
Let $\C$ be a suitable model category. By Proposition~\ref{prop:suitable-mfib-hfib}, we know that $\fib{F_\fb(\C_\fb)}$ is an $h$-fibration, and by \cite[Propositions~12.5~and~13.7]{sentai}, we know that $\fib{HoF_\fb(\C_\fb)}$ is a $\wedgeq$-fibration and that the localization morphism $\gamma\colon\fib{F_\fb(\C_\fb)}\to\fib{HoF_\fb(\C_\fb)}$ is a morphism of $\wedge$-fibrations.

We now claim that $\fib{HoF_\fb(\C_\fb)}$ is in fact an $h^=$-fibration, and that $\gamma$ is a morphism of $h$-fibrations.
\pf
By Proposition~\ref{prop:suitable-hofb-bcc}, we know that the fibers of $\fib{HoF_\fb(\C_\fb)}$ and the functors on the fibers induced by $\gamma\colon\fib{F_\fb(\C_\fb)}\to\fib{HoF_\fb(\C_\fb)}$ are bicartesian closed.

Next, we consider indexed sums. That is, we need to show that the image under $\gamma$ of any cocartesian morphism in $\fib{F_{\fb}(\C_\fb)}$ over a product projection is cocartesian in $\fib{HoF_\fb(\C_\fb)}$. This follows from \cite[Proposition~13.1]{sentai} and the fact that every isomorphism is a weak equivalence.

We next consider indexed products $\prod_{\pi_2}P$. Let $\pi_2\colon{}A\times{}B\to{}B$ be a product projection in $\C_\fb$, and let
\[
\begin{tikzcd}
\pi_2^*\prod_{\pi_2}P\ar[r, "\ct"]\ar[d, "\expev"']&\prod_{\pi_2}P\\
P\\[-10pt]
A\times{}B\ar[r, "\pi_2"]&B
\end{tikzcd}
\]
be a $\prod_{\pi_2}$-diagram in $\fib{F_\fb(\C_\fb)}$. We need to see that its image in $\fib{HoF_\fb(\C_\fb)}$ is a also a $\prod_{\pi_2}$-diagram.

We know already that the image of $\ct$ is cartesian. Hence, it remains to show that for each $Q\in(\C/B)_\fb$ and cartesian $\ct\colon{}f^*Q\to{}Q$ over $\pi_2$, the composite
\[
\pi(Q,\tprod_{\pi_2}P)
\tox{\pi_2^*}
\pi(\pi_2^*Q,\pi_2^*\tprod_{\pi_2}P)
\tox{(\expev\circ\text{--})}
\pi(\pi_2^*Q,P)
\]
is a bijection. As in the proof of Proposition~\ref{prop:suitable-fb-bcc}, it is immediate that it is surjective, and injectivity follows by a similar argument to the one there.

It remains to check the various ``stability'' conditions for $\fib{HoF_\fb(\C_\fb)}$. These are proven in the same way as the stability of products in \cite[Proposition~12.5]{sentai}. Namely, in each case, we reduce to showing the stability of \emph{some} (rather than every) diagram of the appropriate kind, and then we choose the diagram coming from $\fib{F_\fb(\C_\fb)}$, where we already know that stability holds.
\qed

\subsection{Interpreting logic in fibrations}\label{subsec:interpretations}
In this section, we recall the usual notion of interpreting first-order logic in an $h^=$-fibration, and then describe the variant of it giving the homotopical semantics. We will then prove that, up to homotopy equivalence, the homotopical semantics in a suitable model category $\C$ agrees with the (ordinary) semantics in the $h^=$-fibration $\fib{HoF_\fb(\C_\fb)}$.

We note that in the case of a \emph{posetal} $h^=$-fibration $\fib{C}$ (i.e., one in which the fibers are Heyting algebras), the interpretation in $\fib{C}$ unambiguously assigns to each formula an object in one of the fibers of $\fib{C}$. However, in the general case, there are choices involved, and this assignment is only determined up to isomorphism.

\subsubsection{}\label{defn:signature}\defn
A \emph{(multi-sorted) algebraic signature} $\sigma$ is given by a set $\Ob\sigma$ of \emph{sorts} and, for each finite sequence $\vec{A}$ of sorts and each sort $B$, a set $\sigma(\vec{A},B)$ of \emph{function symbols} (with \emph{arity} $\vec{A}$ and \emph{codomain sort} $B$). We denote the set of finite (possibly empty) sequences in a set $X$ by $X^{<\omega}$, and write $\len{\vec{A}}$ for the length of a finite sequence. Also, we denote concatenation of sequences (or of a sequence with a single element) by juxtaposition. Given a sequence with the name $\vec{X}$, we will denote its entries by $X_1,\ldots,X_{\len{X}}$.

Given a finite product category $\C$, an \emph{interpretation $M$ of $\sigma$ in $\C$} consists of the following data (i)-(iii):
\begin{enumerate}[(i),itemsep=0pt,topsep=-3pt]
\item A function $M\colon\Ob\sigma\to\Ob\C$
\item A choice of object $M\vec{A}$ and product diagram $\set{\pi_i^M\colon{}M\vec{A}\to{}MA_i}_{i=1}^{\len{\vec{A}}}$ on $MA_1,\ldots,MA_{\len{\vec{A}}}$ for each sequence $\vec{A}\in(\Ob\sigma)^{<\omega}$ (where we require that $M\seq{A}=MA$ and $\pi_1^{M}=\id_{MA}\colon{}M\seq{A}\to{}MA$ for each $A\in\Ob\sigma$)
\item A morphism $Mf\colon{}M\vec{A}\to{}MB$ for each $f\in\sigma(\vec{A},B)$.
\end{enumerate}
We write $M\colon\sigma\to\C$ to indicate that $M$ is an interpretation of $\sigma$ in $\C$.

Next, given an interpretation $M\colon\sigma\to\C$ and an f.p.\ functor $F\colon\C\to\D$, we obtain an interpretation $F\circ{}M\colon\sigma\to\D$ by setting
$(F\circ{}M)(A)=F(MA)$ for $A\in\Ob\sigma$;
$(F\circ{}M)(\vec{A})=F(M\vec{A})$ and $\pi_i^{F\circ{}M}=F\pi_i^M\colon(F\circ{}M)(\vec{A})\to(F\circ{}M)A_i$ for $\vec{A}\in(\Ob\sigma)^{<\omega}$ and $1\le{}i\le\len{\vec{A}}$; and
$(F\circ{}M)(f)=F(Mf)$ for $f\in\sigma(\vec{A},B)$.

For the rest of \S\ref{subsec:interpretations}, fix an algebraic signature $\sigma$.

\subsubsection{}\defn\label{defn:syntax}
We recall the (usual) syntax of first-order logic over $\sigma$.

We fix, once and for all, an arbitrary infinite set $\Varn$ of ``variable names'' (for definiteness, we could take $\Varn=\N$).

We next declare the \emph{symbols} to be used in the syntax. These consist of: (i) the \emph{(sorted) variable symbols}, which are given by the set $\Varn\times\Ob\sigma$, (ii) the \emph{functions symbols}, which are just the function symbols of $\sigma$ (i.e., they are given by the disjoint union of all the sets $\sigma(\vec{A},B)$), and (iii) the additional symbols $\wedge$, $\vee$, $\To$, $\forall$, $\exists$, $\top$, and $\bot$ (negation $\neg \phi$ is taken as an abbreviation of $\phi\To\bot$).

Now, on the basis of these, we define the set of \emph{$\sigma$-terms} (or just \emph{terms}), and then of \emph{$\sigma$-formulas} (or just \emph{formulas}) as follows.

Each $\sigma$-term $t$ will have a sort $\tp(t)\in\Ob\sigma$ associated with it (the sort \emph{of} $t$), and the set of $\sigma$-terms is given as follows: (i) each sorted variable symbol $(v,A)$ is a term of sort $A$, and (ii) $ft_1\ldots t_{\len{\vec{A}}}$ is a term of sort $B$ whenever $f\in\sigma(\vec{A},B)$ is a function symbol and $t_1,\ldots,t_{\len{\vec{A}}}$ are terms, with $\tp t_i=A_i$.

Given a sequence of terms $\vec{t}=\br{t_1,\ldots,t_n}$, we will write $\tp(\vec{t})$ for the sequence of sorts $\br{\tp(A_1),\ldots,\tp(A_n)}$.

Next, the $\sigma$-formulas are given as follows: (i) $s=t$ is a formula whenever $s$ and $t$ are terms of the same sort, (ii) $\top$ and $\bot$ are formulas, (iii) $\phi\wedge\psi$, $\phi\vee\psi$, and $\phi\To\psi$ are formulas whenever $\phi$ and $\psi$ are, and (iv) $\forall v\phi$ and $\exists v\phi$ are formulas whenever $\phi$ is a formula and $v$ is a variable symbol.

We take for granted the notions \emph{free} and \emph{bound} variables in a term or formula (all variables in a term being free), as well as the notion of one formula being obtained from another by the \emph{renaming of bound variables}. In fact, we will henceforth identify two formulas when they differ only by renaming of bound variables (this identification is needed, for example, in order to make the operation of \emph{capture-avoiding substitution} -- which we also take for granted -- well defined in general).

We also take for granted the \emph{principles of structural induction and recursion}. That is, to prove something about all terms or formulas, we can proceed by induction, the base cases being the atomic terms or formulas, and the induction step being to prove statement for a given formula after assuming it for its constituents. Similarly, to define a function on the set of all terms or formulas, we can likewise proceed by recursion.

Finally, we define a \emph{formula-in-context} to be a pair $(\phi,\vec{x})$, where $\vec{x}$ is a sequence of distinct variable symbols containing all the free variables of $\phi$, and we define a \emph{term-in-context} similarly.

\subsubsection{}\defn\label{defn:interpretation-over}
Let $\B$ be a finite product category $\B$, and $M\colon\sigma\to\B$ an interpretation. Given a sequence of variables $\vec{x}$ we write $M(\vec{x})$ for $M(\tp\vec{x})$.

We now define a function taking each term-in-context $(t,\vec{x})$ to a morphism $M_{\vec{x}}t\colon M(\vec{x})\to M(\tp t)$ by recursion on $t$ as follows. If $t$ is a variable $x_i$, we set $M_{\vec{x}}(t)=\pi_i^M\colon M(\vec{x})\to M(\tp x_i)$, and if $t=ft_1\ldots t_n$ with $f\in\sigma(\vec{A},B)$ (so that $\tp(t)=B$ and $\tp(t_i)=A_i$), then we take $M_{\vec{x}}(t)$ to be the composite
\[
  M_{\vec{x}}(t)\tox{\br{M_{\vec{x}}(t_1),\ldots,M_{\vec{x}}(t_n)}}
  M(\vec{A})\tox{Mf}MB.
\]
(note here that the free variables of each $t_i$ are among those in $\vec{x}$ since this is the case for $t$).

Next, given an $h^=$-fibration $\fibr{C}CB$ over $B$, an \emph{interpretation $\widehat M$ in $\fib{C}$ over $M$} is a function assigning to each formula-in-context $(\phi,\vec{x})$ an object $\widehat M_{\vec{x}}(\phi)\in\fib{C}^{M(\vec{x})}$, satisfying the following conditions:
\begin{enumerate}[(i),itemsep=0pt,topsep=0pt]
\item\label{item:sigma-interp-over-FIRST-topbot} If $\phi$ is $\top$ or $\bot$, then $\widehat M_{\vec{x}}(\phi)$ must be a terminal object $\top_{M(\vec{x})}$ or initial object $\bot_{M(\vec{x})}$, respectively.
\item If $\phi$ is $\psi\square\chi$, where $\square$ is one of $\wedge$, $\vee$, or $\To$, then $\widehat M_{\vec{x}}(\phi)$ must be, respectively, a product, coproduct, or exponential object $\widehat M_{\vec{x}}(\psi)\square \widehat M_{\vec{x}}(\chi)$ of the objects $\widehat M_{\vec{x}}(\psi)$ and $\widehat M_{\vec{x}}(\chi)$ in $\fib{C}^{M(\vec{x})}$.\footnote{It is probably better to say, here and in the remaining clauses, that $M_{\vec{x}}\phi$ is a certain \emph{diagram} in $\C$ with a specified object (lying in $\fib{C}^{M(\vec{x})}$), so that, for instance $M_{\vec{x}}(\psi\wedge\chi)$ is not only \emph{a} product $M_{\vec{x}}(\psi)\wedge M_{\vec{x}(\psi)}$, but in fact remembers'' \emph{how} it is a product; however, this will not make a difference for us.}
\item\label{item:sigma-interp-over-2LAST-quant} If $\phi$ is $\forall z\psi$ or $\exists z\psi$, then $\widehat M_{\vec{x}}(\phi)$ (where we can assume, by our convention on bound variables, that $z\notin\vec{x}$) must be, respectively, an indexed product $\prod_{\pi}\widehat M_{\vec{x}z}(\psi)$ or sum $\sum_{\pi}\widehat M_{\vec{x}z}(\psi)$ of $\widehat M_{\vec{x}z}(\psi)$, where $\pi=\langle{\pi_1^M,\cdots,\pi_{\len{\vec{x}}}^M}\rangle$.
\item\label{item:sigma-interp-over-eq} If $\phi$ is $s=t$, with $s$ and $t$ terms of sort $B$, then $\widehat M_{\vec{x}}(\phi)$ must be a pullback $\br{M_{\vec{x}}(s),M_{\vec{s}}(t)}^*\Eq_{MB}$ along $\br{M_{\vec{x}}(s),M_{\vec{s}}(t)}\colon M(\vec{x})\to MB\times MB$ of an equality object $\Eq_{MB}$ over $MB\times MB$.
\end{enumerate}

We also say that the pair $(\widehat M,M)$ is an \emph{interpretation of $\sigma$ in} $\fib{C}$ and write $(\widehat M,M)\colon\sigma\to\fib{C}$.

If $\phi$ is a closed formula (i.e., has no free variables), we say that $M$ \emph{satisfies} $\phi$, and write $M\vDash\phi$, if there exists a morphism $\tm\to\widehat M_{\emptyset}(\phi)$ in $\fib{C}^{M(\emptyset)}=\fib{C}^{\tm_\B}$ (this notion really depends on $\fib{C}$ and not just on $M$, but by Proposition~\ref{prop:sigma-interp-over-exist-unique} below it does not depend on $\widehat M$).

Given an interpretation $\widehat M$ in $\fib{C}$ over $M$ and a morphism of $h^=$-fibrations $(\Phi,\phi)\colon\fib{C}\to\fibr{C'}{C'}{B'}$, we obtain an interpretation in $\fib{C'}$ over $\phi\circ M$, denoted $\Phi\circ \widehat M$, given by $(\Phi\circ\widehat M)_{\vec x}(\psi):=\Phi(\widehat M_{\vec x}(\psi))$.

\subsubsection{}\prop\label{prop:sigma-interp-over-exist-unique}
Given an $h^=$-fibration $\fibr{C}{C}B$ and an interpretation $\sigma\colon M\to\C$ in $\C$, any two interpretations in $\fib{C}$ over $M$ are isomorphic -- i.e., if $\widehat M$ and $\widehat M'$ are two interpretations in $\fib{C}$ over $M$, then $\widehat M_{\vec{x}}\phi$ and $\widehat M'_{\vec{x}}\phi$ are isomorphic objects in $\fib{C}^{M(\vec{x})}$ for every $(\phi,\vec{x})$. Moreover (using the axiom of choice), there always exists a interpretation in $\fib{C}$ over $M$.
\pf
The first claim follows by induction on $\phi$ since for each possible $\phi$, one of the clauses in the definition of interpretation over $M$ determines $\widehat M_{\vec{x}}(\phi)$ up to isomorphism based on the constituent formulas of $\phi$.

For the second claim, we just choose fiberwise product, fiberwise coproducts, etc., in $\fib{C}$ and then recursively define an interpretation over $M$.
\qed

\subsubsection{}
We now consider some examples of interpretations in $h^=$-fibrations, all of which were mentioned in \S\ref{subsec:propositions-as-spaces} in the introduction.

If $\C=\Set$, then an interpretation $M\colon\sigma\to\Set$ is just a $\sigma$-structure in the usual, set-theoretic sense.\footnote{The one difference, perhaps, is that the ``product'' operation on sets is usually fixed, whereas in this definition, the interpretation $M$ ``chooses'' the relevant products.} Then, if we consider the sub-$h^=$-fibration $\fib{P(\Set)}\subset\fib{F(\Set)}$ of $\fib{F(\Set)}$ consisting of morphisms which are inclusions of subsets -- then there is exactly one interpretation in $\fib{P(\Set)}$ over $M$, and this recovers the classical semantics of first-order logic. Namely, $M_{\vec{x}}(\phi)$ is then given by $\set{\vec{a}\in M(\vec{x}) \mid M\underset{\vec{x}\mapsto\vec{a}}{\vDash}\phi}$.

Next, if we instead take the fibration $\fib{F(\Set)}$, then an interpretation in $\fib{F(\Set)}$ over $M$ recovers the ``propositions as sets'' semantics.

Similarly, taking $\fib{F(\C)}$ for any locally cartesian closed category with finite coproducts $\C$ with, we recover the ``propositions as objects of $\C$'' semantics.

Finally, if we take $\C=\sSet$, we get the ``propositions as objects of $\sSet$'' semantics which, as explained in \S\ref{subsec:propositions-as-spaces}, is \emph{not} our desired homotopical semantics.

We now want to convince ourselves that taking the fibration $\fib{HoF_{\fb}(\Kan)}$ \emph{does} give us the homotopical semantics). But let us first define what the latter is.

\subsubsection{}\defn\label{defn:homotopical-interpretation}
Given an interpretation $M\colon\sigma\to\C$, with $\C$ any suitable model category (in fact, any model category which is locally bicartesian closed), we define \emph{a homotopical interpretation $\widehat M$ over $M$} to be a function taking each $\sigma$-formula-in-context $(\phi,\vec{x})$ to an object $\widehat M_{\vec{x}}(\phi)\in\fib{F(\C)}^{M(\vec{x})}=\C/M(\vec{x})$, satisfying the same conditions as in Definition~\ref{defn:interpretation-over}, except that in \ref{item:sigma-interp-over-eq} the equality object $\Eq_{MB}$ over $MB$ is replaced by a \emph{path object} over $MB$ -- i.e, the second factor in a factorization of $\Delta_{MB}\colon MB\to MB\times MB$ as a weak equivalence followed by a fibration.

Below, we will be concerned only with interpretations $M\colon\sigma\to\C_\fb$ landing in the subcategory $\C_\fb$. In this case, each object $\widehat M_{\vec{x}}(\phi)$ lies in the subfibration $\fib{F_\fb(\C_\fb)}$ of $\fib{F(\C)}$ (this follows by induction using Proposition~\ref{prop:suitable-mfib-hfib} and the fact that $\widehat M_{\vec{x}}(s=t)$ is by definition a fibration).

Again, for a closed formula $\phi$, we write $M\vDash\phi$ to mean that there exists a morphism $\tm_{\C}\to\widehat M_{\emptyset}\varphi$ in $\C/M(\vec{x})\cong\C$ and again (now by Proposition~\ref{prop:htpical-interp-unique}), we note that this is independent of $\widehat M$.

We again have, as in Proposition~\ref{prop:sigma-interp-over-exist-unique}, that a homotopical interpretation over $M$ always exists. However, we no longer have the uniqueness up to isomorphism, the reason being that there exist non-isomorphic path spaces over a given space. Still (for suitable $\C$), it is unique up to \emph{homotopy equivalence} as we will see presently.

\subsubsection{}\defn
Let $\C$ be a model category and $A,B\in\C_\cfb$ be cofibrant-fibrant objects. Recall that in this case, a morphism $f\colon A\to B$ is a weak equivalence if and only if there exists $g\colon B\to A$ such that $gf$ and $fg$ are homotopic to the identity (where we recall that for objects in $\C_\cfb$, ``left-homotopic'' and ``right-homotopic'' are equivalent) -- see, e.g., \cite[Proposition~10.7]{sentai}. For emphasis, we call such an $f$ a \emph{homotopy equivalence} and $g$ a \emph{homotopy inverse}.

Note that $f$ is homotopy equivalence if and only if it is an equivalence (see Definition~\ref{defn:psnt}) with respect to the Quillen 2-categorical structure on $\C_\cfb$ \cite[Definition~15.7]{sentai}.

If such an $f$ exists, we say that $A$ and $B$ are \emph{homotopy equivalent}.

We note that in practice, we will deploy this notion for suitable model categories $\C$, in which every object is cofibrant, so it is only necessary to restrict to fibrant objects.

\subsubsection{}\prop\label{prop:htpical-interp-unique}
Let $M\colon\sigma\to\C_\fb$ be an interpretation with $\C$ a suitable model category.

Now suppose $\widehat M$ is an \emph{arbitrary} function assigning each $\sigma$-formula-in-context $(\phi,\vec{x})$ an object $\widehat M_{\vec{x}}(\phi)\in\C/M(\vec x)_\fb$. Recall the localization morphism $\gamma\colon\fib{F_\fb(\C_\fb)}\to\fib{HoF_\fb(\C_\fb)}$, and define the function $\gamma\circ\widehat M$ (taking values in the fibers of $\fib{HoF_\fb(\C_\fb)}$) as in the end of Definition~\ref{defn:interpretation-over}.

Then $\gamma\circ\widehat M$ is an interpretation in $\fib{HoF_\fb(\C_\fb)}$ over $M$ if and only if $\widehat M$ is homotopy equivalent to a homotopical interpretation over $M$, in the sense that there is some homotopical interpretation $\widehat M'$ over $M$ such that $\widehat M_{\vec x}\phi$ and $\widehat M'_{\vec{x}}\phi$ are homotopy equivalent for each $(\phi,\vec{x})$.

Moreover, any two homotopical interpretations over $M$ are homotopy equivalent.
\pf
Suppose $\widehat M$ is homotopy equivalent to a homotopical interpretation $\widehat M'$. Then $\gamma\circ\widehat M$ and $\gamma\circ\widehat M'$ are isomorphic, and hence it suffices to check that the latter is an interpretation in $\fib{HoF_\fb(\C_\fb)}$ over $M$.

We need to check each of the conditions in Definition~\ref{defn:interpretation-over}. Conditions \ref{item:sigma-interp-over-FIRST-topbot}-\ref{item:sigma-interp-over-2LAST-quant} follow from the fact that they are satisfied by $\widehat M$ and the fact that $\gamma\colon\fib{F_\fb(C_\fb)}\to\fib{HoF_\fb(C_\fb)}$ is a morphism of $h$-fibrations. Condition~\ref{item:sigma-interp-over-eq} follows from \cite[Proposition~13.1]{sentai} and the definition of homotopical interpretation.

Conversely, suppose $\gamma\circ\widehat M$ is an interpretation in $\fib{HoF_\fb(\C_\fb)}$ over $M$, and let $\widehat M'$ be any homotopical interpretation. Then by what we just showed, $\gamma\circ\widehat M'$ is also an interpretation in $\fib{HoF_\fb(\C_\fb)}$ over $M$, and hence isomorphic to $\gamma\circ\widehat M$ by Proposition~\ref{prop:sigma-interp-over-exist-unique}. But this means precisely that $\widehat M$ and $\widehat M'$ are homotopy equivalent (since, by the suitability of $\C$, all the objects $\widehat M_{\vec{x}}(\phi)$ and $\widehat M'_{\vec{x}}(\phi)$ are cofibrant-fibrant).

The last claim follows for the same reason.
\qed

\subsection{Homotopy equivalence}\label{subsec:htpy-equiv}
We now introduce the remaining notions that are present in the statement of the special invariance theorem, Theorem~\ref{thm:special-invariance} (which was also stated in the \hyperref[thm:special-invariance-in-intro]{introduction}).

The first is the definition of homotopy equivalence of $\sigma$ interpretations. We also include the definition of (ordinary) homomorphism of $\sigma$-interpretations, for the sake of motivation, and since we will use it later (in the \hyperref[defn:free-fp-cat]{definition of free finite product category}).

The second is, given two fibrations $X\to A$ and $Y\to B$ in a model category, the notion of a homotopy between morphisms $X\to Y$ \emph{lying over} a given homotopy between morphisms $A\to B$. In the case of topological spaces, this is a familiar notion, and we verify that it is well-behaved in general.

\subsubsection{}\label{defn:homomorphism}\defn
Given two interpretations $M,N\colon\sigma\to\C$, a \emph{homomorphism} $\alpha\colon{}M\to{}N$ consists of morphisms $\alpha_A\colon{}MA\to{}NA$ for each sort $A\in\Ob\sigma$, such that for each $f\in\sigma(\vec{A},B)$, the following diagram commutes, where we write $\alpha_{\vec{A}}$ for $\alpha_{A_1}\times\cdots\times\alpha_{A_{\len{\vec{A}}}}$.
\[
  \begin{tikzcd}[column sep=60pt]
    M\vec{A}\ar[r, "\alpha_{\vec{A}}"]\ar[d, "Mf"']&
    N\vec{A}\ar[d, "Nf"]\\
    MB\ar[r, "\alpha_{B}"]&NB
  \end{tikzcd}
\]
Given a suitable model structure on $\C$ and two interpretations $M,N\colon\sigma\to\C_\fb=\C_{\cfb}$, a \emph{homotopy homomorphism} $\alpha\colon{}M\to{}N$ consists of a morphism $\alpha_A\colon{}MA\to{}NA$ for each $A\in\Ob\sigma$, such that, for each function symbol $f\in\sigma(\vec{A},B)$, the morphisms $\alpha_B\circ{}Mf$ and $Nf\circ\alpha_{\vec{A}}$ are homotopic (note that $M\vec{A}$ is still in $\C_\fb=\C_{\cfb}$ as the latter is close under products).

A homotopy homomorphism $\alpha\colon{}M\to{}N$ is a \emph{homotopy equivalence} if $\alpha_A$ is a homotopy equivalence for each $A\in\Ob\sigma$. If such an $\alpha$ exists, we say that $M$ an $N$ are \emph{homotopy equivalent}.

\subsubsection{}\defn\label{defn:homotopy-over}
Let $\C$ be a model category, let $(X,A,x),(Y,B,y)\in\C^\to$ be fibrations, and let $(p,f),(q,g)\colon{}(X,A,x)\to(Y,B,y)$ be morphisms. A \emph{right homotopy} from $(p,f)$ to $(q,g)$ consists of right-homotopies $k\colon{}A\to B^I$ from $f$ to $g$ and $\hat{k}\colon{}X\to Y^I$ from $p$ to $q$, together with a morphism $y^I\colon Y^I\to B^I$ making the following two diagrams on the right commute.\footnote{Of course, one could give this definition for general morphisms and not just fibrations, but we will not need the general notion, and it would complicate Proposition~\ref{prop:homotopy-over-left-right} below.}
\begin{equation}\label{eq:homotopy-over-diag}
  \begin{tikzcd}
    X\times I\ar[r, "\hat{h}"]\ar[d, "x\times I"']&[-8pt]Y\ar[d, "y"]\\
    A\times I\ar[r, "h"]&B
  \end{tikzcd}
  \hspace{2pt}
  \begin{tikzcd}
    X+X\ar[r, "\br{\partial_1,\partial_2}"]\ar[d, "x"]&X\times I\ar[r, "\sigma"]\ar[d, "x\times I"]
    &[-10pt]X\ar[d, "x"]\\
    A+A\ar[r, "\br{\partial_1,\partial_2}"]&A\times I\ar[r, "\sigma"]&A
  \end{tikzcd}
  \hspace{15pt}
  \begin{tikzcd}
    X\ar[r, "\hat{k}"]\ar[d, "x"']&[-5pt]Y^I\ar[d, "y^I"]\\
    A\ar[r, "k"]&B^I
  \end{tikzcd}
  \hspace{2pt}
  \begin{tikzcd}
    Y\ar[r, "s"]\ar[d, "y"]&[-8pt]Y^I\ar[r, "\br{d_1,d_2}"]\ar[d, "y^I"]
    &[5pt]Y\times Y\ar[d, "y\times y"]\\
    B\ar[r, "s"]&B^I\ar[r, "\br{d_1,d_2}"]&B\times B
  \end{tikzcd}
\end{equation}
We also call $\hat k$ a right homotopy from $p$ to $q$ \emph{over k}. If such a homotopy exists, we say that $(p,f)$ and $(q,g)$ are \emph{right homotopic} and that $p$ and $q$ are \emph{right homotopic over $k$}.

There is a dual notion of a left homotopy $\hat h$ from $p$ to $q$ over a left homotopy $h$ from $f$ to $g$, corresponding to the two diagrams above to the left.

Note that if $f=g$ and $h$ is the trivial homotopy $A\times I\tox{\sigma}A\tox{f=g}B$, then $p$ and $q$ are left-homotopic over $h$ if and only if $(p,f)$ to $(q,g)$ are left-homotopic in the model structure of \cite[Definition~11.1]{sentai} by \cite[Proposition~11.3]{sentai}.

\subsubsection{}\prop\label{prop:homotopy-over-left-right}
With $(p,f)$ and $(q,g)$ as in Definition~\ref{defn:homotopy-over}, if $X$ and $A$ are cofibrant and $Y$ and $B$ are fibrant, then the notions of $(p,f)$ and $(q,g)$ being left or right homotopic are equivalent, and are both equivalence relations.

Moreover, being homotopic over some fixed (left or right) homotopy $h$ from $f$ to $g$ depends only on the homotopy class of $h$ (in the sense of \cite[\S14]{sentai}), and if $p$ and $q$ are homotopic over $h$, and $q$ and $r$ are homotopic over $h'$, then $p$ and $r$ are homotopic over the composite of $h$ and $h'$.
\pf
The first claim follows from the existence of a model structure on $\C^\to$ in which the notions of left and right homotopic are the ones given above, and in which an object of $\C^\to$ is cofibrant iff both domain and codomain are, and is fibrant iff it is a fibration with fibrant codomain (and hence domain).

The model structure in question is characterized by the property that its cofibrations and weak equivalences are the object-wise cofibrations and weak equivalences; see, e.g., \cite[Theorem~5.1.3]{hovey} (this is also known as the ``injective model structure'', and is also a ``Reedy model structure'').

By inspecting the definitions, one finds that a left-homotopy is indeed precisely the notion given above. However, for right-homotopies, there is the additional condition that the induced morphism $Y^I\to B^I\times_{B\times B}(Y\times Y)$ must be a fibration in order for (\ref{eq:homotopy-over-diag}) to be a path object for $(Y,B,y)$. But this can always be arranged by factoring it as a trivial cofibration followed by a fibration $Y^I\to Y^{I'}\to B^I\times_{B\times B}(Y\times Y)$ and then replacing $Y^I$ by $Y^{I'}$.

The last statement follows from the fact that the codomain functor $\C^\to\to\C$ preserves ``all'' of the model structure -- fibrations, cofibrations, weak equivalences, and finite limits and colimits -- and hence takes composites of homotopies to composites of homotopies.

For the second to last claim, we need to show that, given a homotopy $\hat h\colon X\times I\to Y$ from $p$ to $q$ over a homotopy $h\colon A\times I\to B$ (where we denote the auxiliary morphism $X\times I\to A\times I$ by $x\times I$ as above), and a homotopic homotopy $h'\colon A\times I'\to B$, there is a homotopic homotopy $\hat h'\colon X\times I'\to Y$ over $h'$.

So suppose we are given a left-homotopy $H\colon A\times J\to B$ from $h$ to $h'$, so that we have a factorization of $[\sigma,\sigma]\colon(A\times I)+_{A+A}(A\times I')\to A$ as a cofibration $(A\times I)+_{A+A}(A\times I')\tox{[\partial_1,\partial_2]}A\times J$ and a weak equivalence $\sigma\colon A\times J\to A$ (which we may assume is a trivial fibration since $B$ is fibrant).\footnote{Our notation and conventions for pushouts are analogous to those for pullbacks and coproducts, \cite[\S\S5.1,10.2]{sentai}.}

Now let $X\times I'$ be a cylinder object on $X$ with $[\partial_1,\partial_2]\colon X+X\to X\times I'$ a cofibration, and let $x\times I'\colon X\times I'\to A\times I'$ be a morphism as in (\ref{eq:homotopy-over-diag}) (this can be obtained as a diagonal filler assuming $\sigma\colon A\times I'\to A$ is a trivial fibration, which we may by \cite[Proposition~14.4]{sentai}), and factor $[\sigma,\sigma]\colon(X\times I)+_{X+X}(X\times I')\to X$ as a cofibration and trivial fibration $(X\times I)+_{X+X}(X\times I')\tox{[\partial_1,\partial_2]}X\times J\tox{\sigma}X$. Now choose a diagonal filler $x\times J$ in
\[
  \begin{tikzcd}
    X\times I+_{X+X}X\times I'\ar[r, "{[\partial_1(x\times I),\partial_2(x\times I')]}"]
    \ar[d, "{[\partial_1,\partial_2]}"']&[50pt]
    A\times J\ar[d, "\sigma"]\\
    X\times J\ar[r, "x\sigma"]\ar[ru, "x\times J", dashed]&A.
  \end{tikzcd}
\]
Finally, choose a diagonal filler $\hat H$ in the following square (where the morphism on the left is a trivial cofibration).
\[
  \begin{tikzcd}
    X\times I\ar[r, "\hat h"]\ar[d, "\partial_1"']&[35pt]Y\ar[d, "y"]\\
    X\times J\ar[r, "H(x\times J)"']\ar[ru, dashed, "\hat H"]\ar[r]&B
  \end{tikzcd}
\]
The composite $\hat H\partial_2\colon X\times I'\to Y$ is then the desired homotopy.
\qed

\subsubsection{}\defn\label{defn:htpy-equiv-over}
Let $\C$ be a model category and $(p,f)\colon(X,A)\to(Y,B)$ a morphism in $(\C_\cfb)^\to$, and suppose $f$ is a homotopy equivalence. We say that $p$ is a \emph{homotopy equivalence over $f$} if $p$ is a homotopy equivalence, and moreover, there exist homotopy inverses $g\colon B\to A$ and $q\colon Y\to X$ such that the homotopies from $qp$ to $\id_X$ and from $pq$ to $\id_Y$ lie over homotopies from $gf$ to $ \id_A$ and from $fg$ to $\id_B$.

\subsubsection{}
With these definition in place, we can now formulate the special invariance theorem, Theorem~\ref{thm:special-invariance}, which was already restated in the \hyperref[thm:special-invariance-in-intro]{introduction}, and so we won't repeat it here.

\subsection{Topological spaces}\label{subsec:top-spaces}
In \S\ref{subsec:h-fib-of-spaces}, we showed that $\fib{HoF_\fb(\Kan)}$ is a $h^=$-fibration and that $\gamma\colon\fib{F_\fb(\Kan)}\to\fib{HoF_\fb(\Kan)}$ is a morphism of $h$-fibrations.

As mentioned in the introduction, we cannot hope to prove this with $\Top$ ($=\Top_\fb$) in place of $\Kan$, since $\Top$ is not locally cartesian closed.

However, it would also be nice to be able to define homotopical semantics valued in topological spaces. We now briefly discuss some possible solutions to this issue.

\subsubsection{}\sectitle{Composing with $\Sing\colon\Top\to\sSet$}
The simplest way to define $\Top$-valued homotopical semantics is to simply do so via simplicial sets (which is the approach we take in \S\ref{subsec:examples}): we have the singular simplicial set functor $\Sing\colon\Top\to\sSet$, which preserves products, being a right adjoint.

Hence, given an interpretation $M\colon\sigma\to\Top$, we can compose it with $\Sing$ to obtain a $\sigma$-interpretation $\Sing\circ M$ in $\sSet$ and we then define a homotopical interpretation over $\sigma$ to be one over $\Sing\circ M$.

The reason why this is a sensible thing to do is that -- at least when $M$ is valued in spaces homotopy equivalent to a CW-complex -- the properties which can be described using the homotopical semantics (namely, the ``homotopical'' ones) can all be recovered from the associated singular simplicial set -- see \S\ref{subsec:examples} for examples.

\subsubsection{}\sectitle{The fibration $\fib{HoF(\Top)}$}
An important observation in this connection is that, although $\fib{F_\fb(\Top)}$ cannot be a $h$-fibration, it turns out that $\fib{HoF_\fb(\Top)}$ (or the equivalent fibration $\fib{HoF(\Top)}$) \emph{is}, after all, a $h^=$-fibration, by virtue of its close relationship to $\fib{HoF_\fb(\Kan)}$.

The reason for this is that, when $\Top$ is endowed with its standard (``Quillen'') model structure (or better, with its ``mixed'' model structure -- see \cite[p.~356]{may-ponto}), the adjunction
$
\begin{tikzcd}
  \Top
  \ar[r, shift left=3pt, "_{\Sing}"]
  \ar[r, phantom, "_{_\top}", pos=0.5]&
  \sSet
  \ar[l, shift left=3pt]
\end{tikzcd}
$
(the left adjoint being ``geometric realization'') becomes a Quillen equivalence, from which it follows that, for fibrant $X\in\Top$, the induced functor $\Top/X\to\sSet/\Sing(X)$ is also a Quillen equivalence (with respect to the induced model structures), which therefore induces equivalences $\Ho(\Top/X)\to\Ho(\sSet/\Sing(X))$.

It follows that we have a morphism of fibrations $\fib{HoF(\Top)}\to\fib{HoF(\sSet)}$ which is a ``fiberwise equivalence'' (Definition~\ref{defn:fiberwise-equivalence}), from which it follows that $\fib{HoF(\Top)}$ is a $h^=$-fibration (since $\fib{HoF(\Kan)}$ is and $\Sing(X)\in\Kan$ for all $X\in\Top$).

Thus, since homotopical interpretations in a model category $\C$ correspond to ordinary interpretations in $\fib{HoF_\fb(\C_\fb)}$, we \emph{can}, in a sense, define homotopical semantics directly in $\Top$, by simply defining this to be the fibrational semantics in $\fib{HoF(\Top)}$.

The main problem with this is that, unlike in the case of suitable $\C$, the $h^=$-fibration structure on $\fib{HoF(\Top)}$ is given very inexplicitly, and hence we cannot say, for any particular formula, what its interpretation in $\Top$ actually is.

However, it is worth noting that we can get a handle on \emph{some} formulas, namely those not involving implication and universal quantification, as $\fib{F_\fb(\Top)}$ \emph{is} a ``$\exists\vee\!\wedgeq$-fibration'' (and $\fib{F_\fb(\Top)}\to\fib{HoF_\fb(\Top)}$ a morphism of such).

\subsubsection{}\sectitle{Convenient categories of spaces}
One might hope to circumvent the failure of $\Top$ to be locally cartesian by passing to a modified category of spaces; the well-established solution to $\Top$ not being cartesian closed is to pass to a ``convenient'' subcategory of $\Top$ (such as the compactly generated spaces) that is cartesian closed. Unfortunately, the various categories used for this purpose are not (or at least not known to be) \emph{locally} cartesian closed.

However, there are well-known ``enlargements'' of subcategories of $\Top$ which \emph{are} locally cartesian closed, i.e., a (reasonably large) full subcategory $\C\subset\Top$ with a full and faithful functor $i\colon\C\to\D$. An example is the category $\D$ of subsequential spaces \cite{johnstone-topological-topos} (with $\C$ the category of sequential topological spaces).

Hence, one might hope to use such a category $\D$ for the homotopical semantics instead, by extending to it the usual model structure(s) on $\Top$. One can indeed produce such model structures, at least in the case of subsequential spaces. However, the result still fails to be suitable, and as a result suffer from the same problem mentioned in the previous section: it becomes difficult to actually compute the interpretation of any given formula (in particular, it requires taking cofibrant replacements).

\section{The invariance theorem}\label{sec:invariance}
We now formulate and prove the homotopy-invariance property of the homotopical semantics; we first prove the ``abstract'' version, which is purely categorical, and then deduce from it the ``special'' version of the theorem, concerning the homotopical semantics. Let us discuss these is order.

The formulation and proof of the abstract invariance theorem is actually independent of the results of the previous section; these were presented first only in order to give some motivation for the abstract invariance theorem. The latter involves an arbitrary $h^=$-fibration $\fib{C'}$, which for the special invariance theorem we will want to take to be the $h^=$-fibration $\fib{HoF_\fb(\Kan)}$ introduced above.

Moreover, $\fibr{C'}{C'}{B'}$ is assumed to have the structure of a 1-discrete 2-fibration (1D2F), which is a 2-categorical structure on $\C'$ and $\B'$ compatible with the fibration in which the fibers are however still 1-categories (this is the ``1-discreteness'') -- see \cite[\S6]{sentai} for the definition and some motivation. However, the presence of a 1D2F structure on $\fib{C'}$ is not really an extra requirement: in \emph{op. cit.}, it is proven that every $\wedgeq$-fibration admits a ``universal'' 1D2F structure, and that, in the case of $\fib{HoF_\fb(\C_\fb)}$ (for any model category $\C$), this recovers the familiar 2-category structure on $\C_\cfb$ given in terms of homotopies (and we will shown in \S\ref{subsec:special-invariance} that the 2-categorical structure on the total category of $\fib{HoF_\fb(\C_\fb)}$ also has a natural interpretation in terms of homotopies).

The abstract invariance theorem also involves a second fibration $\fibr{C}CB$ which is \emph{free} on a signature $\sigma$ in the sense discussed in the introduction, so that the morphisms $\fib{C}\to\fib{C'}$ correspond precisely to interpretations of $\sigma$ in $\C'$ (the morphism $\C\to\C'$ then capturing the unique-up-to-isomorphism interpretation of $\sigma$ in $\fib{C'}$ over $\B'$); moreover, the \emph{natural isomorphisms} between such morphisms, correspond to isomorphisms of $\sigma$-interpretations $\sigma\to\B$. As discussed in the introduction, this immediately leads to the \emph{isomorphism invariance} property for interpretations in $\fib{C'}$. Note that the 1D2F structure on $\fib{C'}$ plays no role in this.

Now, when the f.p.\ category $\B$ has a 2-categorical structure, there is a notion of \emph{pseudonatural equivalence} of interpretations $\sigma\to\B'$, which correspond to homotopy equivalences of interpretations in the case of $\Kan$. On the other hand, we have a notion of pseudonatural equivalence of morphisms $\fib{C}\to\fib{C'}$, and the existence of one of these between the morphisms $\fib{C}\to\fib{C'}$ induced by two interpretations $\sigma\to\B'$ recovers, in the case of $\fib{HoF_\fb(\Kan)}$ the homotopy invariance property asserted in the special invariance theorem.

The import of the invariance theorem, then, is that the free fibration $\fib{C}$ satisfies with respect to pseudonatural equivalence the same universal property that it satisfies with respect to isomorphisms, i.e., that a pseudonatural equivalence between two interpretations $\sigma\to\B'$ induces a pseudonatural equivalence of the induced morphisms $\fib{C}\to\fib{C'}$. This comprises two parts.

In the first part, we have $\B$ being an \emph{arbitrary} f.p.\ category, and $\fib{C}$ a free $h^=$-fibration \emph{over} $\B$, so that an f.p.\ morphism $\B\to\B'$ induces a morphism $\fib{C}\to\fib{C'}$. The statement is then that a pseudonatural equivalence between two morphisms $\B\to\B'$ induces a pseudonatural equivalence between the resulting morphisms $\fib{C}\to\fib{C'}$.

The second part only involves f.p.\ (2-)categories and makes no reference to fibrations. Namely, here we take $\B$ to be free on $\sigma$, so that each interpretation $\sigma\to\B'$ induces an f.p.\ morphism $\B\to\B'$, and the statement is then that a pseudonatural equivalence between interpretations $\sigma\to\B'$ induces a pseudonatural equivalence between the induced morphisms $\B\to\B'$.

With the abstract invariance theorem in place, we then, in \S\ref{subsec:special-invariance}, prove the ``special'' or ``syntactic'' invariance theorem which, unlike the abstract invariance theorem, involves the interpretation of terms and formulas in a $h^=$-fibration, as introduced in \S\ref{subsec:interpretations}. This is deduced easily from the abstract invariance theorem by using the equivalence between the 2-categorical notions present in the abstract invariance theorem and the corresponding ``homotopical'' notions in the particular 1D2F $\fib{HoF_\fb(\Kan)}$. Most of this equivalence was proven in or is deduced easily from the work in \cite{sentai}; however, there is one essential part (namely, the homotopical interpretation of 2-cells in the total category of $\fib{HoF_\fb(\Kan)}$) that is missing in \emph{op. cit.}, and which is proven in \S\ref{subsec:special-invariance}.

The special invariance theorem is also proven in two ``stages'' (though these do not correspond to the two parts of the abstract invariance theorem): the first is to draw a general ``homotopy-invariance'' (or better, ``pseudonatural-equivalence-invariance'') conclusion about the interpretation of first-order logic in an arbitrary $h^=$-fibration $\fib{C'}$, and the second is to specialize this to the case of $\fib{HoF_\fb(\Kan)}$.

We begin with some preliminaries on 2-categories, and we prove some facts about 1D2Fs that were not considered in \cite{sentai} (but otherwise, we refer to the latter for general background on 1D2Fs).

\subsection{Preliminaries on 2-categories and 1D2Fs}\label{subsec:2-cat-prelims}
We begin with some generalities on 2-categories and 1D2Fs.

First, we will recall the notion of pseudonatural transformation and equivalence, and some related notions. We note that this is a small part of a larger structure; the totality of 2-categories forms a \emph{3-dimensional structure}, which is more or less elaborate depending on the weakness or strictness of the concepts being considered. Our approach here is to strive for the minimal possible generality that suffices for our purposes.

We then prove some basic facts that we will need about 1D2Fs. Each of them is a generalization of a basic property of fibrations, in which equality of morphisms is replaced by existence of a 2-cell.

Finally, we give the notion of pseudonatural equivalence of morphisms of fibrations, which will be essential for the formulation of the abstract invariance theorem.

\subsubsection{}\label{defn:psnt}\defn
Given a category $\C$, a 2-category $\D$ and functors $F,G\colon\C\to\D$, a \emph{pseudonatural transformation} $\alpha\colon{}F\to{}G$ consists of the following data (i)-(ii), subject to the conditions (iii)-(iv):
\begin{enumerate}[(i),itemsep=0pt,topsep=-3pt]
\item a 1-cell $\alpha_A\colon{}FA\to{}GA$ for each $A\in\C$
\item an isomorphism 2-cell
\[
    \begin{tikzcd}[column sep=60pt]
      FA
      \ar[r, "Ff"]\ar[d, "\alpha_A"']&
      FB
      \ar[d, "\alpha_B"]
      \ar[ld, Rightarrow, shorten <=5pt, shorten >=5pt, "\alpha_{f}", pos=0.4, "\sim" {sloped, pos=0.35, inner sep=3pt},
      start anchor={[xshift=-10pt]}, end anchor={[xshift=10pt]}]
      \\
      GA\ar[r, "Gf"']
      &
      GB
    \end{tikzcd}
\]
for each morphism $f\colon{}A\to{}B$ in $\C$
\item For each pair $A\tox{f}B\tox{g}C$ of composable morphisms in $\C$, we have $(\id_{Gg}\circ\alpha_f)(\alpha_g\circ\id_{Ff})=\alpha_{gf}$.
\item For each $A\in\C$, we have $\alpha_{\id_A}=\id_{\alpha_A}$.
\end{enumerate}

Note that a pseudonatural transformation $\alpha$ with $\alpha_f=\id$ for all $f$ is just a natural transformation.

Given pseudonatural transformations $F\tox{\alpha}G\tox{\beta}H$, their \emph{composite} $\beta\circ\alpha$ is defined by the prescriptions $(\beta\circ\alpha)_A=\beta_A\circ\alpha_A$ and $(\beta\circ\alpha)_f=(\beta_f\circ\alpha_A)(\beta_A\circ\alpha_f)$. We leave it to the reader to verify that this is again a pseudonatural transformation.

A 1-cell $f\colon{}A\to{}B$ in a 2-category is an \emph{equivalence} if there exists a 1-cell $g\colon{}B\to{}A$ and isomorphism 2-cells $g\circ{}f\cong{}\id_A$ and $f\circ{}g\cong{}\id_B$ -- such a $g$ is called a \emph{quasi-inverse} to $f$. The pseudonatural transformation $\alpha\colon{}F\to{}G$ is a \emph{pseudonatural equivalence} if $\alpha_A$ is an equivalence in $\D$ for each $A\in\C$.

\subsubsection{}\label{defn:psnt-whisker}\defn
Given categories $\B$ and $\C$, 2-categories $\D$ and $\E$, functors $F\colon\B\to\C$ and $G,H\colon\C\to\D$, a 2-functor $K\colon\D\to\E$ and a pseudonatural transformation $\alpha\colon{}G\to{}H$, as in
\[
  \begin{tikzcd}[row sep=0pt, column sep=20pt]
    &&
    \ar[dd, shorten >=0pt, shorten <=-5pt, Rightarrow, "\alpha"' inner sep=5pt, pos=0.3]
    &
    &\\[-5pt]
    \B\ar[r, "F"]&
    \C\ar[rr, "G", bend left]\ar[rr, "H"', bend right]
    &&
    \D\ar[r, "K"]&\E,\\
    &&{}&&{}&
  \end{tikzcd}
\]
we define (i) the \emph{whiskering of $\alpha$ by $K$}, which we denote by $K\circ\alpha$, to be the pseudonatural transformation $KG\to{}KH$ defined by $(K\circ\alpha)_A=K(\alpha_A)$ and $(K\circ\alpha)_f=K(\alpha_f)$ for $A\in\Ob\C$ and $f\in\Ar\C$; and (ii) the \emph{whiskering of $\alpha$ by $F$}, denoted by $\alpha\circ{}F$ to be the pseudonatural transformation $GF\to{}HF$ defined by $(\alpha\circ{}F)_A=\alpha_{FA}$ and $(\alpha\circ{}F)_f=\alpha_{Ff}$.

We leave to the reader the easy proof that these are in fact pseudonatural transformations as claimed.

We note that the whiskering of a pseudonatural equivalence (on either side) is again a pseudonatural equivalence.

\subsubsection{}\label{prop:1d2f-cart-gen}\prop
In a 1D2F, we have the following generalization of the universal property for cartesian morphisms.

Let $\fibr{C}CB$ be a 1D2F, and suppose we have a solid diagram
\[
  \begin{tikzcd}[row sep=5pt]
    &|[alias=Q]|Q\car[rd, "q"]&\\
    P\ar[ru, dashed, "p"]\ar[rr, "r"'{name=rar}]
    \ar[from=rar, to=Q, Rightarrow, dashed, shorten <=2pt, shorten >=1pt]&&
    R\\[5pt]
    &|[alias=B]|B\ar[rd, "g"]&\\
    A\ar[ru, "f"]\ar[rr, "h"'{name=h}]
    \ar[from=h, to=B, Rightarrow, "\alpha"' inner sep=4pt, shorten <=2pt, shorten >=1pt]&&
    C\\
  \end{tikzcd}
\]
with $q$ cartesian over $g$, $r$ lying over $h$, and $\alpha$ a 2-cell $h\to{}gf$.

Then there exists a unique 1-cell $p\colon{}P\to{}Q$ over $f$ for which there exists a (necessarily unique) 2-cell $r\to{}qp$ over $\alpha$.
\pf
We know that there is a unique 2-cell $\sigma\colon{}r\to{}r'$ over $\alpha$ with domain $r$, and we can (and must) then take $p$ to be unique morphism $P\to{}Q$ over $f$ such that $qp=r'$:
\[
  \begin{tikzcd}[row sep=5pt, baseline=(C.base)]
    &Q\car[rd, "q"]&\\[10pt]
    P\ar[ru, "p", dashed]
    \ar[rr, "r"'{name=rar}]
    \ar[rr, "r'"{name=rarr}, bend left]
    \ar[from=rar, to=rarr, Rightarrow, "\sigma"' inner sep=4pt, shorten <=2pt, shorten >=2pt]&&
    R\\[5pt]
    &|[alias=B]|B\ar[rd, "g"]&\\
    A\ar[ru, "f"]\ar[rr, "h"'{name=h}]
    \ar[from=h, to=B, Rightarrow, "\alpha"' inner sep=4pt, shorten <=2pt, shorten >=1pt]&&
    |[alias=C]|C\\
  \end{tikzcd}
  \tag*{\qed}
\]

\subsubsection{}\label{prop:1d2f-cartuniq}\prop
Next, we have a generalization in 1D2Fs of the uniqueness up to isomorphism of cartesian morphisms in fibrations.\vspace{-5pt}

Given a 1D2F $\fibr{C}CB$ and a diagram
\[
  \begin{tikzcd}[row sep=5pt]
    P\car[rd, "p" name=rar]\ar[dd, "r"']&\\
    &R&\\
    |[alias=Q]|Q\car[ru, "q"']&
    \ar[from=rar, to=Q, Rightarrow, shorten <=6pt, shorten >=6pt,
        "\sigma" pos=0.4, "\sim" {sloped, pos=0.4}]
    \\[10pt]
    A
    \ar[r, bend left, "f"{name=f}]
    \ar[r, bend right, "g"'{name=g}]
    &B
    \ar[Rightarrow, from=f, to=g, "\alpha" inner sep= 3pt, shorten <=3pt, shorten >=3pt,
    "\vsim"' {inner sep=3pt, pos=0.4}]
  \end{tikzcd}
\]
in which $r$, $p$, $q$, $\sigma$ lie over $\id_A$, $f$, $g$, $\alpha$, respectively, and $p$ and $q$ are cartesian: if $\alpha$ (and hence $\sigma$) is an isomorphism 2-cell, then $r$ is an isomorphism 1-cell.
\pf
By Proposition~\ref{prop:1d2f-cart-gen}, there is a unique morphism $r'\colon{}Q\to{}P$ over $A$ for which there exists a (necessarily unique) 2-cell $\sigma'\colon{}q\to{}pr'$ over $\alpha\I$:
\[
  \begin{tikzcd}[row sep=16pt, column sep=40pt]
    P\ar[rd, "p" name=rar]\ar[d, "r"']&\\
    |[alias=Q]|Q\ar[r, "q"' {name=q, pos=0.4}]\ar[d, "r'"']&R
    \ar[from=rar, to=Q, Rightarrow, shorten <=6pt, shorten >=6pt,
        "\sigma" pos=0.4, "\sim" {sloped, pos=0.4}]\\
    |[alias=P]|P\ar[ru, "p"']&
    \ar[from=q, to=P, Rightarrow, shorten <=-2pt, shorten >=9pt,
        "\sigma'" {pos=0.0, inner sep=1pt}, "\sim" {sloped, pos=0},
        start anchor={[xshift=-3pt]}]\\
    \\[-9pt]
    A
    \ar[r, bend left=60pt, "f"{name=f}]
    \ar[r, "g"' {description, name=g}]
    \ar[r, bend right=60pt, "f"'{name=ff}]
    &B.
    \ar[Rightarrow, from=f, to=g, "\alpha" inner sep= 3pt, shorten <=4pt, shorten >=2pt,
        "\vsim"' {inner sep=3pt, pos=0.5}]
    \ar[Rightarrow, from=g, to=ff, "\alpha\I" {inner sep= 3pt, pos=0.3}, shorten <=3pt,
        shorten >=3pt, "\vsim"' {inner sep=3pt, pos=0.3}]
  \end{tikzcd}
\]
Then $r'r$ and $\id_P$ are both the unique morphism $t\colon{}P\to{}P$ over $\id_A$ for which there exists a 2-cell $p\to{}pt$; hence $r'r=\id_P$. Similarly, $rr'=\id_Q$.
\qed

\subsubsection{}\label{prop:1d2f-2-of-3}\prop
We also have a generalization in 1D2Fs of the ``2-of-3'' property of cartesian morphisms in fibrations.

Given a 1D2F $\fibr{C}CB$ and a diagram
\[
  \begin{tikzcd}[row sep=5pt]
    &|[alias=Q]|Q\car[rd, "q"]&\\
    P\ar[ru, "p"]\car[rr, "r"'{name=rar}]
    \ar[from=rar, to=Q, Rightarrow, shorten <=4pt, shorten >=1pt, "\sigma"' inner sep=3pt,
        "\vsim" inner sep=3pt]&&
    R\\[5pt]
    &|[alias=B]|B\ar[rd, "g"]&\\
    A\ar[ru, "f"]\ar[rr, "h"'{name=h}]
    \ar[from=h, to=B, Rightarrow, "\alpha"' inner sep=4pt, shorten <=2pt, shorten >=1pt,
        "\vsim" inner sep=3pt]&&
    C\\
  \end{tikzcd}
\]
with $p$, $q$, $r$, $\sigma$ lying over $f$, $g$, $h$, $\alpha$, respectively, if $\alpha$ (and hence $\sigma$) is invertible, and $q$ and $r$ are cartesian, then so is $p$.
\pf
Choosing a cartesian lift $\ct\colon{}f^*Q\to{}Q$ of $f$, we obtain a factorization
$p=\ \ct\cind{p}$ as in
\[
  \begin{tikzcd}
    f^*Q\car[r, "\ct"]&|[alias=Q]|Q\car[rd, "q"]&\\
    P\ar[ru, "p"]\car[rr, "r"'{name=rar}]\ar[u, "\cind{p}"]
    \ar[from=rar, to=Q, Rightarrow, shorten <=4pt, shorten >=1pt, "\sigma"' inner sep=3pt,
        "\vsim" {inner sep=3pt}, start anchor={[xshift=2pt]}]&&
    R.
  \end{tikzcd}
\]
It then follows from Proposition~\ref{prop:1d2f-cartuniq} that $\cind{p}$ is an isomorphism, and hence that $p$ is cartesian.
\qed

\subsubsection{}\label{defn:psnt-over-psnt}\defn
Let $\fibr{C}CB$ be a fibration, $\fibr{C'}{C'}{B'}$ be a 1D2F, $\phi,\psi\colon\B\to{}\B'$ functors, and $\Phi,\Psi\colon\C\to{}\C'$ functors lying over $\phi$ and $\psi$. We say that a pseudonatural transformation $\beta\colon\Phi\to\Psi$ \emph{lies over} a pseudonatural transformation $\alpha\colon\phi\to\psi$ if $\fib{C'}\circ\beta=\alpha\circ\fib{C}$ (here we are using the ``whiskering'' operations from \ref{defn:psnt-whisker}).

\subsection{The abstract invariance theorem: first part}\label{subsec:abstract-invariance-theorem}
We now come to the proof of the (first part of the) abstract version of the homotopy-invariance property. In this generality, the theorem states that a \emph{free} $h^=$-fibration satisfies with respect to \emph{pseudonatural} transformations the same property which it satisfies with respect to natural transformations.

We make a small digression to discuss \emph{freeness}. In classical (``0-categorical'') algebra, a \emph{free} object (group, ring, etc.) is required to satisfy a certain universal property, which then determines it up to isomorphism. In the case of categorical structures, it is often more natural to impose conditions that determine the object under consideration up to \emph{equivalence}. However, there are usually different conditions which do this.

In the case at hand, namely that of $\fibr{C}CB$ being a \emph{free $h^=$-fibration over $\B$}, the weakest such condition one could impose is that, for any other $h^=$-fibration $\fibr{C'}{C'}B$ over $\B$, there exists a morphism $\fib{C}\to\fib{C'}$ of fibrations over $\B$, and that, any two such morphisms are isomorphic; and this is in fact the definition we use.

A slightly stronger condition one could demand is that, given any $h^=$-fibration $\fibr{C'}{C'}{B'}$, any f.p.\ functor $\phi\colon\B\to{}\B'$ can be extended to a morphism $(\Phi,\phi)\colon\fib{C}\to\fib{C'}$ of $h^=$-fibrations, and for any natural isomorphism $\phi\to{}\psi$ to another f.p.\ functor, there is a natural isomorphism $\Phi\to\Psi$ lying over it; let us call this the \emph{global} universal property, and the previous one the \emph{local} universal property.

One can see that a free $h^=$-fibration (in the above sense -- i.e., satisfying the local universal property) automatically satisfies the global universal property. The invariance theorem then says that, when $\fib{C'}$ is a 1D2F, $\fib{C}$ satisfies the global universal property with respect to \emph{pseudonatural equivalences} and not just natural isomorphisms.

The proof of the global universal property from the local one is more or less the same for natural isomorphisms and for pseudonatural equivalences -- indeed, we will prove both of them simultaneously -- and proceeds (roughly) by pulling back the fibration $\fib{C'}$ along the given functors $\phi,\psi\colon\B\to{}\B'$ and then appealing to the universal property of $\fib{C}$ with respect to $h^=$-fibrations over $\B$.

An essential step here is showing that the natural isomorphism or pseudonatural equivalence $\phi\to\psi$ induces an equivalence of fibrations $\phi^*\fib{C'}\cong\psi^*\fib{C'}$ over $\B$ between the pullbacks. It is easy to see why this should work for pseudonatural transformations and not just natural transformations, if we look at the situation from the ``pseudo-functor to $\Cat$'' perspective on fibrations (see, e.g., \cite[\S6]{sentai}). Namely, from this perspective, the pullback $\phi^*\fib{C'}$ of a fibration $\fib{C'}$ along a morphism $\phi\colon\B\to\B'$ is just given by composing $\phi^\op\colon\B^\op\to\B'^\op$ with the pseudo-functor $\fpsf{C'}\colon\B^\op\to\Cat$. On the other hand, a morphism of fibrations over $\B$ corresponds to a pseudonatural transformation of pseudo-functors from $\B^\op$ to $\Cat$. Hence, given a natural isomorphism $\alpha\colon\phi\to\psi$, the induced equivalence $\phi^*\fib{C'}\simeq\psi^*\fib{C'}$ between the pullbacks along $\phi$ and $\psi$ is obtained as the whiskering
\[
  \begin{tikzcd}[row sep=0pt, column sep=15pt]
    &\ar[from=dd, shorten >=-3pt, shorten <=1pt, Rightarrow, "\ \alpha^\op"', pos=0.7]&
    &&\\[-5pt]
    \B^\op\ar[rr, "\phi^\op", bend left]\ar[rr, "\psi^\op"', bend right]&&
    \B'^\op\ar[r, "\fpsf{C'}"]&\Cat.\\
    &{}&&{}&
  \end{tikzcd}
\]
The point is now that this whiskering can be carried out just as well if $\alpha$ is a \emph{pseudo}-natural transformation.

\subsubsection{}\label{defn:pullback-fib}\defn
Given a prefibration $\fibr{C'}{C'}{B'}$ and a functor $F\colon\B\to\B'$, we define the \emph{pullback $\fibr{\mathit{F^*}C}{\mathit{F^*}C}{B}$ of $\fib{C}$ along $F$} to be the usual pullback
\[
  \begin{tikzcd}
    F^*\C'\ar[r, "\crt{F}{\C'}"]\ar[d, "F^*\fib{C'}"']
    \ar[rd, phantom, "\lrcorner", pos=0.1]
    &\C'\ar[d, "\fib{C'}"]\\
    \B\ar[r, "F"]&\B'
  \end{tikzcd}
\]
in the (1-)category of categories, where we write $\crt{F}{\C'}$, as indicated, for the associated functor $F^*\C'\to\C'$. Explicitly, the objects of $F^*\C'$ are pairs $(A,P)$ with $A\in\B$ and $P\in\fib{C'}^{FA}$, and the morphisms are pairs $(f,p)$ with $f\in\Ar\B$ and $p\in\Ar\C'$ lying over $Ff$.

The prefibration $F^*\fib{C'}$ is a fibration if $\fib{C'}$ is, and inherits many properties from $\fib{C'}$ -- in particular, if $\B$ is an f.p.\ category and $F$ is an f.p.\ functor, then if $\fib{C'}$ is an $h^=$-fibration, $F^*\fib{C'}$ is as well. The reason is that $\crt{F}{\C'}$ induces isomorphisms $(F^*\fib{C'})^A\to\fib{C'}^{FA}$ on fibers for each $A\in\B$, and -- if $\fib{C'}$ is a fibration -- a morphism in $F^*\C'$ is cartesian or cocartesian if and only if its image under $\crt{F}\C'$ is, and similarly the $\prod$-diagrams in $F^*\C'$ are exactly those whose images under $\crt{F}C'$ are $\prod$-diagrams. In particular, $(\crt{F}\C',F)$ is a morphism of $(h^=)$-fibrations.

By the universal property of the pullback, given any other prefibration $\fibr{C}{C}{B}$ over $\B$ and any morphism $(\Phi,F)\colon\fib{C}\to\fib{C'}$ of prefibrations, we have a unique morphism $\cind{\Phi}\colon\fib{C}\to{}F^*\fib{C'}$ of prefibrations over $\B$ such that $\crt{F}\C'\cdot\cind{\Phi}=\Phi$. It follows from the above observations that $\cind{\Phi}$ is a morphism of ($h^=$-)fibrations whenever $\fib{C'}$ and $\fib{C}$ are ($h^=$-)fibrations and $(\Phi,\phi)$ is a morphism thereof.

If $\fib{C'}$ is a 1D2F, then $F^*\fib{C'}$ is defined as the pullback of the underlying fibration of $\fib{C'}$.

\subsubsection{}\label{defn:psnt-induced-stuff}\constr
Let $\B$ be a category, $\fibr{C'}{C'}{B'}$ a cloven 1D2F (i.e., the underlying fibration is cloven, \cite[\S1.3]{sentai}), let $F,G\colon\B\to{}\B'$ be functors, and let $\alpha\colon{}F\to{}G$ be a pseudonatural transformation. We will construct from this a morphism $\pnif\alpha\colon{}G^*\fib{C'}\to{}F^*\fib{C'}$ of fibrations over $\B$ between the associated pullback fibrations, as well as a pseudonatural transformation $\pnipn\alpha\colon\crt{F}\C'\circ\pnif{\alpha}\to\crt{G}\C'$ over $\alpha$, as shown below. We will define these simultaneously.
\[
  \begin{tikzcd}
    |[alias=FC]|F^*\C'\ar[rrrd, "\crt{F}\C'"]\ar[rddd, "F^*\fib{C'}"']
    &[-30pt]&[-30pt]&[20pt]\\[-10pt]
    &&&\C'\ar[dd, "\fib{C'}"]\\[-10pt]
    &&G^*\C'\ar[lluu, dashed, "\pnif\alpha"]
    \ar[ru, "\crt{G}\C'"' name=GC, near start]
    \ar[ld, "G^*\fib{C'}", near start]
    \ar[from=FC, to=GC, Rightarrow, dashed, "\pnipn\alpha", shorten <=10pt, shorten >=5pt]
    &\\[30pt]
    &\B
    \ar[rr, bend left, "F"{name=F}]
    \ar[rr, bend right, "G"'{name=G}]
    \ar[from=F, to=G, Rightarrow, "\alpha", shorten <=5pt, shorten >=5pt]
    &&\B'
  \end{tikzcd}
\]

Given $A\in\B$ and $(A,P)\in(G^*\fib{C'})^A$ (i.e., $P\in\fib{C'}^{GA}$), we set $\pnif\alpha(A,P)=(A,(\alpha_A)^*P)\in(F^*\fib{C'})^A$, and we set $\pnipn\alpha_{(A,P)}$ to be the morphism $\crt{(\alpha_A)}P\colon\alpha_A^*P\to{}P$.

Next, let $(f,p)\colon(A,P)\to(B,Q)$ be a morphism in $G^*\C'$. Since we want $\pnif\alpha(f,p)$ to be a morphism over $f\colon{}A\to{}B$ in $\B$, it must be of the form $(f,p')$ for some $p'\colon(\alpha_A)^*P\to{}(\alpha_A)^*Q$ over $Ff$. Seeing as we want there to be a 2-cell
\[
  \pnipn\alpha_{(f,p)}\colon
  \pnipn\alpha_{(B,Q)}\cdot(\crt{F}\C')(\pnif\alpha(f,p))
  =\crt{(\alpha_B)}Q\cdot{}p'
  \longrightarrow
  p\cdot\crt{(\alpha_A)}P
  =(\crt{G}\C')(f,p)\cdot\pnipn\alpha_{(A,P)}
\]
lying over $\alpha_f$, we see that we are forced to define $p'$ to be the unique morphism over $Ff$ for which there exists a 2-cell
\[
  \begin{tikzcd}
    (\alpha_A)^*P\ar[rr, dashed, "p'"]\ar[rd, "\crt{(\alpha_A)}P"']
    &&(\alpha_B)^*Q\car[rd, "\crt{(\alpha_B)}Q"]
    \ar[ld, Rightarrow, dashed, shorten <=5pt, shorten >=5pt]
    &\\
    &P\ar[rr, "p"']&&Q\\
    FA\ar[rr, "Ff"]\ar[rd, "\alpha_A"']
    &&FB\ar[rd, "\alpha_B"]
    \ar[ld, Rightarrow, shorten <=5pt, shorten >=5pt, "\alpha_f"']
    &\\
    &GA\ar[rr, "Gf"]&&GB
  \end{tikzcd}
\]
lying over $\alpha_f$ (such a $p'$ exists by Proposition~\ref{prop:1d2f-cart-gen} since $\alpha_{f}$ is invertible), and we must take $\pnipn\alpha_f$ to be this 2-cell lying over $\alpha_f$.

We now prove simultaneously that $\pnif\alpha$ is a functor and that $\pnipn\alpha$ is a pseudonatural transformation.

Let $(A,P)\tox{(f,p)}(B,Q)\tox{(g,q)}(C,R)$ be morphisms in $G^*\C'$. Let us write $(f,p')$ and $(g,q')$ for $\pnif\alpha(f,p)$ and $\pnif\alpha(g,q)$, as well as $(gf,(qp)')$ for $\pnif\alpha(gf,qp)$. We must show that  the 2-cells
\[
  \begin{tikzcd}[column sep=10pt, baseline=(afp.base)]
    \alpha_A^*P\ar[rr, "p'"]\ar[rd, "\crt{\alpha_A}P"']
    &&\alpha_B^*Q\ar[rd, "\crt{\alpha_B}Q"' {xshift=5pt}, near start]
    \ar[ld, Rightarrow, shorten <=5pt, shorten >=5pt, "\pnipn\alpha_{(f,p)}"' name=afp]
    \ar[rr,  "q'"]
    &&\alpha_C^*R\ar[rd, "\crt{\alpha_C}R"]
    \ar[ld, Rightarrow, shorten <=5pt, shorten >=5pt, "\pnipn\alpha_{(g,q)}"']
    &
    \\
    &P\ar[rr, "p"']&&Q
    \ar[rr, "q"']&&R\\
    FA\ar[rr, "Ff"]\ar[rd, "\alpha_A"']
    &&FB\ar[rd, "\alpha_B"' {xshift=5pt}, near start]
    \ar[ld, Rightarrow, shorten <=5pt, shorten >=5pt, "\alpha_f"']
    \ar[rr, "Fg"]
    &&FC\ar[rd, "\alpha_C"]
    \ar[ld, Rightarrow, shorten <=5pt, shorten >=5pt, "\alpha_g"']
    &\\
    &GA\ar[rr, "Gf"]&&GC
    \ar[rr, "Gg"]&&GC
  \end{tikzcd}
  \quad\text{and}\quad
  \begin{tikzcd}[column sep=10pt, baseline=(agfqp.base)]
    \alpha_A^*P\ar[rr, "(qp)'"]\ar[rd, "\crt{\alpha_A}P"']
    &&\alpha_C^*R\ar[rd, "\crt{\alpha_C}R"]
    \ar[ld, Rightarrow, shorten <=5pt, shorten >=5pt, "\pnipn\alpha_{(gf,qp)}"' name=agfqp]
    &\\
    &P\ar[rr, "qp"']&&R\\
    FA\ar[rr, "F(gf)"]\ar[rd, "\alpha_A"']
    &&FC\ar[rd, "\alpha_C"]
    \ar[ld, Rightarrow, shorten <=5pt, shorten >=5pt, "\alpha_{gf}"']
    &\\
    &GA\ar[rr, "G(gf)"]&&GC
  \end{tikzcd}
\]
are equal -- i.e., that
$(\id_q\circ{}\pnipn\alpha_{(f,p)})(\pnipn\alpha_{(g,q)}\circ{}\id_{p'})=\pnipn\alpha_{(gf,qp)}$ -- and also that $q'p'=(qp)'$. The first claim follows at once from the second by 1-discreteness. The second claim is true since $(qp)'$ is by definition the unique morphism $\alpha_A^*P\to\alpha_C^*R$ over $F(gf)$ for which there exists a 2-cell as above on the right, and $q'p'$ also has this property.

The proof of the remaining (unitality) property of $\pnif\alpha$ and $\pnipn\alpha$ is similar.

Finally, we must see that $\pnif{\alpha}$ preserves cartesian morphisms. This follows from the definition of $\pnif{\alpha}$, Proposition~\ref{prop:1d2f-2-of-3}, and the fact that a morphism $(f,p)$ in $F^*\C'$ is cartesian if and only if $p$ is.

\subsubsection{}\label{defn:fiberwise-equivalence}\defn
A morphism $(\Phi,\phi)\colon\fibr{C}CB\to\fibr{C'}{C'}{B'}$ is a \emph{fiberwise equivalence} if the induced functor $\fib{C}^A\to\fib{C'}^{\phi A}$ is an equivalence for each $A\in\B$.

\subsubsection{}\label{prop:induced-stuff-equivalence}\prop
With $\fib{C}$, $F$, $G$, and $\alpha$ as in Construction~\ref{defn:psnt-induced-stuff},
if $\alpha$ is a pseudonatural equivalence, then $\pnif{\alpha}$ is a fiberwise equivalence and $\pnipn\alpha$ is a pseudonatural equivalence.

Also, if $\alpha$ is natural transformation or natural isomorphism, then so is $\pnipn{\alpha}$.
\pf
To see that $\pnif{\alpha}$ is a fiberwise equivalence, note that the induced functor $\pnif\alpha\colon(G^*\fib{C})^A\to(F^*\fib{C})^A$ is (with respect to the identifications $(F^*\fib{C})^A\cong\fib{C}^{FA}$ and $(G^*\fib{C})^A\cong\fib{C}^{FA}$) just the pullback functor $(\alpha_A)^*\colon\fib{C}^{GA}\to\fib{C}^{FA}$. This is an equivalence, since, given a quasi-inverse $\beta_A$ for $\alpha_A$, we obtain a quasi-inverse $(\beta_A)^*$ for $(\alpha_A)^*$.

We now prove that $\pnipn\alpha$ is a pseudonatural equivalence. Since each $\alpha_A$ is an equivalence in $\B'$ and each $\pnipn\alpha_A$ is a cartesian lift of $\alpha_A$, it suffices to prove that any cartesian lift of an equivalence in a 1D2F is again an equivalence.

Let $f\colon{}A\to{}B$ be an equivalence with quasi-inverse $g\colon{}B\to{}A$, so that there exist isomorphism 2-cells $\alpha\colon\id_A\toi{}gf$ and $\beta\colon\id_B\toi{}fg$, and let $p\colon{}P\to{}Q$ be a cartesian lift of $f$. By Proposition~\ref{prop:1d2f-cart-gen}, there exists a unique morphism $q\colon{}Q\to{}P$ over $g$ for which there exists a (necessarily invertible) 2-cell $\id_Q\toi{}pq$ over $\beta$. It remains to see that $qp\cong\id_P$. By Proposition~\ref{prop:1d2f-2-of-3}, $q$ is cartesian, and hence, by the argument we just gave, there exists $p'\colon{}P\to{}Q$ with $qp'\cong\id_P$. We then have $p\cong{}pqp'\cong{}p'$ and hence $qp\cong{}qp'\cong\id_P$.

The last statement follows immediately from the construction of $\pnipn\alpha$ and the fact that any cartesian lift of an isomorphism is an isomorphism.
\qed

\subsubsection{}\label{defn:free-hfib}\defn
An $h^=$-fibration $\fibr{C}CB$ is \emph{free over $\B$} if, for any $h^=$-fibration $\fibr{C'}{C'}{B}$ over $\B$, there is up to isomorphism a unique morphism $\fib{C}\to\fib{C'}$ of $h^=$-fibrations over $\B$; i.e., there exists such a morphism, and for any two such, there exists a natural isomorphism over $\B$ between them.\footnote{In fact, a free $h^=$-fibration over $\B$ has the (categorically natural) stronger property that for any two morphisms $\fib{C}\to\fib{C'}$ there is a \emph{unique} isomorphism between them, as our Construction~\ref{constr:free-hfib} will show, but we will not need this.}

In the \hyperref[sec:appendix]{appendix}, it is shown that there exists a free $h^=$-fibration over any f.p.\ category $\B$.

We see from the definition that any two free $h^=$-fibrations over $\B$ are equivalent over $\B$.

We will now show that any such $\fib{C}$ also satisfies a universal property with respect to any $h^=$-fibration $\fibr{C'}{C'}{B'}$ (with possibly different base $\B'$).

\subsubsection{}\prop\label{prop:global-univ-prop}
If $\fibr{C}CB$ is a free $h^=$-fibration over $\B$ and $\fibr{C'}{C'}{B'}$ is any $h^=$-fibration, then for any f.p.\ functor $\phi\colon\B\to\B'$, there exists a morphism of $h^=$-fibrations $(\Phi,\phi)\colon\fib{C}\to\fib{C}'$ over $\phi$. Moreover (though we will not need this), $\Phi$ is determined uniquely up to natural isomorphism over $\phi$.
\pf
Taking the pullback $(\crt{\phi}{\C'},\phi)\colon\phi^*\fib{C'}\to\fib{C}$ of $\fib{C'}$ along $\phi$, we have by the discussion in Definition~\ref{defn:pullback-fib} that $\phi^*\fib{C'}$ is an $h^=$-fibration. Hence, by the freeness of $\fib{C}$, we then have a morphism $\Psi\colon\fib{C}\to\phi^*\fib{C'}$ of $h^=$-fibrations over $\B$. Since the composite of morphism of $h^=$-fibrations is again one, $(\crt{\phi}{\C'}\circ\Psi,\phi)$ is as desired.

For uniqueness: given morphisms $(\Phi,\phi),(\Phi',\phi)\colon\fib{C}\to\fib{C'}$ of $h^=$-fibrations, by the universal property of $\phi^*\fib{C'}$, we have morphisms $\Psi,\Psi'\colon\fib{C}\to\phi^*\fib{C'}$ over $\B$ with $\Phi=\crt{\phi}{\C'}\circ\Psi'$ and $\Phi'=\crt{\phi}{\C'}\circ\Psi'$. By the freeness of $\fib{C}$, we have a natural isomorphism $\Psi\toi\Psi'$ over $\B$, and whiskering with $\crt{\phi}{\C'}$ gives a natural isomorphism $\Phi\toi\Phi'$ over $\phi$ as desired.
\qed

\subsubsection{}\label{thm:abstract-invariance}\thm
Suppose $\fibr{C}CB$ is a free $h^=$-fibration, $\fibr{C'}{C'}{B'}$ is a 1D2F which is also an $h^=$-fibration, $(\Phi,\phi),(\Psi,\psi)\colon\fib{C}\to\fib{C'}$ are morphisms of $h^=$-fibrations, and $\alpha\colon\phi\to\psi$ is a pseudonatural equivalence.

Then there exists a pseudonatural equivalence $\Phi\to\Psi$ lying over $\alpha$ (which, if $\alpha$ is a natural isomorphism, can be taken to be a natural isomorphism).

In particular, for each $P\in\C$ over some $A\in\B$, there is an equivalence $p\colon\Phi{}P\to\Psi{}P$ lying over the equivalence $\alpha_A\colon\phi{}A\to\psi{}A$ (note that the natural additional condition that $p$ have a quasi-inverse $q$ lying over a quasi-inverse of $\alpha_A$, so that the associated 2-cells $\id_{\Phi{}P}\toi{}qp$ and $\id_{\Psi{}P}\toi{}pq$ in $\C$ lie over the corresponding ones in $\B$ -- cf. Definition~\ref{defn:homotopy-over} -- is satisfied automatically) -- and if $\alpha$ (and hence each $\alpha_A$) is an isomorphism, then $p$ can be taken to be one as well.
\pf
To begin, we choose a cleavage of $\fib{C'}$ (which we can do by the axiom of choice -- otherwise, we must assume explicitly that $\fib{C'}$ admits a cleavage). Consider the pullbacks of $\fib{C'}$ along $\phi$ and $\psi$. We then have the situation depicted in Construction~\ref{defn:psnt-induced-stuff}, with $F=\phi$ and $G=\psi$. By Proposition~\ref{prop:induced-stuff-equivalence} the morphism $\pnif\alpha\colon\psi^*\fib{C'}\to{}\phi^*\fib{C'}$ of fibrations is a fiberwise equivalence and hence, by Lemma~\ref{lem:equivs-are-hfib-homos} below, a morphism of $h^=$-fibrations.

By the definition of the pullback, we have morphisms of fibrations $\cind{\Phi}\colon\fib{C}\to\phi^*\fib{C}$ and $\cind{\Psi}\colon\fib{C}\to\psi^*\fib{C}$ over $\B$ such that $\crt{\phi}\C'\circ\cind\Phi=\Phi$ and $\crt{\psi}\C'\circ\cind\Psi=\Psi$ -- and by the discussion in Definition~\ref{defn:psnt-induced-stuff}, these are morphisms of $h^=$-fibrations. Hence, since the composite of morphisms of $h^=$-fibrations over $\B$ is again one, we have that $\pnif\alpha\circ\cind{\Psi}\colon\fib{C}\to\phi^*\fib{C}'$ is a morphism of $h^=$-fibrations. Hence, by the freeness of $\fib{C}$, we have a natural isomorphism $\eta\colon\cind{\Phi}\to\pnif{\alpha}\cind{\Psi}$ of morphisms of $h^=$-fibrations over $\B$:
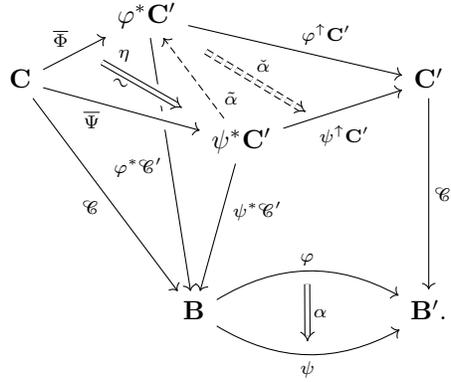
\[
  \begin{tikzcd}
    &|[alias=FC]|\phi^*\C'\ar[rrrd, "\crt{\phi}\C'"]\ar[rddd, "\phi^*\fib{C'}"' pos=0.45]
    &[-30pt]&[-30pt]&[20pt]\\[-10pt]
    \C\ar[ru, "\cind{\Phi}" name=Phi]\ar[rrrd, "\cind{\Psi}"' pos=0.4, crossing over]
    \ar[rrdd, "\fib{C}"']
    &&&&\C'\ar[dd, "\fib{C'}"]\\[-10pt]
    &&&
    |[alias=psiC]|\psi^*\C'\ar[lluu, dashed, "\pnif\alpha"' pos=0.1]
    \ar[from=Phi, to=psiC, Rightarrow, shorten <=12pt, shorten >=8pt, crossing over,
        "\eta" pos=0.3, "\sim"' sloped]
    \ar[ru, "\crt{\psi}\C'"' name=GC, near start]
    \ar[ld, "\psi^*\fib{C'}", near start]
    \ar[from=FC, to=GC, Rightarrow, dashed, "\pnipn\alpha", shorten <=10pt, shorten >=5pt]
    &\\[30pt]
    &&\B
    \ar[rr, bend left, "\phi"{name=F}]
    \ar[rr, bend right, "\psi"'{name=G}]
    \ar[from=F, to=G, Rightarrow, "\alpha", shorten <=5pt, shorten >=5pt]
    &&\B'.
  \end{tikzcd}
\]
Hence, the ``pasting'' of $\eta$ and $\pnipn{\alpha}$ -- i.e. the composite pseudonatural transformation $(\pnipn\alpha\circ\id_{\cind{\Psi}})\circ(\id_{\crt{\phi}{\C'}}\circ\eta)$ -- is as desired (since it is an equivalence and lies over $\alpha$).

In case $\alpha$ is a natural isomorphism, then Proposition~\ref{prop:induced-stuff-equivalence} says that $\pnipn\alpha$ is a natural isomorphism, and hence so is the pasting of $\eta$ and $\pnipn\alpha$.
\qed

\subsubsection{}\label{lem:equivs-are-hfib-homos}\lem
Suppose $\fibr{C}CB$ and $\fibr{C'}{C'}{B'}$ are $h^=$-fibrations, and $(\Phi,\phi)\colon\fibr{C}CB\to\fibr{C'}{C'}{B'}$ is a morphism of fibrations which is a fiberwise equivalence with $\phi$ an f.p.\ functor. Then $(\Phi,\phi)$ is a morphism of $h^=$-fibrations.
\pf
$\phi$ is product-preserving by assumption, and the induced functors $\fib{C}^B\to\fib{C'}^{\phi{B}}$ are clearly bi-cartesian closed since they are equivalences. Hence, it only remains to see that the relevant co-cartesian morphisms and $\prod$-diagrams are preserved.

In fact, all co-cartesian morphisms and $\prod$-diagrams are preserved. This follows easily from the fact that, in the present situation, given a morphism $f\colon{}A\to{B}$ in $\B$ and objects $P$ and $Q$ in $\fib{C}^A$ and $\fib{C}^B$, $\Phi$ induces a bijection between morphisms $P\to{Q}$ lying over $f$ and morphisms $\Phi{P}\to\Phi{Q}$ lying over $\phi{f}$, and moreover that $p\colon{}P\to{Q}$ is cartesian if and only if $\Phi{}p$ is.
\qed

\subsection{The abstract invariance theorem: second part}
The theorem proven in the previous section was stated with respect to an arbitrary free $h^=$-fibration $\fib{C}$ over an f.p.\ category $\B$. However, in practice, we are interested in the particular case of $\B$ being a \emph{free} f.p.\ category on a given signature $\sigma$.

As mentioned above, the definition of the latter is arranged so that the f.p.\ functors $\B\to\B'$ correspond precisely to the interpretations $\sigma\to\B'$, \emph{and} so that the natural isomorphisms between functors $\B\to\B'$ correspond to isomorphisms of interpretations $\sigma\to\B'$.

However, for the purpose of the homotopy invariance for the homotopical semantics, we will want to start with a \emph{homotopy equivalence} of interpretations $\sigma\to\B'$ -- which, in the context of a general 2-category $\B'$, corresponds to a \emph{pseudonatural equivalence} of interpretations $\sigma\to\B'$ (Definition~\ref{defn:2cat-homomorphism} below).

Hence, the second part of the abstract invariance theorem will show that a pseudonatural equivalence of interpretations $\sigma\to\B'$ induces a pseudonatural equivalence between the induced functors $\B\to\B'$. Note that this is analogous to the first part: in both cases, we have a universal property with respect to natural isomorphisms, and want to show that it still holds for pseudonatural equivalences.

The proof proceeds by considering a modified arrow category $\peqat{(\B')}$, such that the f.p.\ functors into $\peqat{(\B')}$ are pseudonatural equivalences of f.p.\ functors into $\B'$, and the $\sigma$-interpretations in $\peqat{(\B')}$ are homotopy-equivalences of $\sigma$-interpretations into $\B'$; this reduces the claim to the original freeness property of $\B$.

\subsubsection{}\defn
Given a 2-category $\C$, we define $\pacat\C$, the \emph{pseudo-arrow category} of $\C$, to have as objects functors $\mathbf{2}\to\C$ -- i.e., 1-cells in $\C$ -- and to have as morphisms pseudonatural transformations, with composition being given by composition of pseudonatural transformations which, as we leave to the reader to verify, is associative and has identities

In other words a morphism $\alpha$ from $f_1\colon{}A_1\to{}B_1$ to $f_2\colon{}A_2\to{}B_2$ in $\pacat\C$ is a triple $(\alpha_A,\alpha_B,\alpha_f)$ with $\alpha_A\colon{}A_1\to{}A_2$, $\alpha_B\colon{}B_1\to{}B_2$, and $\alpha_f\colon\alpha_B\circ{}f_1\tox{\sim}f_2\circ{}\alpha_A$.

There are obvious domain and codomain functors $\dom,\cod\colon\pacat\C\to\C$.

We define $\peqat\C$ to be the full subcategory of $\pacat\C$ with objects the equivalences in $\C$.

\subsubsection{}\prop
Given a 2-category $\C$, if $\C$ has finite 2-categorical products (see \cite[Definition~7.1]{sentai}), then the categories $\pacat\C$ and $\peqat\C$ have finite products. Moreover, the inclusion functor $\peqat\C\hookrightarrow\pacat\C$, as well as the functors $\dom,\cod\colon\pacat\C\to\C$, are f.p.\ functors.
\pf
Given 1-cells $g\colon{}A\to{}C$ and $h\colon{}B\to{}D$ and products $A\times{}B$ and $C\times{}D$ in $\C$, we will show that the morphisms $g\times{}h\to{}g$ and $g\times{}h\to{}h$ in $\pacat\C$ given by $(\pi_1,\pi_1,\id_{g\pi_1})$ and $(\pi_2,\pi_2,\id_{h\pi_2})$ exhibit $g\times{}h$ as a product of $g$ and $h$. It follows that $\dom$ and $\cod$ are f.p.

Given a 1-cell $f\colon{}X\to{}Y$ in $\C$ and morphisms $(s,t,\alpha)\colon{}f\to{}g$ and $(u,v,\beta)\colon{}f\to{}h$ in $\pacat\C$, we must show that there is a unique morphism $(w,x,\gamma)\colon{}f\to{}g\times{}h$ such that $\pi_1w=s$, $\pi_2w=u$, $\pi_1x=t$, $\pi_2x=v$, $\id_{\pi_1}\circ\gamma=\alpha$, and $\id_{\pi_2}\circ\gamma=\beta$. Clearly, we must take $w=\br{s,u}$ and $x=\br{t,v}$. Then, since $C\times{}D$ is a 2-categorical product, there is a unique 2-cell $\gamma\colon\br{gs,hu}\to\br{tf,vf}$ such that $\id_{\pi_1}\circ\gamma=\alpha$ and $\id_{\pi_2}\circ\gamma=\beta$.

Finally, we must see that if $f$ and $g$ are equivalences in $\C$, then $f\times{}g$ is one as well. In fact, if $f\I$ and $g\I$ are quasi-inverses to $f$ and $g$, then $f\I\times{}g\I$ is easily seen to be a quasi-inverse to $f\times{}g$.
\qed

\subsubsection{}\label{prop:htpy-equiv-funs}\prop
Given a category $\C$, a 2-category $\D$, and functors $F,G\colon\C\to\D$, the functors $F$ and $G$ are pseudonaturally equivalent if and only if there exists a functor $H\colon\C\to\peqat\D$ such that $\dom\circ{}H=F$ and $\cod\circ{}H=G$.
\pf
The proof is by inspection, the point being that the data of a functor $H\colon\C\to\peqat\D$ is precisely the data of a pseudonatural equivalence $\dom\circ{}H\to\cod\circ{}H$, and the condition that $H$ be a functor is equivalent to the given data defining a pseudonatural transformation.
\qed

\subsubsection{}\label{defn:2cat-homomorphism}\defn
Given a signature $\sigma$ and a 2-category $\C$ with finite 2-categorical products, the underlying category of $\C$ also has finite products, and we can consider interpretations $\sigma\to\C$. Given two such interpretations $M,N\colon\sigma\to\C$, a \emph{pseudohomomorphism} $\alpha\colon{}M\to{}N$ consists of a 1-cell $\alpha_A\colon{}M\to{}N$ for each $A\in\Ob\sigma$, such that there exists an invertible 2-cell $\alpha_B\circ{}Mf\to{}Nf\circ\alpha_{\vec{A}}$ for each $f\in\sigma(\vec{A},B)$ (where $\alpha_{\vec{A}}$ is defined as in Definition~\ref{defn:homomorphism}).

A pseudonatural homomorphism $\alpha\colon{}M\to{}N$ is a \emph{pseudoequivalence}\footnote{Perhaps ``pseudoisomorphism'' would be more appropriate?} if $\alpha_A$ is an equivalence in $\C$ for each $A\in\Ob\sigma$. If such a $\alpha$ exists, we say that $M$ and $N$ are \emph{pseudoequivalent}.

Note that when $\C$ is the category of cofibrant-fibrant objects in a model category, endowed with the Quillen 2-categorical structure \cite[Definition~15.7]{sentai}, then the notions of pseudohomomorphism and pseudoequivalence specialize to those of homotopy homomorphism and homotopy equivalence (Definition~\ref{defn:homomorphism}).

\subsubsection{}\label{prop:htpy-equiv-sig-ints}\prop
Two interpretations $M,N\colon\sigma\to\C$ are pseudoequivalent if and only if there exists an interpretation $H\colon\sigma\to\peqat\C$ such that $\dom\circ{}H=M$ and $\cod\circ{}H=N$.
\pf
Given a pseudoequivalence $\alpha\colon{}M\to{}N$, we define the interpretation $H\colon\sigma\to\peqat\C$ by setting
$HA=\alpha_A\colon{}MA\to{}NA$ for $A\in\Ob\sigma$;
$H\vec{A}=\alpha_{\vec{A}}\colon{}M\vec{A}\to{}N\vec{A}$ and $\pi_i^H=(\pi_i^M,\pi_i^N,\id)$ for $\vec{A}\in(\Ob\sigma)^{<\omega}$;
and for $f\in\sigma(\vec{A},B)$, we take $Hf$ to be $(Mf,Nf,\beta)$, where $\beta\colon\alpha_B\circ{}Mf\to{}Nf\circ\alpha_{\vec{A}}$ is any invertible 2-cell (which exists since $\alpha_A$ is a pseudohomomorphism).

The proof of the converse (which we will not need) is similar; given $H$, we define $\alpha_A=HA$ for $A\in\Ob\sigma$. The fact that $\alpha$ satisfies the required property then follows from the existence of the morphisms $Hf\in\Ar\pacat\C$ for $f\in\sigma(\vec{A},B)$, and the fact that $H\vec{A}$ and $\alpha_{\vec{A}}$ are \emph{isomorphic} (but not necessarily equal!) for each $\vec{A}\in(\Ob\sigma)^{<\omega}$.
\qed

\subsubsection{}\defn\label{defn:free-fp-cat}
Given a signature $\sigma$ and a f.p.\ category $\C$, the interpretations $\sigma\to\C$ and homomorphisms form a category in an obvious manner, which we denote by $\fpints(\sigma,\C)$.

Given two f.p. functors $F,G\colon\C\to\D$ and a natural transformation $\alpha\colon{}F\to{}G$, we have the interpretations $F\circ M,G\circ M\colon\sigma\to\D$, and we obtain a homomorphism $\alpha\circ{}M\colon{}F\circ{}M\to{}G\circ{}M$ by setting $(\alpha\circ{}M)_A=\alpha_{MA}$.

This defines, for each $M\colon\sigma\to\C$, a functor $(\text{--}\circ M)\colon\FPFun(\C,\D)\to\fpints(\sigma,\D)$.

Given an interpretation $i\colon\sigma\to\C$, we say that the f.p.\ category $\C$ is \emph{free on $\sigma$ (via $i$)} if $(\text{--}\circ i)\colon\FPFun(\C,\D)\to\fpints(\sigma,\D)$ is an isomorphism of categories for each f.p.\ category $\D$.\footnote{It would be more natural, perhaps, to demand that this is only an equivalence, and not an isomorphism, but it will be convenient to assume the stronger property.} Note that this property determines $\C$ up to \emph{isomorphism}.

It is well-known that there exists a free f.p.\ category on any signature $\sigma$ (see Proposition~\ref{prop:free-fp-exists} in the \hyperref[sec:appendix]{appendix}).

\subsubsection{}\label{thm:free-htpy-equiv}\thm
Let $\C$ be a free f.p.\ category on the algebraic signature $\sigma$ via $i\colon\sigma\to\C$, and let $M,N\colon\sigma\to\D$ be two interpretations into an f.p.\ 2-category $\D$. Then $M$ and $N$ are pseudoequivalent if and only if the induced f.p.\ functors $\widetilde{M},\widetilde{N}\colon\C\to\D$ are pseudonaturally equivalent.

Moreover, for any pseudoequivalence $\alpha\colon M\to N$, there is a pseudonatural equivalence $\tilde\alpha\colon\widetilde{M}\to\widetilde{N}$ with $\tilde\alpha_{iA}=\alpha_{A}$ for all $A\in\Ob\sigma$ (and vice versa).
\pf
We have, by definition, that $\widetilde{M}\circ{}i=M$ and $\widetilde{N}\circ{}i=N$.

Given a pseudonatural equivalence $\tilde\alpha\colon\widetilde{M}\to\widetilde{N}$, we thus have a pseudoequivalence
\[
  M=\widetilde{M}\circ{}i\tox{\tilde\alpha\circ{}i}\widetilde{N}\circ{}i=N,
\]
where $\tilde\alpha\circ{}i$ is the pseudoequivalence given by $(\tilde\alpha\circ{}i)_A=\tilde\alpha_{iA}$ for each $A\in\Ob\sigma$. To see that $(\tilde\alpha\circ i)_B\circ Mf\cong Nf\circ(\tilde\alpha\circ i)_{\vec{A}}$ for $f\in\sigma(\vec{A},B)$, first note that $(\tilde\alpha\circ i)_B\circ Mf=\tilde\alpha_{iB}\circ\widetilde{M}(if)\cong
\widetilde{N}(if)\circ\tilde\alpha_{i\vec{A}}=
Nf\circ\tilde\alpha_{i\vec{A}}$ by pseudonaturality of $\tilde\alpha$. But one can also see that $\tilde\alpha_{i\vec{A}}\cong(\tilde\alpha\circ i)_{\vec{A}}$ using the pseudonaturality of $\tilde\alpha$, the product-preservation of $N$, and the 2-categorical universal property of the product $N\vec{A}$, and hence $Nf\circ\alpha_{i\vec{A}}\cong Nf\circ(\alpha\circ i)_{\vec{A}}$.

Conversely, given a pseudoequivalence $\alpha\colon{}M\to{}N$, we have by Proposition~\ref{prop:htpy-equiv-sig-ints} an interpretation $H\colon\sigma\to{}\peqat\D$ with $\dom\circ{}H=M$ and $\cod\circ{}H=N$. Hence, we have an induced f.p.\ functor $\widetilde{H}\colon\C\to\peqat\D$ with $\widetilde{H}\circ{}i=H$ and hence
\[
  (\dom\circ\widetilde{H})\circ{}i=M\quad\text{and}\quad
  (\cod\circ\widetilde{H})\circ{}i=N
\]
(here, we are using that composition of f.p.\ functors is associative with composition of an interpretation and an f.p.\ functor). Hence, using the freeness of $\C$ again, we have $\dom\circ\widetilde{H}=\widetilde{M}$ and $\cod\circ\widetilde{H}=\widetilde{N}$. By Proposition~\ref{prop:htpy-equiv-funs}, we have a pseudonatural equivalence $\tilde{\alpha}\colon\dom\circ\widetilde{H}\to\cod\circ\widetilde{H}$, and hence we obtain the desired pseudonatural equivalence
\[
  \widetilde{M}=
  \dom\circ\widetilde{H}\tox{\tilde\alpha}
  \cod\circ\widetilde{H}=
  \widetilde{N}
  .
\]
The ``moreover'' statement is clear from the two constructions just given.
\qed

\subsection{The special invariance theorem}\label{subsec:special-invariance}
We now show how to derive syntactic conclusions from the abstract invariance theorem. This follows the general pattern in abstract/categorical logic (described, for example, in \cite{makkai-generalized-sketches-1}) of a theorem about syntax and a purely algebraic theorem being related by an appropriate ``translation'' theorem, describing the relationship between the algebraic structures and the syntax being considered.

In the present case, if we were to follow the usual procedure, this would involve showing that the free f.p.\ category $\B$ on a signature $\sigma$ and the free $h^=$-fibration $\fibr{C}CB$ over $\B$ can be constructed ``out of the syntax'' of first-order logic over $\sigma$ -- specifically, in such a way that the morphisms of $\B$ are given by the $\sigma$-terms, and the objects of $\C$ by the formulas (the morphisms then being certain equivalence classes of formal deductions).

In the \hyperref[sec:appendix]{appendix}, we will (mostly) show that $\fib{C}$ indeed admits such a description. \emph{However}, for the purposes of deriving the special invariance theorem from the abstract one, this description is actually not needed; we \emph{only} need the \emph{existence} of the free f.p.\ category $\B$ and $h^=$-fibration $\fib{C}$. That is, the only thing we use about these objects is their defining universal property, and nothing else. We would only need the full syntactic description of $\fib{C}$ if we wanted to go in the other direction and deduce the abstract invariance theorem from the syntactic one (this would also require ``packing more'' into the syntactic invariance theorem than is really natural from the syntactic point of view).

The reason that we don't need the syntactic description of $\fib{C}$ is that, even without it, we can still \emph{interpret} logic over $\sigma$ into $\fib{C}$, and transport this interpretation along any morphism of $h^=$-fibrations $\fib{C}\to\fib{C'}$. Knowing that two such morphisms are pseudonaturally equivalent is then enough to deduce the desired invariance property for the resulting interpretations in $\fib{C'}$.

We first work in the context of a general $h^=$-fibration $\fib{C'}$ (which automatically inherits a 1D2F structure from the results of \cite{sentai}), showing that the abstract invariance theorem implies a syntactic invariance theorem for interpretations in $\fib{C'}$, involving the 2-categorical structure. We then deduce the invariance theorem for the homotopical semantics by specializing to the case $\fib{C'}=\fib{HoF_\fb(\Kan)}$.

This involves translating between the 2-categorical structure on $\fib{HoF_\fb(\Kan)}$ and the model structure on $\Kan$ (or rather $\sSet$). The main point is that the resulting 2-categorical structure on $\Kan$ is the same as the Quillen 2-categorical structure, as was proven in \cite[Part~IV]{sentai}. However, we also need to know that the 2-cells in the total category of $\fib{HoF_\fb(\Kan)}$ are given by homotopies lying over homotopies in the base, as defined in \S\ref{subsec:htpy-equiv}. This result really ``should'' have been included in \emph{op. cit.}, but it was not, so we prove it here.

\subsubsection{}\thm\label{thm:somewhat-abstract}
Let $\sigma$ be a signature, let $\fibr{C'}{C'}{B'}$ be an $h^=$-fibration, and let $(\widehat M,M),(\widehat N,N)\colon\sigma\to\C'$ be interpretations in $\fib{C'}$.

Now suppose that $\fib{C'}$ has a 1D2F structure such that $\B$ has 2-categorical products (for example, it could be the ``universal'' 1D2F structure from \cite{sentai}).

Then given any pseudoequivalence $\alpha\colon M\to N$, there exists, for each formula-in-context $(\phi,\vec{x})$ over $\sigma$, an equivalence $\widehat{M}_{\vec{x}}(\phi)\to\widehat{N}_{\vec{x}}(\phi)$ in $\C'$ lying over the equivalence $\alpha_{\tp\vec{x}}\colon M(\vec{x})\to N(\vec{x})$ in $\B'$.
\pf
Let $\B$ be a free f.p.\ category on $\sigma$ (via $i\colon\sigma\to\B$) and let $\fibr{C}CB$ be a free $h^=$-fibration over $\B$.

By the universal properties of each of these (and Proposition~\ref{prop:global-univ-prop}), we have f.p.\ functors $\phi,\psi\colon\B\to\B'$ with $\phi\circ i=M$ and $\psi f\circ i=N$, and morphisms $(\Phi,\phi),(\Psi,\psi)\colon\fib{C}\to\fib{C'}$ of $h^=$-fibrations over $\phi$ and $\psi$.

Now let $\alpha\colon M\to N$ be a pseudoequivalence. By Theorem~\ref{thm:free-htpy-equiv}, we have a pseudonatural equivalence $\tilde\alpha\colon\phi\to\psi$ with $\tilde\alpha_{iA}=\alpha_A$ for $A\in\Ob\sigma$, from which it follows that $\tilde\alpha_{i\vec{A}}\cong\alpha_{\vec{A}}$ for $\vec{A}\in(\Ob\sigma)^{<\omega}$, as explained in the proof of that theorem.

Next, by Theorem~\ref{thm:abstract-invariance}, we have for each $A\in\B$ and $P\in\C$ an equivalence $\beta_P\colon\Phi P\to\Psi P$ in $\C'$ lying over the equivalence $\alpha_A\colon MA\to NA$ in $\B'$.

Now let $\hat\imath$ be an interpretation in $\fib{C}$ over $i$. Then $\Phi\circ\hat\imath$ and $\Psi\circ\hat\imath$ are interpretations in $\fib{C'}$ over $M$ and $N$, and hence are isomorphic to $\widehat M$ and $\widehat N$, respectively.

Given a formula-in-context $(\chi,\vec{x})$, we have the object $\hat\imath_{\vec{x}}(\chi)$ in $\C$, and hence an equivalence
$\beta_{\hat\imath_{\vec{x}}(\chi)}\colon
\Phi(\hat\imath_{\vec{x}}(\chi))\to\Psi(\hat\imath_{\vec{x}}(\chi))$ over the equivalence $\tilde\alpha_{i(\vec{x})}\colon M(\vec{x})\to N(\vec{x})$. Using the isomorphisms $\Phi(\hat\imath_{\vec{x}}(\chi))\cong\widehat M_{\vec{x}}(\chi)$ and $\Psi(\hat\imath_{\vec{x}}(\chi))\cong\widehat N_{\vec{x}}(\chi)$ over $M(\vec{x})$ and $N(\vec{x})$, and the isomorphism $\tilde\alpha_{i(\vec{x})}\cong\alpha_{\tp\vec{x}}$, it follows that there is also an equivalence $\widehat M_{\vec{x}}(\chi)\to\widehat N_{\vec{x}}(\chi)$ over the equivalence $\alpha_{\tp\vec{x}}$, as desired.
\qed

\subsubsection{}\label{thm:2-cells-in-hofb-tot}\thm
Let $\C$ be a model category, and consider the $h^=$-fibration $\fib{HoF_\fb(\C_\fb)}$ with its 1D2F structure from \cite[Theorem~6.11]{sentai}. Let
$
\begin{tikzcd}[row sep=0pt, column sep=10pt]
  &\ar[dd, shorten >=2pt, shorten <=-2pt, Rightarrow, "\ \alpha", pos=0.4]&\\[-5pt]
  A\ar[rr, "f", bend left=30pt]
  \ar[rr, "g"', bend right=30pt]&&
  B\\
  &{}&
\end{tikzcd}
$
be a 2-cell in $\C_\fb$ with $A$ cofibrant. According to \cite[\S15]{sentai}, such 2-cell is given by a homotopy class of homotopies $f\to{}g$.

Next, let $(p,f),(q,g)\colon{}(X,A,x)\to{}(Y,B,y)$ be morphisms in $(\C^\to)_\fb$ (i.e., morphisms in $\C^\to$ with $x$ and $y$ fibrant -- note that we are using, as usual, the model structure on $\C^\to$ from \cite[Definition~11.1]{sentai}, not the one referenced in the proof of Proposition~\ref{prop:homotopy-over-left-right}) and assume $X$ is cofibrant.
\[
  \begin{tikzcd}
    X\ar[r, "p", shift left]\ar[r, "q"', shift right]\ar[d, "x"']&Y\ar[d, "y"]\\
    A\ar[r, "f", shift left]\ar[r, "g"', shift right]&B.
  \end{tikzcd}\quad\quad
\]
Recall from \cite[Definition~11.7]{sentai} that we write $\gamma(p,f),\gamma(q,g)\colon{}(X,A,x)\to(Y,B,y)$ for the images of $(p,f)$ and $(g,q)$ in $\Ho(\C^\to)_\fb$.

We now claim that there exists a 2-cell $\gamma(p,f)\to{}\gamma(q,g)$ over $\alpha$ if and only if there is a homotopy from $p$ to $q$ over $\alpha$, regarded as a (homotopy class of) homotopy from $f$ to $g$ (this only depends on the homotopy class $\alpha$ by Proposition~\ref{prop:homotopy-over-left-right}).
\pf
Let $k\colon{}A\to{}B^I$ be a representative of $\alpha$.

We first unwind the definition of the 2-cells in the total category of $\fib{HoF_\fb(\C_\fb)}$. Let us write $P$ and $Q$ for $(X,x)$ and $(Y,y)$, respectively. By definition (see \cite[Theorem~6.11~and~Definition~6.4]{sentai}), there exists a 2-cell $\gamma(p,f)\to\gamma(q,g)$ if and only if the diagram
\[
  \begin{tikzcd}
    &\pi_1^*Q\wedge\Eq_B\ar[d, "\nat_B^Q"]\\
    P\ar[r, "\cind{\gamma(q,g)}"']\ar[ru, "\brr{\cind{\gamma(p,f)},\alpha!}"]&\pi_2^*Q\\[-5pt]
    A\ar[r, "\br{f,g}"]&B\times{}B
  \end{tikzcd}
\]
commutes. Now, inserting the definitions of $P$ and $Q$, and of the pullback functors $\pi_1^*$ and $\pi_2^*$ and the equality object $\Eq_B$ in the $\wedgeq$-cloven $\wedgeq$-fibration $\fib{HoF_\fb(\C_\fb)}$ (see \cite[\S15.1]{sentai}), this is the same as diagram on the left below (where, for conciseness, we are identifying objects and morphisms in $\C^\to$ with their domains), where the fiber product $Y\times_B{}B^I$ is taken with respect to $d_1\colon{}B^I\to{}B$.
\begin{equation}\label{eq:2cell-ex-char}
  \begin{tikzcd}
    &Y\times_BB^I\ar[d, "\nat_B^Q"]\\
    X\ar[r, "\gamma(\br{fx,q})"']\ar[ru, "\gamma(\br{p,kx})"]&B\times{}Y
    \\[-5pt]
    A\ar[r, "\br{f,g}"]&B\times{}B
  \end{tikzcd}
  \quad\quad\quad
  \begin{tikzcd}
    &Z\ar[d, "\nu"]\\
    X\ar[r, "\br{fx,q}"']\ar[ru, "z"]&B\times{}Y
    \\[-5pt]
    A\ar[r, "\br{f,g}"]&B\times{}B
  \end{tikzcd}
\end{equation}
We would now like to reformulate the commutativity of this diagram in simpler terms.

Note that the morphism $\br{y,sy}\colon Y\to Y\times_BB^I$ is a weak equivalence, since composing it with the trivial fibration $\pi_1\colon Y\times_BB^I\to Y$ (which is pulled back from the trivial fibration $d_1\colon B^I\to B$) yields $\id_Y$. Hence, we may factor it as a trivial cofibration followed by a trivial fibration $Y\tox{i}Z\tox{\pi}Y\times_BB^I$; we consider $Z$ an object over $B\times B$ via $(\id_B\times y)\pi$. Next, suppose we have a morphism $\nu\colon Z\to B\times Y$ over $B\times B$ such that $\gamma(\nu)\gamma(\pi)\I=\nat_B^Q$ (we will find such a $\nu$ below).

Now choose a lift $z\colon X\to Z$ of $\br{p,kx}$ along $\pi$. Since $(X,A,x)$ is cofibrant and $(B\times Y,B\times B,\id\times y)$ fibrant, the commutativity of the diagram on the left above is then equivalent to the commutativity of the diagram on the right up to homotopy in $\C^\to$ . As remarked at the end of Definition~\ref{defn:homotopy-over}, this is equivalent to $\nu z$ and $\br{fx,q}$ being homotopic over the trivial homotopy from $\br{f,g}$ to itself. Since we must have $\pi_1\nu z=fx$, this is in turn equivalent to $\pi_2\nu z$ and $q$ being homotopic over the trivial homotopy from $g$ to itself.

We conclude that the left side of the equivalence in the theorem statement is equivalent to $\pi_2\nu z$ and $q$ being homotopic over the trivial homotopy from $g$ to itself.

We still need to define $\nu$. To do so, first consider a pullback $(Y\times Y)\times_{B\times B}B^I$ with respect to $y\times y$ and $\br{d_1,d_2}$, and factor $\br{\Delta_Y,sy}\colon Y\to (Y\times Y)\times_{B\times B}B^I$ as a trivial cofibration and fibration $Y\tox{s}Y^I\tox{\br{\br{d_1,d_2},y^I}}(Y\times Y)\times_{B\times B}B^I$ (thus producing the situation in the last diagram of (\ref{eq:homotopy-over-diag}) in Definition~\ref{defn:homotopy-over}) and choose a diagonal filler $j$ in the following square.
\[
  \begin{tikzcd}
    Y\ar[d, "i"']\ar[r, "s"]&Y^I\ar[d, "\br{d_1,y^I}"]\\
    Z\ar[r, "\pi"]\ar[ru, "j", dashed]&Y\times_BB^I
  \end{tikzcd}
\]
Now define $\nu=\br{yd_1,d_2}j\colon Z\to B\times Y$ (note that this is indeed a morphism over $B\times B$). Let us verify that it satisfies the required condition $\gamma(\nu)\gamma(\pi)\I=\nat_B^Q$. By definition of $\nat_B^Q$, this means that the triangle below on the left should commute in $\Ho(\C^\to)_\fb$, or equivalently, that the isomorphic diagram in the center should commute. But in fact, the diagram on the right in $\C^\to$ already commutes before applying $\gamma$.
\[
  \begin{tikzcd}
    &[10pt]Y\times_BB^I\ar[d, "\gamma(\nu)\gamma(\pi)"]\\
    Y\ar[r, "\gamma(\br{y,\id_Y})"']\ar[ru, "\gamma(\br{y,sy})"]&B\times{}Y
    \\[-5pt]
    A\ar[r, "\br{f,g}"]&B\times{}B
  \end{tikzcd}
  \quad\quad
  \begin{tikzcd}
    &[14pt]Z\ar[d, "\gamma(\nu)"]\\
    Y\ar[r, "\gamma(\br{y,\id_Y})"']\ar[ru, "\gamma(i)"]&B\times{}Y
    \\[-5pt]
    A\ar[r, "\br{f,g}"]&B\times{}B
  \end{tikzcd}
  \quad\quad
  \begin{tikzcd}
    &Z\ar[d, "\nu"]\\
    Y\ar[r, "\br{y,\id_Y}"']\ar[ru, "i"]&B\times{}Y
    \\[-5pt]
    A\ar[r, "\br{f,g}"]&B\times{}B
  \end{tikzcd}
\]

Now, since (by definition) $\pi z=\br{p,kx}\colon X\to Y\times_BB^I$, the composite $jz\colon X\to Y^I$ gives a homotopy from $p$ to $d_2jz=\pi_2\nu z$ over $k$. Hence, by Proposition~\ref{prop:homotopy-over-left-right}, we conclude that $p$ and $q$ are homotopic over $k$
if and only if $q$ and $\pi_2\nu z$ are homotopic over the trivial homotopy on $g$, as desired.
\qed

\subsubsection{}\thm\label{thm:special-invariance}
Let $\C$ be a suitable model category and let $M,N\colon{}\sigma\to\C_\fb$ be $\sigma$-interpretations, and suppose we have a homotopy equivalence $\alpha\colon M\to N$.

Then, given homotopical $\sigma$-interpretations $\widehat M$ and $\widehat N$ over $M$ and $N$, there exists, for each formula-in-context $(\phi,\vec{x})$ over $\sigma$, a homotopy equivalence $\widehat M_{\vec{x}}(\phi)\to\widehat N_{\vec{x}}(\phi)$ over the homotopy equivalence $\alpha_{\tp\vec{x}}\colon M(\vec{x})\to N(\vec{x})$.
\pf
As remarked in Definition~\ref{defn:2cat-homomorphism}, a homotopy equivalence of $\sigma$-interpretations in $\C_\fb=\C_\cfb$ is the same thing as a pseudoequivalence with respect to the Quillen 2-category structure. Moreover, it follows from Theorem~\ref{thm:2-cells-in-hofb-tot} that there is a homotopy equivalence $\widehat M_{\vec{x}}(\phi)\to\widehat N_{\vec{x}}(\phi)$ over the homotopy equivalence $\alpha_{\tp\vec{x}}$ if and only if there is an equivalence $\widehat M_{\vec{x}}(\phi)\to\widehat N_{\vec{x}}(\phi)$ over the equivalence $\alpha_{\tp\vec{x}}$ in the 1D2F $\fib{HoF_\fb(\C_\fb)}$.

Since, by Proposition~\ref{prop:htpical-interp-unique}, $\gamma\circ\widehat M$ and $\gamma\circ\widehat N$ are $\sigma$-interpretations in $\fib{HoF_\fb(\C_\fb)}$ over $M$ and $N$, the claim now follows from Theorem~\ref{thm:somewhat-abstract}.
\qed

\section{Examples and further questions}\label{sec:examples-and-questions}
We now give some examples of sentences and their interpretation under the homotopical semantics. In each case, we fix some algebraic signature $\sigma$ and a $\sigma$-interpretation $M$ in $\Kan$ (or in $\Top_\cf$ -- the category of topological spaces which are homotopy equivalent to a CW complex -- and then consider the interpretation $\Sing\circ M$, as discussed in \S\ref{subsec:top-spaces}). We then take some first-order sentence over this signature and see what it means for it to be satisfied in $M$ (or in the case of $\Top_\cf$, in $\Sing\circ M$) -- i.e., for its interpretation to be a non-empty space with respect to any homotopical interpretation $\widehat M$ over $M$ (or $\Sing\circ M$).

Though the notion of satisfaction in $M$ is independent of the choice of homotopical interpretation $\widehat M$, it will be convenient to use this freedom to choose a particular such $\widehat M$, where we recall that (up to isomorphism) or only freedom is in the choice of equality objects; we take these to be the ``standard'' path-spaces $X^I\to X\times X$ in $\sSet$, where $I$ is the simplicial interval $\Delta^1$.

After the examples, we consider some further questions regarding the material of this paper.

\subsection{Examples of interpretations of sentences}\label{subsec:examples}

\subsubsection{}\sectitle{Contractibility}\label{examples-contractibility}
First, we consider the signature $\sigma$ consisting of a single sort $A$ and having no operation symbols, and the sentence in this language
\[
  \exists x\ \forall y\ (x=y).
\]
We claim that this is interpreted under the semantics as ``$A$ is contractible''.

Fix a $\sigma$-interpretation in $\Kan$, i.e., a Kan complex $X$.

Now, the formula $x=y$ in the context $\seq{x,y}$ is interpreted as the path-space fibration $X^I\tox{\br{d_1,d_2}}X\times{}X$. Next, the formula $\forall{}y\,(x=y)$ is interpreted as an indexed product $\prod_{\pi_1}(X^I,\br{d_2,d_2})$. Finally, $\exists{}x\,\forall{}y{}\,(x=y)$ is interpreted (as always, up to isomorphism) as the domain of $\prod_{\pi_1}(X^I,\br{d_2,d_2})$.

Hence, we are interested in when the domain of $\prod_{\pi_1}(X^I,\br{d_1,d_2})$ is non-empty. This will hold if and only if there is a morphism $(\tm_\sSet,x)\to{}\prod_{\pi_1}(X^I,\br{d_1,d_2})$ in $\sSet/X$ for some $x\colon\tm_\sSet\to{}X$. By the adjunction (Remark~\ref{rmk:adjoints}), this is equivalent to having a morphism from $\pi_1^*(\tm_\sSet,x)\cong(X,\br{x!,\id_X})$ to $(X^I,\br{d_1,d_2})$ in $\C/X\times{}X$:
\[
  \begin{tikzcd}
    &X^I\ar[d, "\br{d_1,d_2}"]\\
    X\ar[r, "\br{x!,\id_X}"']\ar[ru, dotted]&X\times{X}.
  \end{tikzcd}
\]
But this is by definition a (right-)homotopy between $X$ and the constant map $x!$, i.e. a contraction of $X$ onto $x$.

If we start with a topological space $X\in\Top_\cf$ instead of a Kan complex, then the above shows that $X$ satisfies the sentence in question if and only if the singular simplicial set of $X$ is contractible. But as is well-known, this holds if and only if $X$ itself is contractible.

\subsubsection{}\sectitle{Homotopies}\label{subsubsec:examples-homotopies}
Now let $\sigma$ be the signature consisting of two sorts $A,B$ and two function symbols $f,g\colon{}A\to{B}$. We consider the sentence
\[
  \forall{x\in{A}}\ (f(x)=g(x)).
\]
We claim that this is interpreted as ``$f$ is homotopic to $g$''.

Suppose we have a $\sigma$-interpretation in $\Kan$; that is, two Kan complexes $X,Y$, and two morphisms $f,g\colon{}X\to{Y}$. The formula $y_1=y_2$ (in the context $\br{y_1,y_2}$) will (again) be interpreted as the path-space fibration $Y^I\tox{\br{d_1,d_2}}Y\times{}Y$. Now, we have the morphism $\br{f,g}\colon{}X\to{Y\times{Y}}$, and $f(x)=g(x)$ (in the context $\br{x}$) will be interpreted as $\br{f,g}^*(Y^I,\br{d_1,d_2})$. Finally, the above sentence will be interpreted as $\prod_!\br{f,g}^*(Y^I,\br{d_1,d_2})$, the points of which are (by the adjunction) in bijection with the sections of $\br{f,g}^*(Y^I,\br{d_1,d_2})$, which are in turn in bijection with the lifts
\[
  \begin{tikzcd}
    &Y^I\ar[d, "\br{d_1,d_2}"]\\
    X\ar[r, "\br{f,g}"']\ar[ru, dotted]&Y\times{Y}
  \end{tikzcd}
\]
which are, of course, by definition (right-)homotopies $f\sim{g}$.

For a $\sigma$-interpretation in $\Top_\cf$, i.e., a pair of maps $f,g\colon{}X\to{}Y$, we thus see that the above sentence is satisfied if and only $\Sing(f)$ and $\Sing(g)$ are homotopic and, again, this is the case if and only if $f$ and $g$ are homotopic.

Now considering a signature with two sorts $A,B$ and two function symbols $f\colon{}A\to{B}$ and $g\colon{}B\to{A}$, we have by the same reasoning as above that
\[
  \forall{x\in{A}}\ (g(f(x))=x)\,\wedge\,
  \forall{y\in{B}}\ (f(g(y))=y)
\]
is interpreted (in both $\Top_\cf$ and $\Kan$) as ``$f$ and $g$ constitute a homotopy equivalence'' (i.e., both composites are homotopic to the identity).

Similarly, for the signature consisting of a single sort $A$ and binary function symbol $f\colon{}A\times{A}\to{A}$,
\[
  \forall{}x,y,z\in{}A\ [f(f(x,y),z)=f(x,f(y,z))]
\]
is interpreted as ``$f$ is homotopy-associative''.

\subsubsection{}\sectitle{Classical logic}
We now give an example showing that the homotopical semantics are not sound with respect to classical logic. By this we mean that there is a formula of the form $\neg\neg{P}\To{P}$ over some signature $\sigma$ and a $\sigma$-interpretation (in $\Top_\cf$ and $\Kan$) under which the interpretation of this formula is empty.

First, we note that it is important that $P$ is not a closed formula. Indeed, a closed formula is interpreted as a Kan complex $X$, and its negation is interpreted as an empty Kan complex or one-point Kan complex according to whether $X$ is non-empty or empty. From this it follows that the interpretation of $\neg\neg{P}\To{P}$ (and similarly $P\vee\neg{P}$) is always non-empty. This circumstance is familiar, for example, from Kleene's realizability semantics for intuitionistic arithmetic.

Now, for our example, we consider, in the signature $\sigma$ consisting of a single sort $A$ and no function symbols, the sentence
\begin{equation}\label{eq:classical-counterexample}
  (\exists x
  \ \forall y
  \ (\neg\neg{x=y}))
  \To
  (\exists x
  \ \forall y
  \ (x=y)).
\end{equation}
Given a structure $X$ (in $\Kan$ or $\Top_\cf$) for $\sigma$ we have already seen that the interpretation of the right side of this implication is inhabited if and only if $X$ is contractible. Let us consider the left side. We first consider the case of $X$ in $\Kan$. We will show that this sentence is satisfied if and only if $X$ is non-empty and path-connected (i.e., for any vertices $x,y\in{X}_0$ there is an edge $e\in{X}_1$ from $x$ to $y$).

We have, again, that $x=y$ is interpreted as the path space $(X^I,\br{d_1,d_2})$. We recall that $\neg\neg x=y$ is an abbreviation of $(x=y\To\bot)\to\bot$. Here, $\bot$ is interpreted as the initial Kan fibration $(\emptyset,\text{!`},X)$.

Now, it easy to see that in any category (such as $\Kan/X$) with a strong initial object $0$ (i.e. every morphism with codomain $0$ is an isomorphism), any exponential object $A\To0$ is a subsingleton (i.e, the morphism $!_{A\To0}$ to the terminal object is a monomorphism). Since there exists a morphism $A\to((A\To0)\To0)$, it follows that the unique morphism $A\to1$ factors through $(A\To0)\To0$.

In the case of a Kan fibration $(E,e)$ in $\Kan/X$, this tells us that $\neg\neg({E},e)$ is a monomorphism into $X$ whose image contains the image of $e$. In particular, if $E$ is surjective onto $X$, then $\neg\neg{(E,e)}$ is an isomorphism.

Now, if $X$ is path-connected, then the path space $X^I$ is clearly surjective on vertices. But as an easy inductive argument shows, any Kan fibration which is surjective on vertices is surjective. Hence, for $X$ path-connected, $\neg\neg{x=y}$ is interpreted as an isomorphism, whence it follows that, for $X$ non-empty and path-connected, the following sentence is satisfied.
\begin{equation}\label{eq:path-conn-sentence}
  \exists x
  \ \forall y
  \ (\neg\neg{x=y}).
\end{equation}
For the other direction, it suffices to see that if a Kan fibration $(E,e)$ in $\Kan/X$ is not surjective, then neither is $\neg\neg(E,e)$, since if $\neg\neg{}x=y$ is interpreted as a non-surjective morphism, the interpretation of (\ref{eq:path-conn-sentence}) must be empty. Suppose $e$ is not surjective and let $p\in{X_0}$ be a vertex not in the image of $e$. Then the minimal sub-simplicial set $\br{p}$ of $X$ containing $p$ is disjoint from the image of $e$. Hence $(\br{p},i)\wedge(E,e)\cong\init$, where $i\colon\br{p}\to{}E$ is the inclusion, so we have a morphism $(\br{p},i)\to\neg({E},e)$. In particular, $(\br{p},i)\wedge\neg({E},e)$ is non-empty, so there cannot be a morphism $(\br{p},i)\to\neg\neg({E},e)$.

Note that by the same kind of argument as in the previous examples, we also have that the interpretation of the sentence (\ref{eq:path-conn-sentence}) in $\Top_\cf$ is ``$A$ is non-empty path-connected''.

Now suppose $X$ is a Kan complex or a topological space in $\Top_\cf$ which is path-connected but not contractible (for instance, the circle). Then the of (\ref{eq:classical-counterexample}) is interpreted as a non-empty Kan complex, whereas the conclusion is interpreted as the empty Kan complex. Hence the implication is empty.

Finally, we note that this implies that (the universal closure of)
\[
  \neg\neg{x=y}\To{x=y}
\]
cannot be satisfied for such an $X$ since this would imply (by the soundness of the interpretation with respect to intuitionistic logic) that (\ref{eq:classical-counterexample}) would also be satisfied.

\subsubsection{}\sectitle{Homotopy-equivalence}
We saw above that we can easily express that two morphisms constitute a homotopy equivalence, in the same way as we would classically express that they constitute a bijection. We can also classically express that a single function $f\colon{}A\to{B}$ is a bijection by
\[
  \forall{b\in{B}}\ \exists{a\in{A}}\ (
  fa=b\wedge
  \forall{a'\in{A}}\ (fa'=b\To{a'=a})
  ).
\]
Let us see that the interpretation of this sentence (in $\Kan$, and hence in $\Top_\cf$) is non-empty if and only if the interpretation of $f$ is a homotopy equivalence.

First of all, the sentence is equivalent to the conjunction of
\[
  \forall{b\in{B}}\ \exists{a\in{A}}\ (
  fa=b)
  \quad\text{and}\quad
  \forall{b\in{B}}\ \exists{a\in{A}}
  \ \forall{a'\in{A}}\ (fa'=b\To{a'=a}
  ).
\]
The same reasoning as in \S\ref{subsubsec:examples-homotopies} shows that the first part is satisfied by a map $f\colon{}X\to{}Y$ if and only if there exists a map $g\colon{}Y\to{}X$ such that $f\circ{}g$ is homotopic to $\id_Y$.

The second part, with the quantifiers removed, is interpreted as a certain fibration over $Y\times{}X\times{}X$. By making use of the relevant adjunctions, we can see that the space which is the interpretation of the quantified sentence is inhabited if and only if there exists a map $g\colon{}Y\to{}X$ and a dotted lift in the following diagram.
\[
  \begin{tikzcd}
    &Y\times{}X^I\ar[d, "\br{\pi_1,d_2\pi_2,d_1\pi_2}"]\\
    (Y\times{}X)\times_{Y\times{}X\times{}X}(X\times(X\times_Y{Y^I}))
    \ar[r]\ar[ru, "k", dashed]&Y\times{}X\times{}X
  \end{tikzcd}
\]
Here, in $X\times_YY^I$, $X$ is mapping to $Y$ via $f$ and $Y^I$ is mapping to $Y$ via $d_1$; and in the object on the bottom-left of the diagram, $Y\times{}X$ is mapping to $Y\times{}X\times{}X$ via $\br{\pi_1,g\pi_1,\pi_2}$, and $X\times(X\times_YY^I)$ is mapping to $Y\times{}X\times{}X$ via $\br{d_2\pi_2\pi_2,\pi_1,\pi_1\pi_2}$.

We claim that such a lift exists if and only if $g\circ{}f$ is homotopic to $\id_X$.

In one direction, we have a map $q\colon{}X\to{}(Y\times{}X)\times_{Y\times{}X\times{}X}(X\times(X\times_Y{Y^I}))$ which is given by $\br{\br{f,\id_X},\br{gf,\br{\id_X,sf}}}$ (where $s$ is, as usual, the canonical map $Y\to{}Y^I$). Hence, given a lift $k$ as above, the composite $\pi_2kq$ gives a homotopy $X\to{}X^I$ from $\id_X$ to $gf$.

In the other direction, suppose we have a homotopy $h\colon{}X\to{}X^I$ from $\id_X$ to $gf$. We then define $k$ as $\br{\pi_1\pi_1,h'}$, where $h'$ is the composite\vspace{-8pt}
\[
  (Y\times{}X)\times_{Y\times{}X\times{}X}(X\times(X\times_Y{Y^I}))
  \tox{\br{h\pi_2\pi_1,g^I\pi_2\pi_2\pi_2}}
  X^I\times_XX^I\to{}X^I
\]
in which the second map is composition of paths.

\subsection{Further problems and questions}\label{subsec:questions}
We mention some possible further directions.

\emph{Completeness}. This is probably the most natural question to ask about the homotopical semantics: are they a complete semantics for intuitionistic logic? I.e., is it the case that, if a sentence $\phi$ over a signature $\sigma$ is satisfied by every $\sigma$-structure in $\Kan$, then $\phi$ is intuitionistically provable?

And if it is not complete for intuitionistic logic, then (as was pointed out to me by Thomas Icard) we can still try to characterize the ``intermediate'' logic between intuitionistic and classical for which it \emph{is} complete.

\emph{Limited expressivity}. We mentioned at the end of \S\ref{subsec:homotopical-structures} that first-order homotopical logic is much less expressive than homotopy type theory. However, we have not \emph{proven} that any particular property is inexpressible, and it would be interesting to do so -- for example, to prove that over the trivial signature with a single sort $A$ and no function symbols, there is no sentence $\phi$ satisfied by exactly those spaces $X$ which are simply connected. Or to give another example, there should be no sentence over the signature consisting of one sort and a single binary operation which is satisfied exactly by those operations satisfying the ``$A_4$'' (or ``Stasheff pentagon'') condition.

\emph{Higher-dimensional generalizations}. Given this lack of expressiveness, it is natural to seek extensions of first-order logic which increase the expressiveness. For example, one could add some, but not all, of what is present in type theory -- say, an additional sort $s=_At$ for any two terms $s$ and $t$ of sort $A$, so that one could express that ``two homotopies are homotopic'': $e=_{s=_At}e'$.

One would hope to have a nice categorical formulation of the corresponding semantics, as we have for first-order homotopical logic. Indeed, it is also natural to seek ``higher-dimensional'' generalizations of the fibrational semantics. For example, the fact that we can only express ``one level'' of homotopies in the language seems to correspond to the fact that the fibrations we are considering are (``only'') two-dimensional. On the ``semantic'' side, there are natural higher-dimensional categories close at hand -- for example, instead of having the fibers of $\fib{HoF(\C)}$ be the \emph{homotopy categories} of the slices $\C/A$, one could try to take the corresponding \emph{$\infty$-categories} (or some truncation thereof). We might then seek a higher-dimensional analogue of the syntactic fibration, morphisms out of which would give the semantics for such ``higher-dimensional'' extensions of first-order logic.
\edef\savingsubsection{\thesubsection}
\appendix

\setcounter{subsection}{\savingsubsection}
\section{The syntactic fibration}\label{sec:appendix}
The main purpose of this appendix is to construct a free $h^=$-fibration $\fib{C}$ over any f.p.\ category $\B$. When $\B$ is itself a free f.p.\ category on some signature $\sigma$, this is the ``syntactic'' fibration associated to $\sigma$, such that the morphisms out of it into some fibration $\fib{C'}$ correspond to interpretations of $\sigma$ in $\fib{C'}$.

The construction of $\fib{C}$ closely parallels the analogous construction for propositional logic from \cite{makkai-harnik-lauchli}, to which we refer for a thorough discussion of and motivation behind the construction.

In \emph{op. cit.}, the construction (which is of a category and not a fibration) is directly ``syntactic''; the objects of the category are propositional formulas, and the morphisms are equivalence classes of logical deductions.\footnote{Though note that their construction is with respect to an arbitrary \emph{theory} over the given propositional language, whereas ours only covers the ``empty theory''.} In an earlier version of this paper, we similarly provided a ``syntactic'' construction of $\fib{C}$ (which in particular could only be carried out in the case that $\B$ is free), but in the present version, we provide a more direct construction, in line with the usual construction of free structures.

We note that in \emph{op. cit}, the ``direct'' construction and the ``syntactic'' construction are one and the same. This is because of the observation due to Lambek that the inference rules of intuitionistic propositional logic precisely correspond to the universal properties of the operations in a bicartesian closed category.

In the case of intuitionistic first-order logic, there is an analogous situation (this being an observation of Lawvere), but with $h^=$-fibrations instead of bicartesian closed categories. The reason that the two constructions of the free $h^=$-fibration are no longer identical is that there is a slight discrepancy between the syntax of first-order logic (as it is usually construed) and the operations in a fibration; in the latter, ``substitution'' (i.e., pullbacks) is a ``primitive'' operation, whereas in the syntax it is not, and also, the fibrational and syntactic treatments of equalities and quantifiers differs slightly as well.

This complicates considerably the task of showing that the ``syntactically'' constructed fibration is in fact an $h^=$-fibration, and that it is free. This -- in addition to the advantage of it working over arbitrary $\B$ -- is why we opt here for the direct construction.

On the other hand, the advantage of the syntactic construction is that it makes it clear that the morphisms are precisely equivalence classes of deductions, and this fact is obscured in the present construction. It is still ``almost'' obvious from the construction (the only problem being the discrepancies mentioned above between the syntax and the $h^=$-fibration operations), but to actually prove it is non-trivial, and we do not attempt to do so, though we discuss the matter in \S\ref{subsec:free-over-free}.

We will, however, show the weaker result that, for two objects $P$ and $Q$ in the free $h^=$-fibration corresponding to given formulas $\phi$ and $\psi$, there exists a morphism $P\to Q$ if and only if $\phi\To\psi$ is intuitionistically provable.

\subsection{The free $h^=$-fibration over a f.p.\ category $\B$}
We now construct the free $h^=$-fibration over an arbitrary finite product category $\B$. As we have mentioned, this notion of freeness is given by a certain universal property which determines the fibration up to equivalence over $\B$.

However, the construction will proceed by producing (in the standard way) a structure which is free in the stricter (``0-categorical'') sense (thus determining it up to isomorphism -- though we will not make use of this fact).

To perform the latter construction, we will have to specify what kind of free structure we want to construct, and for this we will deploy the notion of multi-sorted algebraic signature in a different way than in the rest of the paper (but in the same way as in \cite{makkai-harnik-lauchli}); in particular, in this case, we will be concerned with ordinary set-based structures for the signature, rather than ones valued in simplicial sets or in other categories.

\subsubsection{}\label{subsubsec:initial-structures}
We recall the standard construction of an initial structure (or free structure with no generators) on a multi-sorted algebraic signature $\sigma$, which we will use twice in what follows.\footnote{We note that one can also prove the existence of such structures using the adjoint functor theorem, but we will need to make use of the details of the explicit construction.}

Recall from Definition~\ref{defn:syntax} that the set of terms over $\sigma$ was defined relative to a fixed infinite set $\Varn$ of variable names. In fact, the set thus constructed -- or better, the $\Ob\sigma$-indexed family of sets -- is precisely the \emph{free $\sigma$-structure on the ($\Ob\sigma$-indexed) set $\Var_{\sigma}=\Ob\sigma\times\Varn$}. If we instead take $\Varn=\emptyset$ -- or, what amounts to the same thing -- if we take only the \emph{closed} $\sigma$-terms (those not having any free variables) -- we obtain the definition of the initial $\sigma$-structure. This structure is characterized up to isomorphism by an obvious universal property -- this is precisely the ``principal of structural recursion'' mentioned in Definition~\ref{defn:syntax}. As was said there, we will freely make use of this -- as well as the corresponding principle of structural induction -- when dealing with the initial $\sigma$-structure.

We will also need to make use of the \emph{initial model of an algebraic theory} -- i.e., of an algebraic signature $\sigma$ together with a set of \emph{identities}, each identity being given by a pair of $\sigma$-terms of the same sort. Concretely, this initial model is obtained as the quotient of the initial $\sigma$-structure by a certain equivalence relation: to begin with, we have a relation $R$ consisting of all pairs of closed terms obtained by substituting arbitrary closed terms for the variables in each of the given identities. We then take the least congruence relation containing $R$, i.e., the least equivalence relation $\sim$ satisfying $ft_1\ldots t_n\sim ft_1'\ldots t_n'$ whenever $t_1\sim t_1',\ldots,t_n\sim t_n'$ for any function symbol $f$ of $\sigma$.

Again, this structure has a universal property, now with respect to all $\sigma$-structures satisfying the given identities.

\subsubsection{}\constr\label{constr:free-hfib}
Let $\B$ be a finite product category. We will define the free $h^=$-fibration over $\B$ in terms of initial structures on certain signatures associated to $\B$.

To begin with, we describe a certain multi-sorted algebraic signature $\sigma_\B$ associated to $\B$. For the set of sorts, we set $\Ob\sigma_\B:=\Ob\B$, and it has function symbols:
\begin{enumerate}[(i), topsep=0pt]
\item $\top_A,\bot_A\in\sigma_\B(\emptyset,A)$ for each $A\in\B$ and $\Eq_\Delta\in\sigma_\B(\emptyset,B)$ of sort $B$ for each diagonal morphism $\Delta\colon A\to B$ in $\B$
\item $\wedge_A,\vee_A,\To_A\in\sigma_\B(\seq{A,A},A)$ for each $A\in\B$
\item $\prod_\pi,\sum_\pi\in\sigma_\B(A,B)$ for each product projection $\pi\colon A\to B$ in $\B$ and $f^*\in\sigma_\B(B,A)$ for each morphism $f\colon A\to B$ in $\B$.
\end{enumerate}
Let $\Ob\C$ be the initial structure on the above signature, as described in \S\ref{subsubsec:initial-structures} and, for $A\in\B$, denote by $\Ob\fib{C}^A$ its set of elements of sort $A$ (the notation will be justified soon).

Next, We describe a second algebraic signature $\sigma_\B'$ whose set $\Ob\sigma_{B}'$ of sorts is the set of triples $(P,Q,f)$, with $P\in\fib{C}^A$ and $Q\in\fib{C}^B$ for some $A,B\in\B$, and with $f\colon A\to B$ a morphism in $\B$.

The function symbols of this signature are given schematically below (compare \cite[p.~210]{makkai-harnik-lauchli}). Each figure indicates a set of operations, one for each of the possible values of the relevant parameters $A,B,C,f,g,P,Q,R,S,T,\Delta,A',B',\pi,\pi',h,k$. The range of these parameters is given as follows: $A,B,C\in\Ob\B$; $f\colon A\to B$ and $g\colon B\to C$ are morphisms in $B$; $P,Q,R\in\Ob\fib{C}^A$, $S\in\Ob\fib{C}^B$, and $T\in\Ob\fib{C}^C$; $\Delta\colon A\to B$ is a diagonal morphism; and finally, $A',B',\pi,\pi',h,k$ are objects and morphisms in $\B$ forming a pullback square
\[
  \begin{tikzcd}
    A'\ar[r, "\pi'"]\ar[d, "h"']\ar[dr, "\lrcorner", pos=0, phantom]&B'\ar[d, "k"]\\
    A\ar[r, "\pi"]&B
  \end{tikzcd}
\]
with $\pi$ (and hence $\pi'$) a product projection.

Each figure below displays, above the horizontal line, the arity of the operation, and below, the codomain sort, and also introduces a notation for the function symbol (i.e., recalling that normally, for a function symbol $f$ and terms $t_1,\ldots,t_k$ of the appropriate sorts, we write $ft_1\ldots t_k$ for the resulting term, we introduce an alternative notation for $ft_1\ldots t_k$). To indicate that a term $p$ is of sort $(P,Q,f)$, we write $p\colon P\tox[f]{}Q$. If the subscript $f$ is omitted, it is assumed to be $\id_A$.
{
\newcommand\gap{\hspace{20pt}}
\newcommand\gop{\hspace{40pt}}
\newcommand\blk{\vphantom{P}}
\newcommand\inv{\mathrm{inv}}

{\bf Category and fibration structure:}\label{fig:first-hfib-table}
\[
  \frac{\blk}{1_P\colon{}P\tox{}P}\gap
  \frac{p\colon{}P\tox[f]{}S\quad{}q\colon{}S\tox[g]{}T}{q\circ{}p\colon{}P\tox[gf]{}T}\gap
  \frac{\blk}{\crt{f}{S}\colon f^*S\tox[f]{} S}\gap
  \frac{p\colon{}P\tox[gf]{}T}{\cind{p}^{f,g}\colon{}P\tox[f]{}g^*T}
\]

{\bf Finite products and coproducts:}
\[
  \frac{}{!_P\colon{}P\to{}\top_{A}}\gop
  \frac{\blk}{\ex_P\colon\bot_{A}\to{}P}
\]\[
  \frac{\blk}{\pi_{PQ}\colon{}P\wedge{}Q\to{}P}\gap
  \frac{\blk}{\pi'_{PQ}\colon{}P\wedge{}Q\to{}Q}\gap
  \frac{p\colon{}P\to{Q}\quad{}q\colon{}P\to{R}}{\br{p,q}\colon{}P\to{Q\wedge{R}}}
\]\[
  \frac{\blk}{\kappa_{PQ}\colon{}P\to{}P\vee{}Q}\gap
  \frac{\blk}{\kappa'_{PQ}\colon{}Q\to{}P\vee{}Q}\gap
  \frac{p\colon{}P\to{R}\quad{}q\colon{}Q\to{R}}{\pl[]{p,q}\colon{}P\vee{}Q\to{R}}
\]

{\bf Exponentials:}
\[
  \frac{\blk}{\varepsilon_{PQ}\colon(P\To{}Q)\wedge{}P\to{}Q}\gap
  \frac{p\colon{}P\wedge{}Q\to{R}}{p^\simm\colon{}P\to{}Q\To{R}}
\]

{\bf Indexed products and sums:}
\[
  \frac{}
  {\varepsilon^\Pi_{\pi,P}\colon\pi^*\prod_\pi P\tox[\id_A]{}P}\gap
  \frac{p\colon\pi^*S\tox[\id_A]{}P}{\lambda{}p\colon{}S\tox[\id_B]{}\prod_\pi P}
  \gop
  \frac{}
  {\cort{\pi}{P}\colon{}P\tox[\pi]{}\sum_\pi{}P}\gap
  \frac{p\colon{}P\tox[\pi]{}S}{\underline p\colon\sum_\pi P\tox[\id_B]{}S}
\]

{\bf Equality objects:}
\[
  \frac{\blk}{\rho_\Delta\colon\top_{A}\tox[\Delta]{}\Eq_{\Delta}}\gap
  \frac{f\colon\top_A\tox[\Delta]{}S}
  {\xi{}f\colon\Eq_{\Delta}\tox[\id_B]{}S}
\]

{\bf Stability of the operations under pullbacks}
\[
  \frac{\blk}{\inv^\top_A\colon \top_{A'}\tox[\id_{A'}]{}h^*\top_A}\gop
  \frac{\blk}{\inv^\bot_A\colon h^*\bot_A\tox[\id_{A'}]{}\bot_{A'}}
\]\[
  \frac{\blk}{\inv^\wedge_{PQ}\colon h^*P\wedge h^*Q\tox[\id_{A'}]{}h^*(P\wedge Q)}\gop
  \frac{\blk}{\inv^\vee_{PQ}\colon h^*(P\vee Q)\tox[\id_{A'}]{}h^*P\vee h^* Q}
\]\[
  \frac{\blk}{\inv^\To_{PQ}\colon h^*P\To h^*Q\tox[\id_{A'}]{}h^*(P\To Q)}
\]\[
  \frac{\blk}{\inv^\Sigma_{\pi\pi'hk}\colon k^*\sum_\pi P\to\sum_{\pi'}h^*P}\gop
  \frac{\blk}{\inv^\Pi_{\pi\pi'hk}\colon\prod_{\pi'}h^*P\to k^*\prod_\pi P}
\]
}

Next, we will define certain identities defined over $\sigma_\B'$.

Each figure below represents a set of identities, one (or more) for each possible value of the relevant parameters $A,B,C,f,g,P,Q,R,S,T,\Delta,A',B',\pi,\pi',h,k,D,g',T'$. The range of these parameters is as above, and in addition, $D\in\Ob\B$, $g'$ is a morphism $C\to D$, and $T'\in\Ob\fib{C}^D$.

We introduce certain abbreviations: under ``Exponentials'', we write $x\wedge y$ for $\br{x\circ\pi_{PQ},y\circ\pi'_{PQ}}$; under ``Stability of the operations under pullbacks'', for an expression $x\colon X\tox[\id_A]{}Y$, we write $h^*x$ for $\cind{x\circ\crt{h}Y}^{\id_{A'},h}$, and similarly with $\pi^*x$ under ``Indexed products and sums''; finally, under ``Stability of the operations under pullbacks'', we write ``$X=Y\I$'' as a shorthand for the two identities $Y\circ X=\id_W$ and $X\circ Y=\id_Z$ with appropriate $W$ and $Z$.
{
\newcommand\gap{\hspace{20pt}}
\newcommand\gop{\hspace{40pt}}
\newcommand\blk{\vphantom{P}}
\newcommand\inv{\mathrm{inv}}

{\bf Category:}
\[
  \frac{p\colon{}P\tox[f]{}S\gap{}q\colon{}S\tox[g]{}T\gap{}r\colon{}T\tox[h]{}T'}
  {p\circ{1_P}=p\gap{}1_S\circ{}p=p\gap
  {(r\circ{}q)\circ{}p=r\circ(q\circ{}p)}}
\]

{\bf Fibration:}
\[
  \frac{p\colon{}P\tox[gf]{}T}
  {\crt{g}T\circ \cind{p}^{f,g}=p}\gop
  \frac{p\colon{}P\tox[f]{}g^*T}
  {\cind{\crt{g}T\circ p}^{f,g}=p}
\]

{\bf Finite products and coproducts:}
\[
  \frac{p\colon{}P\to\top_{A}}
  {p=!_P}\gop
  \frac{p\colon\bot_{A}\to{}P}
  {p=\ex_P}
\]\[
  \frac{p\colon{}P\to{}Q\gap{}q\colon{}P\to{}R}
  {\pi_{QR}\circ\br{p,q}=p\gap
  \pi'_{QR}\circ\br{p,q}=q}\gap
  \frac{p\colon{}P\to{}Q\wedge{}R}
  {\seq{{\pi_{QR}\circ{}p,\pi'_{QR}\circ{}p}}=p}
\]\[
  \frac{p\colon{}P\to{}R\gap{}q\colon{}Q\to{}R}
  {[p,q]\circ\kappa_{PQ}=p\gap
  [p,q]\circ\kappa'_{PQ}=q}\gap
  \frac{p\colon{}P\vee{}Q\to{}R}
  {[p\circ\kappa_{PQ},p\circ\kappa'_{PQ}]=p}
\]

{\bf Exponentials:}
\[
  \frac{p\colon{}P\wedge{}Q\to{}R}
  {\varepsilon_{QR}\circ{}(p^\sim\wedge{}1_Q)=p}\gop
  \frac{p\colon{}P\to{}(Q\Rightarrow{}R)}
  {(\varepsilon_{QR}\circ(p\wedge1_Q))^\sim=p}
\]
{\bf Indexed products and sums:}
\[
  \frac{p\colon\pi^*S\tox[\id_A]{}P}
  {\varepsilon^\Pi_{\pi,P}\circ\pi^*\lambda{}p=p}\gap
  \frac{p\colon{}S\tox[\id_B]{}\prod_\pi{}P}
  {\lambda(\varepsilon^\Pi_{\pi,P}\circ\pi^*p)=p}
\gop
  \frac{p\colon{}P\tox[\pi]{}S}
  {\underline{p}\circ\cort{\pi}P=p}\gap
  \frac{p\colon\sum_\pi P\tox[\id_B]{}S}
  {\underline{p\circ\cort{\pi}P}=p}
\]

{\bf Equality objects:}
\[
  \frac{p\colon\top_{A}\tox[\Delta]{}S}
  {\rho_\Delta\circ(\xi p)=p}\gop
  \frac{p\colon\Eq_B\tox[\id_B]{}S}
  {\xi(p\circ{}\rho_\Delta)=p}
\]

{\bf Stability of the operations under pullbacks}
\[
  \inv^\top_A\circ!_{h^*\top_A}=\id_{h^*\top_A}\gap
  \ex_{h^*\bot_A}\circ\inv^\bot_A=\id_{h^*\bot_A}
\]
\[
  \br{h^*\pi_{PQ},h^*\pi'_{PQ}}=(\inv^\wedge_{PQ})\I\gap
  [h^*\kappa_{PQ},h^*\kappa'_{PQ}]=(\inv^\vee_{PQ})\I
\]\[
  h^*\varepsilon_{QR}\circ \inv^\wedge_{(Q\to R)Q}=(\inv^\To_{QR})\I
\]\[
  \underline{\cind{(\cort{P}{\pi})\circ(\crt{h}P)}^{\pi',k}}=(\inv^\Sigma_{\pi\pi'hk})\I\gap
  \lambda\Big(
  h^*\varepsilon^{\Pi}_{\pi,P}\circ
  \cind{\cind{(\crt{k}{\tprod_\pi P})\circ(\crt{\pi'}{(k^*\tprod_\pi P)})}^{h,\pi}}^{\id_{A'},h}
  \Big)
  =(\inv^\Pi_{\pi\pi'hk})\I
\]
}

We now consider the initial model of the algebraic theory given by the above signature and identities, and we define a fibration $\fibr{C}CB$ by taking the set of objects of $\C$ to be $\Ob\C$ (as defined above) -- and in particular taking the objects in the fiber over $A\in\B$ to be $\Ob\fib{C}^A$ (as defined above). Next, for a morphism $f\colon A\to B$ in $\B$ and objects $P\in\Ob\fib{C}^A$ and $Q\in\Ob\fib{C}^B$, we take the set of morphisms $P\to Q$ lying over $f$ to be the set of elements of this initial model of sort $(P,Q,f)$ (i.e., $P\tox[f]{}Q$ in the above notation).

For composition of morphisms in $\C$, we take the operation $\circ$. It follows from the ``Category'' identities that this makes $\C$ into a category, and by the definition of the operation $\circ$, we have that this makes $\fib{C}$ a prefibration (i.e., a functor).

The remaining operations and equations were chosen in precisely such a way so as to ensure that $\fib{C}$ is a $h^=$-fibration.

\subsubsection{}\prop
The $h^=$-fibration $\fibr{C}CB$ defined above is free over $\B$.
\pf
Let $\fibr{C'}{C'}{B}$ be a second $h^=$-fibration over $\B$. We must show that there is up to isomorphism a unique functor $\fib{C}\to\fib{C'}$ of $h^=$-fibrations over $\B$ (for general $\fib{C'}$, this will require the axiom of choice).

By ``specifying all the operations'' in $\fib{C'}$ -- i.e., choosing cartesian lifts (a cleavage), fiberwise products, coproducts, and exponentials, indexed products and sums along product projections, and equality objects -- we obtain a $\sigma_\B$-structure with underlying ($\Ob\sigma_\B=\Ob\B$-indexed-)set $\Ob\C'$.

By the initiality of the $\sigma_\B$-structure $\Ob\C$, this gives us a function $\Phi\colon\Ob\C\to\Ob\C'$ over $\B$ and preserving the specified operations.

We now consider a structure over $\sigma_\B'$, where the set of elements of sort $(P,Q,f)$ is the set $\Hom_f^{\fib{C'}}(\Phi P,\Phi Q)$ of morphisms $\Phi P\to\Phi Q$ in $\fib\C'$ lying over $f$. The operations are given in the obvious way. For instance, given $f\colon A\to B$ in $\B$ and $P\in\fib{C}^A$, the element $\crt{f}P$ of sort $(f^*P,P,f)$ should be a morphism $\Phi(f^*P)\to\Phi P$ lying over $f$. But by the definition of $\Phi$, we know that $\Phi(f^*P)$ is equal to the specified pullback $f^*(\Phi P)$, and hence we can take the morphism in question to be the specified cartesian lift $\crt{f}{(\Phi P)}\colon f^*(\Phi P)\to \Phi P$.

It follows from the fact that $\fib{C'}$ is an $h^=$-fibration that this structure satisfies all the identities listed above. Hence, by the initiality (with respect to the structures satisfying those identities) of the $\sigma_\B'$-structure given by the morphisms of $\C$, we obtain a function $\Hom_{f}^{\fib C}(P,Q)\to\Hom_f^{\fib C'}(\Phi P,\Phi Q)$ preserving all of the operations. In particular, preservation of $\circ$ and $\id$ implies that this is a functor, and it is then by definition a morphism of prefibrations.

The preservation of the remaining operations prove that this morphism of prefibrations preserve the \emph{specified} $h^=$-fibration structure (the specified fiberwise product and coproducts, etc.), from which it follows that \emph{all} instances of such operations (arbitrary fiberwise products and coproducts, etc.) are preserved.

This gives us the existence in the freeness statement.

Now suppose that we have two morphisms of $h^=$-fibrations $\Phi,\Phi'\colon\fib C\to\fib C'$ over $\B$. We would like to show that they are isomorphic over $\B$.

We define isomorphisms $\alpha_P\colon\Phi P\toi \Phi'P$ by structural recursion on $P$. In each case, the isomorphism is determined by the universal properties of $\Phi P$ and $\Phi' P$ using the isomorphisms already established. For example, suppose $P$ is of the form $Q\wedge R$, and we have already isomorphisms $\Phi Q\toi \Phi'Q$ and $\Phi R\toi \Phi'R$. We know that $\Phi P$ is a product of $\Phi Q$ and $\Phi R$ and that $\Phi'P$ is a product of $\Phi'Q$ and $\Phi'R$, and so we have an induced isomorphism $\Phi P\toi\Phi'P$.

It remains to see that $\alpha$ is a natural isomorphism over $\B$. Note first of all that the isomorphisms all lie over identity morphisms in $\B$, as required.

We check the naturality condition for morphisms $p$ in $\C$ by structural induction on $p$ (that is, recalling that the morphisms of $\C$ are equivalence classes of closed $\sigma_\B'$-terms, we do an induction on the $\sigma_\B'$-terms themselves). The base cases of this induction correspond to the 0-ary operations in $\sigma_\B'$. In each case, the naturality condition follows from the definition of $\alpha$. For instance, the isomorphism $\alpha_{P\wedge Q}$ is defined precisely in such a way that the naturality squares for $\pi_{PQ}$ and $\pi'_{PQ}$ commute. (The base cases also include the various ``$\mathrm{inv}$'' morphisms under ``Stability of the operations under pullbacks''; in this case, naturality follows from the fact that the naturality condition for a morphism implies the naturality condition for its inverse.)

The various induction steps all amount to showing that the naturality condition for a morphism induced by a certain universal property follows from the naturality conditions for the morphisms it is induced from. This can be checked straightforwardly in each case.
\qed

\subsubsection{}
We now explain a possible second, more conceptual (though not necessarily simpler) approach to proving the existence of a free $h^=$-fibration over an f.p.\ category $\B$, using the ``2-dimensional monad theory'' of \cite{bkp-2-monad-thy}.

We first observe that the $h^=$-fibrations over $\B$ in which all of the operations (cartesian lifts fiberwise products and coproducts, etc.) have been specified -- let us call these ``$h^=$-cloven'' $h^=$-fibrations -- form a category, in which the morphisms are morphisms of $h^=$-fibrations over $\B$ that preserve all of the specified operations; and in fact, this category has a natural 2-categorical structure, in which we take as 2-cells all natural transformations over $\B$. It also comes with a forgetful 2-functor to the 2-category $\Cat/\B$ of categories over $\B$.

We now would like to say that this forgetful functor is \emph{monadic}; i.e., that the $h^=$-cloven $h^=$-fibrations are precisely the $T$-algebras for a (in fact finitary) 2-monad $T$ on $\Cat/\B$. One can imagine constructing this monad by hand; for instance, the total category of $T(\fibr{C}CB)$ should have as objects terms over the multi-sorted signature $\sigma_\B$ introduced in Construction~\ref{constr:free-hfib}, with $\Ob\C$ as the set of (sorted) variables, and as morphisms (equivalence classes of) terms over the signature $\sigma_\B'$ with $\Ar\C$ as the set of sorted variables. Alternatively, one might hope to appeal to a $\Cat$-enriched version of the Beck monadicity theorem.

In any case, if this monadicity is established -- and assuming, as one would expect, that the morphisms of $T$-algebras (and 2-cells between them) in the sense of \cite[p.~3]{bkp-2-monad-thy} are precisely the morphisms of $h^=$-fibrations (and natural transformations) -- then the existence of the free $h^=$-fibration over $\B$ follows easily from the results of \emph{op. cit.}

Namely, we have (as usual) that the free $h^=$-cloven $h^=$-fibration $\fib{C}$
produced in Construction~\ref{constr:free-hfib} is the free $T$-algebra on the empty category $\emptyset\to\B$ over $\B$ (i.e., the value at $\emptyset\to\B$ of the left adjoint to the forgetful functor from $h^=$-cloven $h^=$-fibrations). The desired universal property is then established as follows. According to \cite[Corollary~5.6]{bkp-2-monad-thy}, $\fib{C}$ (and more generally any free $T$-algebra) is ``flexible'', which implies by \cite[Theorem~4.7]{bkp-2-monad-thy} that every morphism out of $\fib{C}$ is isomorphic to a strict one. Since, for each $h^=$-fibration $\fib{C'}$ over $\B$ (once we make $\fib{C'}$ into a $T$-algebra by arbitrarily specifying operations on it), there is a unique strict morphism $\fib{C}\to\fib{C'}$, it follows that any two morphisms $\fib{C}\to\fib{C'}$ are isomorphic, as desired.

\subsection{Deductions}\label{subsec:deductions}
We will now formally introduce the notion of ``intuitionistic validity'' of first-order formulas; i.e., we introduce the (or rather a) deductive system for intuitionistic first-order logic.

The axiomatization we present is the one due to Lambek (for the propositional part) and Lawvere, in which the rules of inference corresponding to each logical connective correspond to the universal property of the corresponding categorical operation. As a result, we do not even need to state most of the rules, but rather just refer to the corresponding $h^=$-fibration operations as listed in Construction~\ref{constr:free-hfib}.

As we mentioned in the \hyperref[sec:appendix]{introduction to the appendix}, this makes it ``almost'' obvious that the morphisms in the free $h^=$-fibration are given by equivalence classes of deductions. But, as we also mentioned there, there are exceptions: namely, the rules for quantification, equality and substitution take on a different form.

For quantification, the difference is that in the syntax, we only quantify over a single variable at a time, whereas in the fibration, we take indexed sums and products along \emph{arbitrary} product projections (corresponding to quantifying over several variables at once).

For equality, there is a similar circumstance, that we consider cocartesian lifts of \emph{arbitrary} diagonal morphisms (corresponding to equating several pairs of terms at once). In addition, these cocartesian lifts correspond to equalities between variables, whereas in the syntax, we have equalities between arbitrary terms.

Finally, the substitution rule below corresponds to the pullback operation $p\mapsto f^*p$, which is not primitive among the operations listed in Construction~\ref{constr:free-hfib}, but rather is derived from the cartesian lifts and their universal property.

Our definition of ``deduction'' will be given in terms of an initial structure for a certain algebraic signature, just as were the morphisms for the free $h^=$-fibration in Construction~\ref{constr:free-hfib}. Hence, the deductions are given by closed terms over this signature (``proof terms'' as they are often called in the context of dependent type theory) rather than the traditional ``trees''; but of course, the tree corresponding to a certain deduction is recovered as the ``syntax tree'' of the corresponding term.

After introducing the notion of deduction, we prove the fundamental \emph{soundness} property: that provable formulas are satisfied by all interpretations into a fibration; and for that purpose, we need another fundamental fact about fibrations which we have not mentioned explicitly, and which we state in Lemma~\ref{lem:interp-subs-pullbacks}, namely the one relating pullbacks in a fibration to the syntactic operation of substitution.

\subsubsection{}\defn
Let $\sigma$ be a signature. We now define the set of \emph{deductions} of first-order formulas over $\sigma$. Associated to each premise is a \emph{context} $\vec{x}$, which is (as usual) a tuple of sorted variables, as well as a \emph{premise} and \emph{conclusion}, which are $\sigma$-formulas whose free variables are among those in $\vec{x}$. To indicate that $p$ is a deduction with the context, premise, and conclusion $\vec{x}$, $\phi$, and $\psi$, we write $p\colon \phi\le_{\vec{x}}\psi$.

Now, the set of deductions is defined inductively, and is in fact the initial structure for a certain multi-sorted algebraic signature, whose sorts are the triples $\phi\le_{\vec{x}}\psi$ as above, which we call \emph{sequents}. We must now describe the function symbols -- which we will call ``rules of inference'' -- and their arities and output sorts. For the most part, these are given as a subset of those listed in the first table in Construction~\ref{constr:free-hfib} (on page \pageref{fig:first-hfib-table}) -- namely, those under ``Finite products and coproducts'' and ``Exponentials'' and the first two under ``Category and fibration structure'' -- with certain modifications. The modifications are: (i) we now allow $P,Q,R$ to be arbitrary \emph{formulas} (so that for each figure in the table, we get one rule of inference for each choice of such formulas), (ii) we ignore the subscripts $_{f}$, $_g$, $_{gf}$, and $_A$ in the figures for $q\circ p$, $!_P$ and $\ex_P$ and (iii) we replace each arrow $\to$ with $\le_{\vec{x}}$, where $\vec{x}$ is an arbitrary context containing all the free variables of the formulas $P,Q,R$ involved (so that, again, we get one rule of inference for each possible value of $\vec{x}$). The arities and output sorts of the function symbols are to be understood in the same way as in Construction~\ref{constr:free-hfib} (e.g., the arity of the rule of inference corresponding to the function symbol ``$\circ$'' is $\br{P\le_{\vec{x}}Q,Q\le_{\vec{x}}R}$ and its output sort is $P\le_{\vec{x}}R$).

Finally, we have the following rules of inference (corresponding to the rest of the operations in the table on page \pageref{fig:first-hfib-table}, except for the ``stability'' ones), where again the arities and output sorts are to be understood in the same way as above. Note that unlike in the table on page \pageref{fig:first-hfib-table}, we are not introducing any notation for these function symbols (i.e., rules of inference).

Again, each figure below presents a set of inference rules, one for each value of the relevant parameters. Here, $\vec{x}$, $\vec{y}$ are arbitrary sequences of variables, and $y,z$ are arbitrary variables, with $z$ and $y$ of the same sort in the final rule; $P,Q,R$ are arbitrary formulas, with certain restrictions on the free variables; namely, they should be among those in $\vec{y}$ in the ``Substitution rule'', among those in $\vec{x}z$ in the ``Quantifiers'' rules, and among those in $\vec{x}yz$ in the ``Equality rule'', and moreover $z$ should not be free in $P$ and $Q$, respectively, in the first and third quantifier rules; and $\vec{t}$ is an arbitrary tuple of terms with free variables among those in $\vec{x}$ and with $\tp\vec{t}=\tp\vec{y}$.

We denote by $X[\vec{a}:=\vec{b}]$ the result of (simultaneously) substituting the terms $\vec{b}$ for the variables $\vec{a}$ in a term or formula $X$ (where we may assume that the none of the variables in $\vec{a}$ are bound in $X$ -- see Definition~\ref{defn:syntax}).

{
\newcommand\gap{\hspace{20pt}}
\newcommand\gop{\hspace{40pt}}
\newcommand\blk{\vphantom{P}}

{\bf Substitution:}
\[
  \frac{P\le_{\vec{y}}Q}{P[\vec{y}:=\vec{t}]\le_{\vec{x}}Q[\vec{y}:=\vec{t}]}
\]

{\bf Quantifiers:}
\[
  \frac{P\le_{\vec{x}z}Q}{P\le_{\vec{x}}\forall z Q}\gap
  \frac{\blk}{\forall z P\le_{\vec{x}z} P}\gop
  \frac{P\le_{\vec{x}z}Q}{\exists z P\le_{\vec{x}}Q}\gap
  \frac{\blk}{P\le_{\vec{x}z}\exists z P}
\]

{\bf Equality:}
\[
  \frac{\blk}{\top\le_{\vec{x}z}z=z}\gop
  \frac{\top\le_{\vec{x}y}P[z:=y]}{y=z\le_{\vec{x}yz}P}
\]
}

We say that a sequent $\phi\le_{\vec{x}}\psi$ is \emph{provable} if it is the sort of some deduction, and that a formula $\phi$ is \emph{provable} if $\top\le_{\vec{x}}\phi$ is provable for some (and hence, by the substitution rule, every) context $\vec{x}$ containing the free variables of $\phi$.

Generalizing Definition~\ref{defn:interpretation-over}, given an $h^=$-fibration $\fibr{C}CB$ and an interpretation $(\widehat M,M)\colon\sigma\to\fib{C}$, we say that $M$ \emph{satisfies} the sequent $\phi\le_{\vec{x}}\psi$ if there exists a morphism $\widehat M_{\vec{x}}(\phi)\to\widehat M_{\vec{x}}(\psi)$ over $M(\vec{x})$.

\subsubsection{}\lem\label{lem:interp-subs-pullbacks}
Let $\B$ be an f.p.\ category, and let $M\colon\sigma\to\B$ be an interpretation. Then for any tuples-of-terms-in-context $(\vec{t},\vec{y})$ and $(\vec{s},\vec{x})$ with $\tp(\vec{s})=\vec{y}$, we have $M_{\vec{y}}\vec{t}\cdot M_{\vec{x}}\vec{s}=M_{\vec{x}}(\vec{t}[\vec{y}:=\vec{s}])$.

Next, given an $h^=$-fibration $\fibr{C}CB$ over $\B$ and an interpretation $\widehat M$ in $\fib{C}$ over $M$, then for any tuple-of-terms-in-context $(\vec{t},\vec{x})$ and any formula-in-context $(\phi,\vec{y})$ with $\tp(\vec{t})=\tp(\vec{y})$, we have $(M_{\vec x}\vec{t})^*(\widehat M_{\vec y}\phi)\cong\widehat M_{\vec{x}}(\phi[\vec{x}:=\vec{t}])$.
\pf
This is well-known (and as we mentioned, is fundamental to the role of fibrations in logic). Both claims are proven by induction. The second claim uses the first claim for the base case of equality, and the remaining induction steps all make use of the various stability conditions in an $h^=$-fibration. We leave the details to the reader.
\qed

\subsubsection{}\prop\label{prop:soundness}
Given any interpretation $(\widehat{M},M)\colon\sigma\to\fib{C}$ of a signature $\sigma$, $M$ satisfies every provable sequent $\phi\le_{\vec{x}}\psi$ over $\sigma$.
\pf
We need to show that for every deduction $p\colon\phi\le_{\vec{x}}\psi$, there exists a morphism $\widehat M_{\vec{x}}(\phi)\to \widehat M_{\vec{x}}(\psi)$ over $M(\vec{x})$ and of course, we prove this by induction on $p$.

But in each induction step, it follows immediately from the definition of $\widehat M$ that such a morphism exists, where for the ``Substitution'', ``Quantifier'', and ``Equality'' rules, we must use Lemma~\ref{lem:interp-subs-pullbacks} concerning pullbacks and substitution (and the special case of it that says that if $z$ is not free in $\phi$, then $M_{\vec{x}z}(\phi)=\pi^*M_{\vec{x}}(\phi)$, with $\pi=\br{\pi^M_1,\ldots,\pi^M_{\len{\vec{x}}}}$).

For the ``Equality'' rule, we also need to use the stability of equality morphisms (i.e., cocartesian lifts of generalized diagonal morphisms -- see \cite[\S\S1.6-1.7]{sentai}), which holds in $h^=$-fibrations by Proposition~\ref{prop:hfib-is-wedgeq}.
\qed

\subsection{The free $h^=$-fibration over a free f.p.\ category}\label{subsec:free-over-free}
We now specialize to the case of a free $h^=$-fibration over the free f.p.\ category on some signature $\sigma$, and draw certain ``syntactic'' conclusions. Namely, in this case, as we have said, the free $h^=$-fibration can be given a direct, syntactic construction, in which the objects of the total category are formulas and the morphisms are equivalence classes of deductions, and we would like to recover this description from our construction.

Here, we will recover only part of this description: we will show that the objects of the total category ``are'' the formulas (more precisely, that each object is isomorphic to the interpretation of some formula), and that a morphism exists between two of these if and only if the corresponding implication is provable (as we explain in the proof, this amounts to a completeness theorem for the fibrational semantics).

We note that the ``posetal'' version of the syntactic construction of the free $h^=$-fibration is much simpler, and in fact, we will make use of it in the proof of said completeness theorem.

We would like to mention a second way to recover the morphisms-as-deductions description of the free $h^=$-fibration. Namely, rather than adapting the construction of the fibration to conform to the syntax, we can observe that the ``syntax of first-order logic'' is not absolute; one could also just use a different syntax and deductive system so as to be adapted to the above construction of the free $h^=$-fibration. This would involve slightly modifying the implementation of equality and quantification in the syntax, as well as introducing ``explicit substitutions'', i.e., substitution as a primitive operation.

To begin with, we recall the well-known construction, due to Lawvere, of a free f.p.\ category associated to an algebraic signature $\sigma$, cf. \cite[p.~473]{makkai-lauchli2}.

\subsubsection{}\prop\label{prop:free-fp-exists}
Given an algebraic signature $\sigma$, there exists a free f.p.\ category $\B$ over $\sigma$ via an interpretation $i\colon\sigma\to\B$. Moreover any such $\B$ and $i$ have the following property: every object in $\B$ is of the form $i(\vec{A})$ for a unique $\vec{A}\in(\Ob\sigma)^{<\omega}$, and every morphism $i(\vec{A})\to i(\vec{B})$ is of the form $i_{\vec{x}}(\vec{t})$ for a tuple-of-terms-in-context $(\vec{t},\vec{x})$, unique up to substituting $\vec{x}$ for a different sequence $\vec{y}$ of variables with $\tp(\vec{y})=\vec{A}$ (i.e., replacing $(\vec{t},\vec{x})$ with $(\vec{t}[\vec{x}:=\vec{y}],\vec{y})$.
\pf
The stated property together with the first part of Lemma~\ref{lem:interp-subs-pullbacks} determine $\B$ and $i\colon \sigma\to\B$ up to isomorphism, and it is straightforward to verify that this indeed produces an f.p.\ category and has the required universal property.
\qed

\subsubsection{}\defn
Proposition~\ref{prop:free-fp-exists} gives an explicit description of a free f.p.\ category on an algebraic signature $\sigma$: the set of objects is $(\Ob\sigma)^{<\omega}$, the morphisms are equivalence classes of tuples-of-terms-in-context up to renaming of variables, and composition is given by substitution as in Lemma~\ref{lem:interp-subs-pullbacks}.

We denote the resulting category by $\TM_\sigma$ and the associated $\sigma$-interpretation (which is the identity on $(\Ob\sigma)^{<\omega}$ and takes each tuple-of-terms-in-context to its equivalence class) by $i\colon\sigma\to\TM_\sigma$.

\subsubsection{}\prop\label{prop:free-fib-interp-ess-surj}
Given a free $h^=$-fibration $\fibr{C}{C}{Tm_\sigma}$ over $\TM_\sigma$ and an interpretation $\hat\imath$ in $\fib{C}$ over $i$, each object in $\C$ is isomorphic to $\hat\imath_{\vec{x}}\phi$ for some $(\vec{x},\phi)$.
\pf
Let $\mathcal{A}\subset\Ob\C$ be the set of objects which are isomorphic to some $\hat\imath_{\vec{x}}\phi$. It suffices to see that $\mathcal{A}$ is closed under all of the operations which inductively define the set $\Ob\C$.

It is obvious from the definition of $\hat\imath_{\vec{x}}\phi$ that $\mathcal{A}$ is closed under fiberwise products, coproducts, and exponentials and since every object of $\TM_\sigma$ is of the form $i(\vec{A})$, it follows that $\mathcal{A}$ includes initial and terminal objects in each fiber. By Lemma~\ref{lem:interp-subs-pullbacks} and since every morphism in $\TM_\sigma$ is of the form $i_{\vec{x}}(\vec{t})$, $\mathcal{A}$ is also closed under pullbacks.

Next, we have that $\mathcal{A}$ is closed under \emph{certain} indexed products and sums, namely the ones along projections of the form $\br{\pi_1^i,\ldots,\pi_{\len{\vec{A}}}^i}\colon i(\vec{A}B)\to i(\vec{A})$. We also know that $\mathcal{A}$ is closed under indexed product and sums along isomorphisms, as these are the same as pullbacks. Also, if $\mathcal{A}$ is closed under indexed products and sums of any two composable morphisms, then the same is true of the composite (this is because $\prod_{gf}P\cong\prod_g\prod_fP$ and similarly with $\sum_{gf}$, since a composite of (ana-)adjoints is an (ana-)adjoint of the composite).

But now it is easy to see that every product projection in $\TM_\sigma$ is a composite of the above particular product projections and isomorphisms.

Finally, we need to see that $\mathcal{A}$ contains an equality object $\Eq_A$ for every diagonal morphisms $\Delta\colon A\to A\times A$. Firstly, it suffices to do this for \emph{some} diagonal morphism for each object $A$ (since, for example, $\mathcal{A}$ is closed under isomorphisms, as every isomorphism is cartesian). Next, for a product $A\times B$, we have that $\Eq_{A\times B}$ is a fiberwise product of pullbacks of $\Eq_A$ and $\Eq_B$ (see \cite[p.~10]{lawvereequality} or \cite[Theorem~7.2]{sentai}), so the claim for $A\times B$ follows from that of $A$ and $B$. Hence, by induction, it suffices to show the claim for each $i(\br{A})$, for which we have $\Eq_{\br{A}}\cong \hat\imath_{\br{x,y}}(x=y)$ (with $\tp(x)=\tp(y)=A$), and for the terminal object $\tm=i(\emptyset)$, for which $\Eq_\tm\cong\top_{\fib{C}^{\tm\times\tm}}$.
\qed

\subsubsection{}\thm\label{thm:completeness}
Let $\fibr{C}C{Tm_\sigma}$ be a free $h^=$-fibration over $\TM_\sigma$, and fix an interpretation $\hat\imath$ in $\fib{C}$ over $i$. Then a sequent $\phi\le_{\vec{x}}\psi$ is satisfied by $i$ if and only if it is provable.
\pf
The $(\Longleftarrow)$ implication follows from soundness (Proposition~\ref{prop:soundness}). We now consider the other implication.

Fix a sequent $\phi\le_{\vec{x}}\psi$, and suppose it is satisfied by $i$, i.e., there exists a morphism $\hat\imath_{\vec{x}}(\phi)\to\hat\imath_{\vec{x}}(\psi)$ over $i(\vec{x})$. Note that it is then satisfied by \emph{any} interpretation in any $h^=$-fibration $\fibr{C'}{C'}{B'}$. Indeed, given an interpretation $(\widehat M,M)\colon\sigma\to\fib{C'}$, we have, by the freeness of $\fib{C}$, a morphism $(\Xi,\xi)\colon\fib{C}\to\fib{C'}$ with $\xi\circ i=M$. Then $\Xi\circ\hat\imath$ and $\widehat M$ will be isomorphic interpretations over $M$, and the morphism $\hat\imath_{\vec{x}}(\phi)\to\hat\imath_{\vec{x}}(\psi)$ then gives a morphism
\[
  \widehat M_{\vec{x}}(\phi)
  \cong
  (\Xi\circ\hat\imath)_{\vec{x}}(\phi)\to
  (\Xi\circ\hat\imath)_{\vec{x}}(\psi)
  \cong
  \widehat M_{\vec{x}}(\psi)
\]
over $M(\vec{x})$.

Hence, it suffices to show the \emph{completeness of the fibrational semantics} -- i.e., that if a sequent is satisfied by \emph{every} interpretation in an $h^=$-fibration, then it is provable.

A heavy-handed way to do this is to appeal, for example, to Kripke's completeness theorem for intuitionistic logic. Kripke models are precisely interpretations in the $h^=$-fibrations $\fib{P(\Set^{\mathnormal P})}$ with $P$ an arbitrary poset. Kripke's completeness theorem then says that that if a sequent is satisfied by every interpretation in such an $h^=$-fibration, then it is provable.

Besides being overkill, this proof has the following disadvantage: we might like to \emph{deduce} Kripke's completeness theorem from its categorical formulation in \cite{makkai-lauchli1,makkai-lauchli2}, and this proof would make such a deduction circular. Hence, we now give a more direct proof.

Namely, we will produce a \emph{single} fibration $\fibr{C'}{C'}B$ over $\B$ such that any sequent satisfied by $i$ (with respect to $\fib{C'}$) is provable.

$\fib{C'}$ will be a posetal fibration (i.e., the fibers are posets), where each fiber is a ``Lindenbaum-Tarski (Heyting) algebra''. Namely, to define the fiber over $\vec{A}$, consider the preorder whose objects are formulas-in-context $(\phi,\vec{x})$ with $\tp(\vec{x})=\vec{A}$, with the ordering given by provability (i.e., $(\phi,\vec{x})\le(\psi,\vec{y})$ if $\phi\le_{\vec{x}}\psi[\vec{y}:=\vec{x}]$ is provable), and take its ``posetal collapse (or reflection)'' (i.e., identify objects that are equal in the ordering).

Organizing these fibers into a fibration amounts to defining a functor $\TM_\sigma^\op\to\cat{Poset}$ to the category of posets with the given action on objects, and we do this by sending $[\vec{t},\vec{x}]$ (square brackets denotes the equivalence class) to the function $(\fib{C'})^{\tp\vec{t}}\to(\fib{C'})^{\tp\vec{x}}$ taking $[\phi,\vec{y}]$ to $[\phi[\vec{y}:=\vec{t}],\vec{x}]$. That this is well-defined and a morphism of posets follows from the substitution rule, and it follows from the fact that $\phi[\vec{y}:=\vec{t}][\vec{x}:=\vec{s}]=\phi[\vec{y}:=(\vec{t}[\vec{x}:=\vec{s}])]$ that this defines a functor.

That the resulting fibration is fiberwise bicartesian-closed (i.e., has Heyting algebra fibers) follows immediately from the rules of inference for the propositional connectives. The rules of inference for equality and quantifiers give the existence of equality objects and indexed products and sums over \emph{certain} (generalized) diagonal and projection morphisms, but as in the proof of Proposition~\ref{prop:free-fib-interp-ess-surj}, these morphisms generate all of the diagonal morphisms and projections. We leave the details to the reader.

The resulting fibration is in fact the posetal reflection of the free fibration $\fib{C}$, and is the free posetal $h^=$-fibration over $\fib{C}$ (see \cite[p.~349]{makkai-lauchli1}), though we will need neither of these facts.

In fact, we are done, as it is immediate from the definition that, under the obvious interpretation $\hat\imath$ in $\fib{C'}$ over $i$, that a sequent is provable if and only if it is satisfied by $i$.
\qed

\subsubsection{}\cor
For any signature $\sigma$, there exists a free $h^=$-fibration $\fibr{Pf}{Pf}{Tm_\sigma}$ and an interpretation $\hat\imath$ in $\fib{C}$ over $i$ such that
\begin{itemize}[topsep=0pt]
\item The objects in $\fib{Pf}^{\vec{A}}$ are equivalence classes $[\phi,\vec{x}]$ of formulas-in-context with $\tp\vec{x}=\vec{A}$ (the equivalence relation being renaming of variables -- i.e., $(\phi,\vec{x})$ and $(\psi,\vec{y})$ are equivalent iff $\psi=\phi[\vec{x}:=\vec{y}]$).
\item We have $\hat\imath_{\vec{x}}\phi=[\phi,\vec{x}]$ for all $(\phi,\vec{x})$.
\item There is a morphism $[\phi,\vec{x}]\to[\psi,\vec{x}]$ in $\fib{Pf}^{\tp\vec{x}}$ if and only if $\phi\le_{\vec{x}}\psi$ is provable.
\end{itemize}
\pf
Given an arbitrary free $h^=$-fibration over $\TM_\sigma$, Proposition~\ref{prop:free-fib-interp-ess-surj} shows that there is an equivalent (and hence still free) $h^=$-fibration satisfying the first two properties.

The third property then follows from Theorem~\ref{thm:completeness}.
\qed

\bibliographystyle{alpha}
\bibliography{fohl}

\end{document}